\documentclass[a4paper]{article}

\usepackage{amsfonts,amssymb,amsthm}

\catcode`\@=11\relax
\newwrite\@unused
\def\typeout#1{{\let\protect\string\immediate\write\@unused{#1}}}
\def\@nnil{\@nil}
\def\@empty{}
\def\@psdonoop#1\@@#2#3{}
\def\@psdo#1:=#2\do#3{\edef\@psdotmp{#2}\ifx\@psdotmp\@empty \else
    \expandafter\@psdoloop#2,\@nil,\@nil\@@#1{#3}\fi}
\def\@psdoloop#1,#2,#3\@@#4#5{\def#4{#1}\ifx #4\@nnil \else
       #5\def#4{#2}\ifx #4\@nnil \else#5\@ipsdoloop #3\@@#4{#5}\fi\fi}
\def\@ipsdoloop#1,#2\@@#3#4{\def#3{#1}\ifx #3\@nnil 
       \let\@nextwhile=\@psdonoop \else
      #4\relax\let\@nextwhile=\@ipsdoloop\fi\@nextwhile#2\@@#3{#4}}
\def\@tpsdo#1:=#2\do#3{\xdef\@psdotmp{#2}\ifx\@psdotmp\@empty \else
    \@tpsdoloop#2\@nil\@nil\@@#1{#3}\fi}
\def\@tpsdoloop#1#2\@@#3#4{\def#3{#1}\ifx #3\@nnil 
       \let\@nextwhile=\@psdonoop \else
      #4\relax\let\@nextwhile=\@tpsdoloop\fi\@nextwhile#2\@@#3{#4}}
\def\psdraft{
        \def\@psdraft{0}
}
\def\psfull{
        \def\@psdraft{100}
}
\psfull
\newif\if@prologfile
\newif\if@postlogfile
\newif\if@noisy
\def\pssilent{
        \@noisyfalse
}
\def\psnoisy{
        \@noisytrue
}
\psnoisy
\newif\if@bbllx
\newif\if@bblly
\newif\if@bburx
\newif\if@bbury
\newif\if@height
\newif\if@width
\newif\if@rheight
\newif\if@rwidth
\newif\if@clip
\newif\if@verbose
\def\@p@@sclip#1{\@cliptrue}
\def\@p@@sfile#1{
                   \def\@p@sfile{#1}
}
\def\@p@@sfigure#1{\def\@p@sfile{#1}}
\def\@p@@sbbllx#1{
                \@bbllxtrue
                \dimen100=#1
                \edef\@p@sbbllx{\number\dimen100}
}
\def\@p@@sbblly#1{
                \@bbllytrue
                \dimen100=#1
                \edef\@p@sbblly{\number\dimen100}
}
\def\@p@@sbburx#1{
                \@bburxtrue
                \dimen100=#1
                \edef\@p@sbburx{\number\dimen100}
}
\def\@p@@sbbury#1{
                \@bburytrue
                \dimen100=#1
                \edef\@p@sbbury{\number\dimen100}
}
\def\@p@@sheight#1{
                \@heighttrue
                \dimen100=#1
                \edef\@p@sheight{\number\dimen100}
}
\def\@p@@swidth#1{
                \@widthtrue
                \dimen100=#1
                \edef\@p@swidth{\number\dimen100}
}
\def\@p@@srheight#1{
                \@rheighttrue
                \dimen100=#1
                \edef\@p@srheight{\number\dimen100}
}
\def\@p@@srwidth#1{
                \@rwidthtrue
                \dimen100=#1
                \edef\@p@srwidth{\number\dimen100}
}
\def\@p@@ssilent#1{ 
                \@verbosefalse
}
\def\@p@@sprolog#1{\@prologfiletrue\def\@prologfileval{#1}}
\def\@p@@spostlog#1{\@postlogfiletrue\def\@postlogfileval{#1}}
\def\@cs@name#1{\csname #1\endcsname}
\def\@setparms#1=#2,{\@cs@name{@p@@s#1}{#2}}
\def\ps@init@parms{
                \@bbllxfalse \@bbllyfalse
                \@bburxfalse \@bburyfalse
                \@heightfalse \@widthfalse
                \@rheightfalse \@rwidthfalse
                \def\@p@sbbllx{}\def\@p@sbblly{}
                \def\@p@sbburx{}\def\@p@sbbury{}
                \def\@p@sheight{}\def\@p@swidth{}
                \def\@p@srheight{}\def\@p@srwidth{}
                \def\@p@sfile{}
                \def\@p@scost{10}
                \def\@sc{}
                \@prologfilefalse
                \@postlogfilefalse
                \@clipfalse
                \if@noisy
                        \@verbosetrue
                \else
                        \@verbosefalse
                \fi
}
\def\parse@ps@parms#1{
                \@psdo\@psfiga:=#1\do
                   {\expandafter\@setparms\@psfiga,}}
\newif\ifno@bb
\newif\ifnot@eof
\newread\ps@stream
\def\bb@missing{
        \if@verbose{
                \typeout{psfig: searching \@p@sfile \space  for bounding box}
        }\fi
        \openin\ps@stream=\@p@sfile
        \no@bbtrue
        \not@eoftrue
        \catcode`\%=12
        \loop
                \read\ps@stream to \line@in
                \global\toks200=\expandafter{\line@in}
                \ifeof\ps@stream \not@eoffalse \fi
                \@bbtest{\toks200}
                \if@bbmatch\not@eoffalse\expandafter\bb@cull\the\toks200\fi
        \ifnot@eof \repeat
        \catcode`\%=14
}       
\catcode`\%=12
\newif\if@bbmatch
\def\@bbtest#1{\expandafter\@a@\the#1
\long\def\@a@#1
\long\def\bb@cull#1 #2 #3 #4 #5 {
        \dimen100=#2 bp\edef\@p@sbbllx{\number\dimen100}
        \dimen100=#3 bp\edef\@p@sbblly{\number\dimen100}
        \dimen100=#4 bp\edef\@p@sbburx{\number\dimen100}
        \dimen100=#5 bp\edef\@p@sbbury{\number\dimen100}
        \no@bbfalse
}
\catcode`\%=14
\def\compute@bb{
                \no@bbfalse
                \if@bbllx \else \no@bbtrue \fi
                \if@bblly \else \no@bbtrue \fi
                \if@bburx \else \no@bbtrue \fi
                \if@bbury \else \no@bbtrue \fi
                \ifno@bb \bb@missing \fi
                \ifno@bb \typeout{FATAL ERROR: no bb supplied or found}
                        \no-bb-error
                \fi
                \count203=\@p@sbburx
                \count204=\@p@sbbury
                \advance\count203 by -\@p@sbbllx
                \advance\count204 by -\@p@sbblly
                \edef\@bbw{\number\count203}
                \edef\@bbh{\number\count204}
}
\def\in@hundreds#1#2#3{\count240=#2 \count241=#3
                     \count100=\count240        
                     \divide\count100 by \count241
                     \count101=\count100
                     \multiply\count101 by \count241
                     \advance\count240 by -\count101
                     \multiply\count240 by 10
                     \count101=\count240        
                     \divide\count101 by \count241
                     \count102=\count101
                     \multiply\count102 by \count241
                     \advance\count240 by -\count102
                     \multiply\count240 by 10
                     \count102=\count240        
                     \divide\count102 by \count241
                     \count200=#1\count205=0
                     \count201=\count200
                        \multiply\count201 by \count100
                        \advance\count205 by \count201
                     \count201=\count200
                        \divide\count201 by 10
                        \multiply\count201 by \count101
                        \advance\count205 by \count201
                     \count201=\count200
                        \divide\count201 by 100
                        \multiply\count201 by \count102
                        \advance\count205 by \count201
                     \edef\@result{\number\count205}
}
\def\compute@wfromh{
                \in@hundreds{\@p@sheight}{\@bbw}{\@bbh}
                \edef\@p@swidth{\@result}
}
\def\compute@hfromw{
                \in@hundreds{\@p@swidth}{\@bbh}{\@bbw}
                \edef\@p@sheight{\@result}
}
\def\compute@handw{
                \if@height 
                        \if@width
                        \else
                                \compute@wfromh
                        \fi
                \else 
                        \if@width
                                \compute@hfromw
                        \else
                                \edef\@p@sheight{\@bbh}
                                \edef\@p@swidth{\@bbw}
                        \fi
                \fi
}
\def\compute@resv{
                \if@rheight \else \edef\@p@srheight{\@p@sheight} \fi
                \if@rwidth \else \edef\@p@srwidth{\@p@swidth} \fi
}
\def\compute@sizes{
        \compute@bb
        \compute@handw
        \compute@resv
}
\def\psfig#1{\vbox {
        \ps@init@parms
        \parse@ps@parms{#1}
        \compute@sizes
        \ifnum\@p@scost<\@psdraft{
                \if@verbose{
                        \typeout{psfig: including \@p@sfile \space }
                }\fi
                \special{ps::[begin]    \@p@swidth \space \@p@sheight \space
                                \@p@sbbllx \space \@p@sbblly \space
                                \@p@sbburx \space \@p@sbbury \space
                                startTexFig \space }
                \if@clip{
                        \if@verbose{
                                \typeout{(clip)}
                        }\fi
                        \special{ps:: doclip \space }
                }\fi
                \if@prologfile
                    \special{ps: plotfile \@prologfileval \space } \fi
                \special{ps: plotfile \@p@sfile \space }
                \if@postlogfile
                    \special{ps: plotfile \@postlogfileval \space } \fi
                \special{ps::[end] endTexFig \space }
                \vbox to \@p@srheight true sp{
                        \hbox to \@p@srwidth true sp{
                                \hss
                        }
                \vss
                }
        }\else{
                \vbox to \@p@srheight true sp{
                \vss
                        \hbox to \@p@srwidth true sp{
                                \hss
                                \if@verbose{
                                        \@p@sfile
                                }\fi
                                \hss
                        }
                \vss
                }
        }\fi
}}
\catcode`\@=12\relax

\textwidth=16cm
\textheight=23cm
\oddsidemargin = -0.7truecm
\topmargin =-1truecm

\newcommand{\dr}{\partial}
\newcommand{\C}{{\bf C}}
\newcommand{\R}{{\bf R}}
\newcommand{\N}{{\bf N}}
\newcommand{\Z}{{\bf Z}}
\newcommand{\II}{I\hspace{-0.1cm}I}
\newcommand{\III}{I\hspace{-0.1cm}I\hspace{-0.1cm}I}
\newcommand{\tr}{\mbox{tr}}
\newcommand{\cat}{\mbox{CAT}}
\newcommand{\ric}{\mbox{ric}}
\newcommand{\dev}{\mbox{dev}}
\newcommand{\sh}{\mbox{sh}}
\newcommand{\ch}{\mbox{CH}}
\newcommand{\arctg}{\mbox{arctg}}
\newcommand{\SO}{\mbox{SO}}
\newcommand{\CP}{\mbox{{\bf C}P}}
\newcommand{\can}{\mbox{can}}
\newcommand{\argth}{\mbox{argth}}
\newcommand{\hess}{\mbox{Hess}}
\newcommand{\dmt}{\tilde{\dr M}}

\newcommand{\ihm}{ideal hyperbolic manifold }

\newcommand{\deltab}{\overline{\delta}}
\newcommand{\gammab}{\overline{\gamma}}
\newcommand{\sigmab}{\overline{\sigma}}
\newcommand{\Sib}{\overline{\Sigma}}
\newcommand{\Thetab}{\overline{\Theta}}
\newcommand{\Omegab}{\overline{\Omega}}

\newcommand{\cb}{\overline{c}}
\newcommand{\eb}{\overline{e}}
\newcommand{\fb}{\overline{f}}
\newcommand{\gb}{\overline{g}}
\newcommand{\hb}{\overline{h}}
\newcommand{\nb}{\overline{n}}
\newcommand{\tb}{\overline{t}}
\newcommand{\vb}{\overline{v}}
\newcommand{\xb}{\overline{x}}
\newcommand{\Cb}{\overline{C}}
\newcommand{\Gb}{\overline{G}}
\newcommand{\Kb}{\overline{K}}
\newcommand{\Pb}{\overline{P}}
\newcommand{\Mb}{\overline{M}}
\newcommand{\Sb}{\overline{S}}
\newcommand{\ricb}{\overline{\ric}}
\newcommand{\Db}{\overline{D}}
\newcommand{\Rb}{\overline{R}}
\newcommand{\thetab}{\overline{\theta}}
\newcommand{\omegab}{\overline{\omega}}

\newcommand{\Ct}{\tilde{C}}
\newcommand{\Ft}{\tilde{F}}
\newcommand{\Mt}{\tilde{M}}
\newcommand{\St}{\tilde{S}}
\newcommand{\dMt}{\tilde{\dr M}}
\newcommand{\phit}{\tilde{\phi}}
\newcommand{\gammat}{\tilde{\gamma}}
\newcommand{\sigmat}{\tilde{\sigma}}
\newcommand{\Sigmat}{\tilde{\Sigma}}

\newcommand{\Hsh}{{\mathcal H}_{\mbox{sh}}}
\newcommand{\Hsi}{{\mathcal H}_{\mbox{si}}}
\newcommand{\Hsm}{{\mathcal H}_{\mbox{sm}}}

\newcommand{\Thetasm}{\Theta_{\mbox{sm}}}

\newcommand{\bM}{{\bf M}}

\newcommand{\Rr}{\stackrel{\circ}{R}}

\newtheorem{prop}{Proposition}[section]
\newtheorem{lemma}[prop]{Lemma}
\newtheorem{sublemma}[prop]{Sub-lemma}

\newtheorem{thm}[prop]{Theorem}
\newtheorem{cor}[prop]{Corollary}
\newtheorem{remark}[prop]{Remark}

\newtheorem{df}[prop]{Definition}
\newtheorem{pty}[prop]{Property}
\newtheorem{question}[prop]{Question}

\newcommand{\pg}[1]{\subsection{#1}}

\newenvironment{thn}[1]{\vskip 0.2cm \noindent{\bf Theorem #1.} \it}{\rm
\vspace{0.2cm}} 
\newenvironment{crn}[1]{\vskip 0.2cm \noindent{\bf Corollary #1.} \it}{\rm
\vspace{0.2cm}} 
\newenvironment{lmn}[1]{\vskip 0.2cm \noindent{\bf Lemma #1.} \it}{\rm
\vspace{0.2cm}} 
\newenvironment{qn}[1]{\vskip 0.2cm \noindent{\bf Question #1.} \it}{\rm
\vspace{0.2cm}} 
\newenvironment{sketch}{\vskip 0.2cm \noindent{\bf Brief sketch of the
    proof~. ~~~}}{$\qed$ \vspace{0.2cm}} 

\newcommand{\btm}{\begin{thm}}
\newcommand{\etm}{\end{thm}}
\newcommand{\bpt}{\begin{pty}}
\newcommand{\ept}{\end{pty}}
\newcommand{\blm}{\begin{lemma}}
\newcommand{\elm}{\end{lemma}}
\newcommand{\bsl}{\begin{sublemma}}
\newcommand{\esl}{\end{sublemma}}
\newcommand{\bcr}{\begin{cor}}
\newcommand{\ecr}{\end{cor}}
\newcommand{\bdf}{\begin{df}}
\newcommand{\edf}{\end{df}}
\newcommand{\bprop}{\begin{prop}}
\newcommand{\eprop}{\end{prop}}
\newcommand{\bas}{\begin{asser}}
\newcommand{\eas}{\end{asser}}
\newcommand{\beq}{\begin{equation}}
\newcommand{\eeq}{\end{equation}}
\newcommand{\bpv}{\begin{proof}}
\newcommand{\epv}{\end{proof}}
\newcommand{\bpvs}{\begin{sketch}}
\newcommand{\epvs}{\end{sketch}}
\newcommand{\bit}{\begin{itemize}}
\newcommand{\eit}{\end{itemize}}
\newcommand{\bpn}{\begin{pfn}}
\newcommand{\epn}{\end{pfn}}
\newcommand{\btn}{\begin{thn}}
\newcommand{\etn}{\end{thn}}
\newcommand{\bcn}{\begin{crn}}
\newcommand{\ecn}{\end{crn}}
\newcommand{\bqn}{\begin{qn}}
\newcommand{\eqn}{\end{qn}}
\newcommand{\bln}{\begin{lmn}}
\newcommand{\eln}{\end{lmn}}
\newcommand{\brk}{\begin{remark}}
\newcommand{\erk}{\end{remark}}
\newcommand{\bq}{\begin{question}}
\newcommand{\eq}{\end{question}}

\newenvironment{pfn}[1]{\vskip 0.2cm \noindent{\it Proof #1.}}{$\square$
\vspace{0.2cm}}

\newcommand{\cA}{{\mathcal A}}
\newcommand{\cC}{\mathcal{C}}
\newcommand{\cD}{\mathcal{D}}
\newcommand{\cG}{\mathcal{G}}
\newcommand{\cH}{{\mathcal H}}
\newcommand{\cM}{{\mathcal M}}
\newcommand{\cP}{{\mathcal P}}
\newcommand{\cT}{{\mathcal T}}

\newcommand{\cAb}{\overline{\mathcal A}}

\newcommand{\Met}{\mathcal{M}et}
\newcommand{\Imm}{\mathcal{I}mm}
\newcommand{\CMet}{\mathcal{CM}et}
\newcommand{\CImm}{\mathcal{CI}mm}
\newcommand{\gab}{\overline{\gamma}}
\newcommand{\hyp}{\mathbf{H}^3}
\newcommand{\dhyp}{\partial\hyp}
\newcommand{\cL}{\mathcal{L}}
\newcommand{\isom}{\mathrm{Isom}}
\newcommand{\db}{\overline{\partial}}

\newcommand{\hbu}{\stackrel{\bullet}{h}}
\newcommand{\hbbu}{\stackrel{\bullet}{\hb}}
\newcommand{\Vb}{\stackrel{\bullet}{V}}
\newcommand{\alphab}{\stackrel{\bullet}{\alpha}}
\newcommand{\thetabu}{\stackrel{\bullet}{\theta}}
\newcommand{\thetabbu}{\stackrel{\bullet}{\thetab}}

\begin{document}

\title{Hyperbolic manifolds with polyhedral boundary}

\author{Jean-Marc Schlenker\thanks{
Laboratoire Emile Picard, UMR CNRS 5580,
Universit{\'e} Paul Sabatier,
118 route de Narbonne,
31062 Toulouse Cedex 4,
France.
\texttt{schlenker@picard.ups-tlse.fr; http://picard.ups-tlse.fr/\~{
}schlenker}. }}

\date{Nov. 2001; revised, Sept. 2002}

\maketitle

\begin{abstract}

Let $(M, \partial M)$ be a compact 3-manifold with boundary which admits
a complete, convex co-compact hyperbolic metric. 
For each hyperbolic metric $g$ on $M$ such that $\dr M$ is smooth and
strictly convex, the induced metric on $\dr M$ has curvature $K>-1$, and
each such metric on $\dr M$ is obtained for a unique choice of $g$. 
A dual statement is that, for each $g$ as above, the third fundamental
form of $\dr M$ has curvature $K<1$, and its closed geodesics which are
contractible in $M$ have length $L>2\pi$. Conversely, any such metric on
$\dr M$ is obtained for a unique choice of $g$. 

We are interested here in the similar situation where $\partial
M$ is not smooth, but rather looks locally like an ideal polyhedron in
$H^3$. We can give a 
fairly complete answer to the question on the third fundamental form
--- which in this case concerns
the dihedral angles --- and some partial results about the induced
metric. 

This has some by-products, like an affine piecewise flat structure on
the Teichm{\"u}ller space of a surface with some marked points, or an
extension of the Koebe circle packing theorem to many 3-manifolds with
boundary. 

\bigskip

\begin{center} {\bf R{\'e}sum{\'e}} \end{center}

Soit $(M, \partial M)$ une vari{\'e}t{\'e} compacte de dimension 3 {\`a} bord, qui
admet une m{\'e}trique compl{\`e}te convexe co-compacte. Pour chaque m{\'e}trique
hyperbolique $g$ sur $M$ telle que $\dr M$ est r{\'e}gulier et strictement
convexe, la m{\'e}trique induite sur $\dr M$ est {\`a} courbure $K>-1$;
r{\'e}ciproquement, chaque m{\'e}trique {\`a} courbure $K>-1$ sur $\dr M$ est
obtenue pour un unique choix de $g$. Un {\'e}nonc{\'e} dual est que, pour ces
m{\'e}triques $g$ sur $M$, la troisi{\`e}me forme fondamentale de $\dr M$ est {\`a}
courbure $K<1$, et ses g{\'e}od{\'e}siques ferm{\'e}es qui sont contractiles dans
$M$ sont de longueur $L>2\pi$; r{\'e}ciproquement, chaque m{\'e}trique de ce
type est obtenue pour un unique choix de $g$. 

Nous nous int{\'e}ressons au cas similaire o{\`u} $\dr M$ n'est pas r{\'e}guli{\`e}re,
mais ressemble au contraire localement {\`a} un poly{\`e}dre id{\'e}al 
dans $H^3$. On donne un {\'e}nonc{\'e} assez complet concernant
la troisi{\`e}me forme fondamentale --- qui dans ce cas se formule en termes
d'angles di{\`e}dres --- et un {\'e}nonc{\'e} partiel pour la m{\'e}trique induite sur
le bord. 

Ceci a comme cons{\'e}quence l'existence d'une m{\'e}trique affine plate par
morceaux sur l'espace de Teichm{\"u}ller d'une surface munie de points
marqu{\'e}s, ou une extension du th{\'e}or{\`e}me de Koebe sur les empilements de
cercles {\`a} beaucoup de 3-vari{\'e}t{\'e}s {\`a} bord. 

\end{abstract}

\tableofcontents

\section{Introduction}

\subsection{Hyperbolic manifolds with convex boundary}

We will be motivated by the following results from \cite{hmcb}, which
show that the 
hyperbolic metrics on $M$ such that $\dr M$ is smooth and strictly
convex can be well understood from quantities induced on the boundary.

\btm \label{tm:hmcb-I}
Let $g$ be a hyperbolic metric on $M$ such that $\dr M$ is smooth and
strictly convex. Then the induced metric $I$ on $\dr M$ has curvature
$K>-1$. Each smooth metric on $\dr M$ with $K>-1$ is induced on $\dr M$
for a unique choice of $g$. 
\etm

\btm \label{tm:hmcb-III}
Let $g$ be a hyperbolic metric on $M$ such that $\dr M$ is smooth and
strictly convex. Then the third fundamental form $\III$ of $\dr M$ has
curvature $K<1$, and its closed geodesics which are contractible in $M$
have length $L>2\pi$. Each such metric is obtained for a unique choice
of $g$. 
\etm

Theorem \ref{tm:hmcb-I} was conjectured by Thurston, and its existence
part was proved by Labourie \cite{L4}. Note that the hypothesis of
theorem \ref{tm:hmcb-III} is that the metric on $\dr M$, when
lifted to the boundary of the universal cover of $M$, is globally
$\cat(1)$; this implies a local "curvature" condition, but also a global
condition, namely that all closed geodesics have length larger than
$2\pi$.  

A striking feature of both those
statements is that they appear to hold, at least in some cases, without
any smoothness assumption (beyond 
the convexity). In the simple case where $M$ is a ball, both results can
be stated in terms of convex surfaces in $H^3$:

\btm[Pogorelov \cite{Po}] \label{tm:po}
Let $h$ be a smooth metric with curvature $K>-1$ on $S^2$. Then $(S^2,
h)$ has a unique isometric embedding into $H^3$.
\etm

\btm[\cite{these,cras}] \label{tm:these}
Let $h$ be a smooth metric with curvature $K<1$ on $S^2$. Then $(S^2,
h)$ has an embedding in $H^3$ such that the third fundamental form of
the image is $h$ if and only if all closed geodesics of $(S^2, h)$ have
length $L>2\pi$. This embedding is then unique. 
\etm

Both of those statements also hold without any smoothness assumption; it
is a result of \cite{Po} for theorem \ref{tm:po}, and has been announced
recently by G. Moussong for theorem \ref{tm:these} (see
\cite{davis-moussong}). In particular, the polyhedral case correspond to
earlier results of Aleksandrov \cite{Al} for theorem \ref{tm:po}, and
to results of Andreev \cite{andreev}, Rivin and Hodgson \cite{Ri,RH} for
theorem \ref{tm:these}. 

It also appears that theorems \ref{tm:hmcb-I} and \ref{tm:hmcb-III} hold, at
least in some cases, for convex surfaces which are not complete. This is
again most apparent when $M$ is a ball. In particular:
\begin{itemize}
\item direct analogs of theorems \ref{tm:hmcb-I} and \ref{tm:hmcb-III} hold
  for smooth, non-compact convex surfaces in $H^3$ \cite{rsc}.
\item the analog of theorem \ref{tm:hmcb-I} is true for ideal polyhedra, a
  result of Rivin \cite{rivin-comp}.
\item the analog of theorem \ref{tm:hmcb-III} also holds for ideal
  polyhedra, a result of Andreev \cite{andreev-ideal} and Rivin
  \cite{rivin-annals}. Thurston \cite{thurston-notes} noted that those
  results are also strongly related to questions on circle packings, in
  particular the Koebe theorem. 
\item the analog of theorem \ref{tm:hmcb-I} holds for hyperideal polyhedra
  in $H^3$, see \cite{shu}.
\item theorem \ref{tm:hmcb-III} also holds for hyperideal polyhedra, see
  \cite{shu,bao-bonahon}
\end{itemize}
\medskip

When $M$ is topologically more complicated than a ball, known results
are limited to the "fuchsian" case, i.e. when $M$ is the product of a
surface of genus at least 2 by an interval, with an isometric
involution.  In this case, a version of theorem \ref{tm:hmcb-III},
concerning fuchsian manifolds whose boundary is locally like an ideal
polyhedron, can be found hidden behind work of Thurston \cite[chapter
13]{thurston-notes} and Colin de Verdi{\`e}re \cite{CdeV} on circle packings
on hyperbolic surfaces; the relationship (which is well known) should be
clear from section 3 below. 

The main goal of this paper is to extend this result to the situation
where $M$ is topologically "general" --- but its boundary will always be
supposed to be locally like an ideal polyhedron, in a sense which is
precisely defined in section 3. In this case, we will call $M$ an {\bf
  ideal hyperbolic manifold}. 

We need the following definition. A {\bf
  cellulation} of $\dr M$ is a decomposition of $\dr M$ into the images
by diffeomorphisms of convex polygons in $\R^2$. We
will call the images of the polygons the {\bf 2-cells} of the
cellulations, and the  {\bf 1-cells} will be the images of the edges
of the polygons. We demand that the cellulations are well behaved in
the sense that: 
\begin{itemize}
  \item the 2-cells have disjoint interiors.
  \item two 1-cells which are distinct have disjoint interiors.
  \item the intersection of two 2-cells is a disjoint union of 1-cells and
    vertices.  
  \item at least 3 images of polygons meet at the image of each vertex. 
\end{itemize}

\bdf 
Let $\Gamma$ be the 1-skeleton of a cellulation of $\dr M$. A {\bf
  circuit} in $\Gamma$ is a sequence $e_0, e_1, \cdots, e_n=e_0$ of
edges of $\gamma$ such that the dual edges $e_0^*, e_1^*, \cdots, 
e_n^*$ are the
successive edges of a closed path in the dual graph $\Gamma^*$, which
is contractible in $M$. A
circuit is {\bf elementary} if the dual path bounds a face
(i.e. containing no edge or vertex). 
\edf

We can now describe the dihedral angles of ideal hyperbolic manifolds;
it extends the results of \cite{andreev-ideal,rivin-annals} on ideal
polyhedra.  

\btn{\ref{tm:dihedral}}
Suppose that $M$ has incompressible boundary. 
Let $\Gamma$ be the 1-skeleton of a cellulation of $\dr M$.
Let $w$ be a function from the edges of $\Gamma$ to
$(0,\pi)$ such that:
\begin{enumerate}
\item for each elementary circuit in $\Gamma$, the sum of the values of
  $w$ is $2\pi$;
\item for each non-elementary circuit in $\Gamma$, the sum of the values
  of $w$ is strictly larger than $2\pi$. 
\end{enumerate}
Then there a unique hyperbolic metric $g$ on $M$, such that $(M,g)$ is
an ideal hyperbolic manifold, with exterior dihedral angles given by
$w$.
\etn

The hypothesis that $M$ has incompressible boundary is necessary for
technical reasons (see section 8). Theorem \ref{tm:hmcb-I} would seem to
indicate that the result also holds without this hypothesis. 

Note that the hypothesis of this theorem on the function $w$, implies
topological properties on $\Gamma$. In particular:

\blm \label{lm:topolo}
Let $\sigma$ be a cellulation of $\dr M$, such that their is a function
$w$ on the edges of $\sigma$, satisfying the hypothesis of theorem
\ref{tm:dihedral}. Then:
\begin{itemize}
\item each 1-cell, when lifted to the universal cover of $\dr M$, is
  topologically a segment. 
\item each 2-cell, when lifted to the universal cover of $\dr M$, is
  topologically a disk. 
\end{itemize}
\elm

The proof is in section 3.

\subsection{Circle packings}

Theorem \ref{tm:dihedral} has a consequence in terms of circle packing
on the boundary of 3-manifolds. Remember that a {\bf circle packing} on the
sphere $S^2$ is a set of closed disks with disjoint interiors in
$S^2$. Given a circle packing, one can define its {\bf incidence graph}
as a graph on $S^2$ which has one vertex for each disk, and an edge
between two vertices if and only if the corresponding disk are tangent. 

The classical Koebe circle packing theorem states that, for each graph
$\gamma$ in $S^2$ which is the 1-skeleton of a triangulation, there exists
a unique circle packing whose incidence graph is $\Gamma$. It was proved
by Koebe \cite{koebe} for triangulations, and extended by Thurston
\cite{thurston-notes} using the Andreev
theorem on ideal polyhedra \cite{andreev-ideal}; in this more general
case one should demand that, for any connected component of the
complement of the disks, there is  a circle which is orthogonal to all
the neighboring disks. 

Thurston also extended this theorem to circle packings on hyperbolic
surfaces. Theorem \ref{tm:dihedral} provides a further extension to the
boundary of our 3-manifold $M$, when it is provided with the
$CP^1$-structure on the boundary coming from a hyperbolic metric on
$M$. 

\btn{\ref{tm:koebe}}
Suppose that $M$ has incompressible boundary. 
Let $\Gamma$ be the 1-skeleton of a triangulation of $\dr M$. There is a
unique couple $(g, c)$, where $g$ is a complete, convex co-compact
hyperbolic metric on $M$, and $c$ is a circle packing on $\dr M$ (for
the $CP^1$-structure defined on $\dr M$ by $g$) whose incidence graph
is $\Gamma$. 
\etn

If one considers not only triangulation but more generally cellulations,
the same result holds, with the additional condition that, for each
connected component of the complement of the closed disks, there exists
a circle which is orthogonal to all the 
adjacent circles --- this condition is automatically satisfied for a
triangulation. 

The proof is given in section 10.

\subsection{Induced metrics on the boundary}

It would be interesting to know whether theorem \ref{tm:hmcb-I} extends to
ideal hyperbolic manifolds. Unfortunately I do not know the answer to
this question; but an infinitesimal version does hold:

\blm \label{lm:01}
For each ideal hyperbolic manifold $M$, each infinitesimal variation of
the hyperbolic structure on $M$ (among ideal hyperbolic manifolds) induces
a non-trivial infinitesimal variation of the induced metric on $\dr M$.
\elm

This assertion is worth mentioning because, in other questions
related to theorem \ref{tm:hmcb-I}, it is precisely this
infinitesimal rigidity which is lacking (see section 1).

However, in this case it is not sufficient to obtain a global result
concerning the induced metrics on the boundary of ideal hyperbolic
manifolds, like in theorem \ref{tm:dihedral} for the third fundamental
form. The reason is that, to understand the induced metrics completely,
one would have to go beyond the category of ideal hyperbolic manifolds,
and obtain an infinitesimal rigidity result also for "bent" manifolds
(as defined in section 3). We give more details on this in section 11.

We will also give below a more precise result in the special case of
manifolds which we call "fuchsian"; they are ideal hyperbolic manifolds
such that $\dr M$ has two connected component, with an isometric
involution exchanging the connected components of $\dr M$.

\subsection{Fuchsian polyhedra}

In the course of the proof of this theorem, we are led to study
fuchsian equivariant polyhedra; they are the objects
arising as the universal covers of the 
boundaries of the fuchsian manifolds mentioned above. As a consequence
of theorem \ref{tm:dihedral}, we find a characterization of the dihedral
angles of fuchsian ideal polyhedra. In addition, we will give in
section 4 results concerning other fuchsian polyhedra, having some
non-ideal vertices. 

\bdf
A {\bf fuchsian polyhedron} is a triple $(S, \phi, \rho)$, where:
\begin{itemize}
\item $S$ is a surface of genus $g\geq 2$, with $N$ marked points $x_1,
  x_2, \cdots, x_N$, $N\geq 1$.
\item $\phi$ is a polyhedral map from the universal cover $\St$ of $S$
  into $H^3\cup \dr_\infty H^3$ such that the image is locally like a
  polyhedron in $H^3$, with vertices the images by $\phi$ of the inverse
  images in $\St$ of the $x_i$.
\item $\rho$ is a morphism from $\pi_1S$ into the subgroup of the
  isometry group of $H^3$ of isometries fixing a totally geodesic
  2-plane $P_0$.
\item for any $x\in \St$ and any $\gamma\in \pi_1S$, $\phi(\gamma
  x)=\rho(\gamma)\phi(x)$. 
\end{itemize}
\edf

Some special cases are of interest. Let $I$ be the set of points in
$\St$ which are sent by the canonical projection $\St\rightarrow S$ to a
marked point. We say that $(S, \phi, \rho)$ is:
\begin{itemize}
\item ideal, if each point in $I$ is sent to an ideal point.
\item finite, if no point in $I$ is sent to an ideal point. 
\item semi-ideal, if the points in $I$ can be sent to either ideal or
non-ideal vertices (this includes the finite and the ideal cases).
\end{itemize}
We obtain an existence result for the third fundamental forms of
fuchsian polyhedra. 

\btn{\ref{tm:semi-ideal}}
Let $\Sigma$ be a surface of genus $g\geq 2$, and let $h$ be a
spherical cone-metric on $\Sigma$, with negative singular curvature at
the singular points. Suppose that all contractible 
closed geodesics of $(\Sigma, h)$
have length $L>2\pi$, except when they bound a hemisphere. 
Then there is a unique fuchsian polyhedral
embedding of $(\Sigma, h)$ into $H^3$ whose third
fundamental form is $h$.
\etn

\medskip

As a consequence of the analysis of the dihedral angles of ideal
hyperbolic manifolds in the fuchsian case, we will obtain the following.

\btn{\ref{tm:teichmuller-affine}}
For each $g\geq 2$ and each $N\geq 1$, there is a natural unimodular
piecewise affine structure $\cA_{g,N}$ on
the Teichm{\"u}ller $\cT_{g,N}$ space of the genus $g$ surface with $N$
marked points.
\etn

This affine structure has some singularities, but also a number of
interesting properties; for instance there is a natural function on
$\cT_{g,N}$, which is defined by taking the volume of a fuchsian ideal
hyperbolic manifold, and it is concave on the maximal dimension cells of
$\cA_{g,N}$. There are also
some questions which remain open concerning $\cA_{g,N}$.

\medskip

Finally, in the case of fuchsian polyhedra, the infinitesimal rigidity
result for the induced metrics (lemma \ref{lm:01}) happens to be
sufficient to obtain a satisfactory global result on the induced
metrics. 

\btn{\ref{tm:metrique-fuchsian}}
Let $S$ be a surface of genus $g\geq 2$, and let $N\geq 1$. 
For each complete, finite area
hyperbolic metric $h$ on $S$ with $N$ cusps, there is a unique ideal
fuchsian hyperbolic manifold $M$ such that the induced metric on each
component of the boundary is $h$.
\etn

The proof is in section 9. 

\subsection{Some hints on the techniques}

The main idea behind the results obtained here is already a few years
old: it is the use
of the striking properties of the volume function on 
the space of hyperbolic structures. The first property is that the
volume is a concave function when the space of hyperbolic structures is
parametrized by the dihedral angles (this comes from the elementary
lemma \ref{lm-vol} concerning the ideal simplex) and the second is the
Schl{\"a}fli formula, which provides a fundamental link between the dihedral
angles and the induced metric on ideal polyhedra, or on hyperbolic
manifolds with ideal-like boundary.

Those ideas can be traced back to the treatment of the Andreev
theorem in \cite{thurston-notes}; the existence of a "generating
function" for the problem was then pointed out in \cite{CdeV}, and this
function was identified as the volume in \cite{bragger} and \cite{Ri2}. 

Following this path for hyperbolic manifolds with ideal-like boundary,
however, leads to some technical difficulties. The first is that it is
not clear how one can find an ideal triangulation of the manifold; by
the way, a similar 
difficulty appears in other situations where one wants to deform
hyperbolic structures, see e.g. \cite{porti-petronio}. This is treated in
section 5, where it is shown that some finite cover of the manifolds
considered have an ideal triangulation. All the arguments can then be
given in this finite cover, and an equivariance argument brings the
result down to the manifolds we want to study.

Another technical point is to prove that the set of possible dihedral
angle assignations (as it appears in
theorem \ref{tm:dihedral}) is connected; this is necessary because the
proofs rely on a deformation argument. The solution chosen here uses the
fact that the conditions are "independent" on each of the connected
components of the boundary. Moreover, if $M$ has incompressible
boundary, then the conditions on each connected component of $\dr M$ is
the same as for the corresponding fuchsian ideal polyhedra; so that the
connectedness in the general 
case follows from the understanding of the ideal fuchsian polyhedra. So
those fuchsian polyhedra are studied in section 4. The methods used
there are quite different from those of the  
other sections; an existence and uniqueness result is proved for the
third fundamental forms of  fuchsian manifolds whose boundary locally
looks like a compact polyhedron in $H^3$, 
and from there one deduces an existence result for manifolds with
ideal-like boundary by an approximation argument. 

Although some other technically interesting details appear at different
points in the paper, it does not seem necessary to describe them here ---
the reader will probably enjoy finding them as he reads along.

Even in the case of ideal polyhedra in $H^3$, I think that the approach
used here is different, and I believe slightly more direct, than the one
used in previous papers --- although no new result is achieved. 
Moreover, I hope that the methods used --- as mentioned above, 
mostly developed earlier
to study ideal polyhedra --- could hint at some approaches to the
questions concerning the convex cores of complete, convex co-compact
hyperbolic manifolds. Of course many difficulties
remain in this direction, although recent works of Bonahon
\cite{bonahon,bonahon2} might provide useful tools.

Note that some of the results presented here --- notably theorem
\ref{tm:dihedral} --- have a non-empty intersection with some of the
results obtained by Rivin in \cite{rivin-comp}. The methods used, and
their scope, however, are quite different.

\subsection{What follows}

Section 2 contains a reminder of the classical Schl{\"a}fli formula,
as well as some remarks and interpretations in terms of symplectic
geometry. Section 3 then gives some basic definitions and 
results on the geometry of
ideal hyperbolic manifolds. 

The case of fuchsian manifolds is described in section 4, using methods
that are much closer to those used in the classical theory of hyperbolic
polyhedra (as developed in particular by Aleksandrov \cite{Al}), and
some methods from \cite{iie}.

Section 5 deals with some technical
questions on triangulations of ideal hyperbolic manifolds; the point is
that, although I do not know how to construct a well-behaved
triangulation of such a manifold, it is not too difficult to
construct one on a finite cover, and this will be sufficient for this
paper. 

The geometrical constructions start in section 6, where more general
hyperbolic structures on triangulated manifolds are investigated. The
idea --- which basically comes from earlier works, see
\cite{thurston-notes,CdeV,bragger,rivin-annals} is to use variational
properties of the volume functional on a larger class of hyperbolic
structures to obtain hyperbolic metrics. This is continued in section 7,
where first-order properties of the volume are described.

This leads in section 8 to results concerning dihedral angles, and in
section 9 to results on induced metrics on the boundary. Some
applications to circles packings are mentioned in section 10, and other
considerations stand in the last section.

\section{The Schl{\"a}fli formula}

We recall in this section the classical Schl{\"a}fli formula. It is the main
tool used in the sequel, so we also describe some interesting
interpretations of it and some related properties of the volume of
simplices. 

We first state the Schl{\"a}fli formula for compact polyhedra; for a proof,
see e.g. \cite{milnor-schlafli} or \cite{geo2}.

\blm \label{lm-schlafli}
Let $(P_t)$ be a one-parameter family of convex polyhedra in $H^3$. Let
$I$ be the set of its edges, $(L_i)_{i\in I}$ be its edge lengths, and
$(\theta_i)_{i\in I}$ the corresponding (interior) dihedral angles. Then:
\beq \label{schlafli}
dV = -\frac{1}{2}\sum_{i} L_i d\theta_i~. \eeq
\elm

\subsection{A symplectic viewpoint}
There is an amusing symplectic interpretation of this formula. 
To explain it simply we choose an abstract polyhedron $P_0$, with $v$
vertices, $e$ edges and $f$ 2-faces, and
consider only polyhedra with the same combinatorics.
Call $\cL:=\R^e$ and $\Theta:=(0,\pi)^e$ the sets containing the possible
lengths and exterior dihedral angles, respectively, of hyperbolic polyhedra
having the same combinatorial type as $P_0$.
Now consider
the symplectic vector space $(\cL\times \Theta, \omega)$, 
with $\omega := \sum_i
dL_i \wedge d\theta_i$. Let $\cP$ be the subset corresponding to the
edge lengths and dihedral angles of convex polyhedra of the same
combinatorial type as $P_0$. It is well known, and not difficult to
prove, that:

\bprop 
$\cP$ is a submanifold (with boundary) of $\cL\times \Theta$ of dimension $e$.
\eprop

We leave the proof to the reader, it uses the Euler formula and the fact
that the number of constraints on the positions of the vertices is the
weighted number of non-triangular faces. 

The Schl{\"a}fli formula is essentially equivalent to the following:

\bcn{to the Schl{\"a}fli formula}
$\cP$ is lagrangian in $(\cL\times \Theta, \omega)$.
\ecn

\bpv
Define the 1-form:
$$ \beta := \sum_i L_i d\theta_i~. $$
Then $d\beta=\omega$.
Since $\beta_{|\cP}=dV_{|\cP}$, $d\beta_{|\cP}=0$ so that $\omega$
vanishes on $\cP$.
\epv

\subsection{Ideal polyhedra}
Now let $\cP_\infty$ be the set of ideal polyhedra having the
combinatorial type of $P_0$; they are obtained by letting the vertices
go to infinity -- this is possible for many combinatorial types. A basic
remark is that those polyhedra still have finite volume. Their edge
lengths, however, are not defined as such (they are infinite). To define
an analog of the edge lengths for ideal polyhedra, we choose for each
vertex $V$ a horosphere $H_V$ ``centered'' on $V$, that is, a level set
for the Busemann function associated to $V$. We then define the length of
the edge joining 2 vertices $V$ and $W$ as the distance, along the edge,
between $H_V$ and $H_W$; we use the signed length, so that the length is
negative if the horospheres overlap. The set of edge lengths of the
elements of $\cP_\infty$ is thus defined up to the addition of a
constant for each vertex, and is contained in the set
$\cL_\infty:=\R^e/\R^v$. 

On the other hand, the dihedral angles of the ideal polyhedra are
constrained. That is because the link of each vertex is a Euclidean
polygon, so the sum of its exterior angles is $2\pi$, so that the sum of
the exterior dihedral angles of the edges containing a given vertex of
an ideal polyhedron is always $2\pi$. Therefore, the set of dihedral
angles of the polyhedra in $\cP_\infty$ stays in a $(e-v)$-dimensional
space $\Theta_\infty$.

The amusing fact that we announced is then:

\brk \label{rk:ideal}
\begin{enumerate}
\item $\cL_\infty\times \Theta_\infty$ is obtained from $\cL\times 
  \Theta$ by symplectic reduction, and it therefore still carries a
  symplectic form $\omega_\infty$;
\item the Schl{\"a}fli formula (\ref{schlafli}) still has a meaning, and
  still holds, for $\cP_\infty$; 
\item it still implies that $\cP_\infty$ is a lagrangian submanifold of
  $(\cL_\infty\times \Theta_\infty, \omega_\infty)$.
\end{enumerate}
\erk

\bpv
The second point is probably classical, and we leave it to the reader. 
For the first point, consider the group $G:=\R^v$, and its action $\phi$
on $\cL\times \Theta$ by:
\begin{eqnarray}
  \label{eq:action}
  \phi : & G=\R^v \times (\cL\times \Theta) & \rightarrow \cL\times \Theta
  \nonumber \\
  & ((\alpha_j)_{j=1, \cdots, v}, (L_i, \theta_i)_{i=1,\cdots, e}) &
  \mapsto (L_i+\alpha_{i_+}+\alpha_{i_-}, \theta_i)_{i=1,\cdots, e} 
  \nonumber
\end{eqnarray}
where $i_+$ and $i_-$ are the two ends of the edge $i$.

It is quite easy to check that $\phi$ has a moment map $\mu $ defined by:
\begin{eqnarray}
  \label{eq:moment}
  \mu : & \cL\times \Theta & \rightarrow \cG^*=\R^v \nonumber \\
  & (L_i, \theta_i) & \mapsto \left(\sum_{k \in
  E_j}\theta_k\right)_{j=1,\cdots,v} \nonumber 
\end{eqnarray}
where $E_j$ is the set of edges containing a vertex $j$.

Then $\cL_\infty\times \Theta_\infty \simeq \mu ^{-1}(2\pi,\cdots, 2\pi)/G$,
where $G$ acts by $\phi$, so that $\cL_\infty\times \Theta_\infty$ is obtained
by symplectic reduction from $\cL\times \Theta$ as announced, with a
symplectic form $\omega_\infty$.

For the last point note that $\beta$ determines a well-defined 1-form
$\beta_\infty$ on $\cL_\infty\times \Theta_\infty$, so again
$\omega_\infty=d\beta_\infty$. But the Schl{\"a}fli formula for ideal
polyhedra shows that $\beta_\infty$ vanishes on $\cP_\infty$, which is thus
lagrangian. 
\epv

\subsection{The volume function}
We will also need some elementary and well-known properties of the
volume of hyperbolic simplices, which we recall here for the reader's
convenience. More details can be found for instance in
\cite{thurston-notes}, chapter 7.

\bdf \label{df-loba}
The Lobachevsky function is defined as:
$$ \Lambda(\theta) := - \int_0^\theta \log |2\sin u | du~. $$
\edf

Now recall that there is a 2-parameter family of ideal simplices in
$H^3$ (up to global isometries), which can be parametrized for instance
by the complex cross-product of the four vertices in
$\dr_\infty H^3\simeq \CP^1$. For each ideal simplex, the dihedral angles
of two opposite edges are equal, and the sum of the exterior dihedral
angles of the edges containing a given vertex is $2\pi$. An ideal
simplex is completely determined -- again up to global isometry -- by
its three interior dihedral angles $\alpha, \beta$ and $\gamma$, under the
condition that their sum is $\pi$.

The volume of an ideal simplex is given by a simple formula (see
e.g. \cite{thurston-notes}, chapter 7):

\blm \label{lm-vol}
The volume of the ideal simplex with dihedral angles $\alpha, \beta$ and
$\gamma$ is $\Lambda(\alpha)+\Lambda(\beta)+\Lambda(\gamma)$. 
\elm

As a consequence, one finds (see e.g. \cite{Ri2}):

\bcr \label{cr:concavity}
$V$ is a concave function of the dihedral angles $\alpha, \beta,
\gamma$, which vanishes when one of the angles goes to $0$.
\ecr

\bpv 
This is proved by an elementary computation of the Hessian of $V$. 
\epv

Note that the Schl{\"a}fli formula is not restricted to hyperbolic
polyhedra; it is also valid in the other Riemannian or pseudo-Riemannian
space-forms (see \cite{suarez}) and also has an interesting extension to
the setting of Einstein manifolds with boundary, see
\cite{sem,sem-era,bernig}. I do not know whether any analog of the methods
described in this paper -- where the volume is used as a "generating
function" to prove geometric results -- can be found in this more
general context.

There is another, completely elementary proof of corollary
\ref{cr:concavity}, which uses the Schl{\"a}fli formula instead of explicit
computations of the volume in terms of the Lobachevsky function. It is
also more general. The key point is that, according to the
Schl{\"a}fli formula, the infinitesimal deformations of an ideal simplex
which do not change the lengths of the edges (i.e. the isometric
deformations) are exactly the deformations which are in the kernel of
the Hessian of the volume, seen as a function of the dihedral
angles. Now it is easy to see --- and classical --- that ideal simplices
are rigid, i.e. they have no non-trivial infinitesimal isometric
deformation. Therefore, the Hessian of the volume has constant
signature. But an elementary argument (with the Schl{\"a}fli formula) 
shows that the regular simplex has maximal volume, so that $V$ is a
concave function. 

This line of reasoning can of course be used in different contexts; for
instance, it shows that the signature of the Hessian of the volume, as a
function of the dihedral angles, has constant signature on the space of
compact hyperbolic simplices, since they are also known to be rigid.

\section{Convex hulls of ideal points}

The goal of this paper is to understand hyperbolic manifolds whose
boundaries look locally like ideal polyhedra. We will introduce a class
of manifolds with boundary obtained by taking the convex hull of a
finite number of ideal points in a complete, convex co-compact
hyperbolic manifold. The main point is that the class of those manifolds
with boundary separates into two sub-class, according to whether their
"convex core" intersects their boundary or not. We will then indicate
why both cases 
actually exist in a non-trivial way, since it is not completely obvious
at first sight.

\subsection{Some definitions}
We first define a notion of convexity --- it is classical but other
definitions are sometimes used.

\bdf \label{df:convexite}
Let $M$ be a hyperbolic manifold. We say that a subset $C\subset M$ is
{\bf convex} if, 
for any points $x$ and $y$ in $C$, any geodesic segment in $M$
with endpoints $x$ and $y$ remains in $C$. For any subset $E$ of $M\cup
\dr_\infty M$,
the {\bf convex hull} of $E$ is the smallest convex set containing $E$.
\edf

\bdf \label{horns}
$M$ is a {\bf hyperbolic manifold with horns} if there is a complete
co-compact hyperbolic manifold $N$ and a finite family $x_1, \cdots,
x_p$ of points in $\dr_\infty N$ such that $M$ is isometric to the
convex hull in $N$ of $\{ x_1, \cdots, x_p\}$.
\edf

For instance, the interior of an ideal polyhedron in $H^3$ is a
hyperbolic manifold with horns. So is the convex core of a convex co-compact
hyperbolic manifold; in this case $p=0$.

\brk
Let $M$ be a compact hyperbolic manifold with convex boundary, such that the
induced metric on the boundary is a complete hyperbolic metric of finite
area on the boundary minus a finite number of points. Then $M$ is
a hyperbolic manifold with horns.
\erk

\bpv
Since $M$ is compact with convex boundary it is isometric to a subset of
a complete convex co-compact hyperbolic manifold $N$. Let $\Mt$ be the
universal cover of $M$, then $\dr \Mt$ is a convex surface in $H^3$ and
its induced metric is hyperbolic; it is therefore the convex hull in
$H^3$ of its boundary points in $\dr_\infty H^3$. Taking the quotient by
$\pi_1M$, we see that $\dr M$ is the convex hull of its boundary points
in $\dr_\infty M$. Those points clearly correspond to the cusps of $\dr
M$, so $M$ is the convex hull in $H^3/\pi_1M$ of a finite number of ideal
points. 
\epv

Consider a hyperbolic metric on $M$ for which $\dr M$ is convex. 
There is then a unique convex co-compact hyperbolic 3-manifold $N$ in which
$M$ admits an isometric embedding which is surjective on the
$\pi_1$; we call it the {\bf extension} of $M$, and denote it by $E(M)$. 
The convex core of $E(M)$ can be defined as the smallest convex subset
of $E(M)$, so it is contained in $M$; we will denote it by $C(M)$.

\bdf \label{df:ideal}
Let $M$ be a hyperbolic manifold with horns. $M$ is an {\bf ideal hyperbolic
  manifold} if $C(M)\cap\dr M=\emptyset$, otherwise $M$ is a {\bf bent
  hyperbolic manifold}.
\edf
 
For instance, the interiors of ideal polyhedra in $H^3$ are of the ideal
kind, while the convex cores of convex co-compact manifolds are of the
bent type. Most of what follows concerns ideal hyperbolic manifolds
only, while the bent case appears as the problem lurking in the
background. The relationship between the two will be explored after the
next subsection, using the tools that it contains. 

\subsection{Relation with circle packings}

To understand the behavior of the boundary of a hyperbolic
manifold with horns $M$, it is relevant to consider its universal cover $\Mt$,
which can naturally be identified with a convex subset of $H^3$. The
vertices of $\Mt$ correspond to ideal points in $H^3$, i.e. to points in
$\dr_\infty H^3$. The combinatorics of the faces of $\dr \Mt$  is then
described in terms of Delaunay cellulations of the sphere.

\bdf \label{df:delaunay}
Let $\sigma$ be a cellulation of $S^2$. 
$\sigma$ is {\bf Delaunay} if, for each cell $c$ of $\sigma$:
\begin{enumerate}
\item the vertices of $c$ are co-cyclic, i.e. lie on a circle $C(c)$ in
  $S^2$; 
\item one of the closed disks bounded by $C(c)$ contains no other
  vertex of $\sigma$.
\end{enumerate}
\edf

Note that this definition involves only the vertices of $\sigma$, not
its edges. It depends on the M{\"o}bius structure of $S^2$, but not on its
metric structure. 

\blm \label{lm:36}
Let $E$ be a discrete subset of $S^2$. There is a unique maximal
Delaunay cellulation of a subset of $S^2$ whose vertices are the
elements of $E$. It combinatorics is that of the convex hull of $E$ in
$\R^3$. 
\elm

The proof is elementary; the circles appearing in the definition of a
Delaunay cellulation are the boundaries at infinity of support
planes of the convex hull of $E$.

Let $M$ be a hyperbolic manifold with horns. By definition, $M$ is the
convex hull of a finite set of ideal points $S=\{ x_1, \cdots,
x_p\}\subset \dr_\infty E(M)$. Consider the universal cover
$\tilde{E(M)}=H^3$ of $E(M)$; $S$ lifts to a set of points $\St$ which
is invariant under the action of $\pi_1(M)$, so that the accumulation set
of $\St$ is the limit set $\Lambda\subset S^2$ of the action of $\pi_1M$
on $H^3$. 


Consider a cellulation $\sigma$ of a subset $\Omega$ of $S^2$; we
consider it as of a partly geometric and partly combinatorial nature, with
fixed vertices in $S^2$ but edges defined up to isotopy. 

The following statement is a consequence of lemma \ref{lm:36}.

\bprop 
There is a unique Delaunay cellulation of $S^2\setminus \Lambda$ with
vertices the elements of $\St$; it is obtained as the combinatorial
structure of the convex hulls of $\St$ in $H^3$.
\eprop

The circles appearing in the definition of a Delaunay cellulation are
simply the traces on $\dr_\infty H^3$ of the 2-planes which are the
faces of the convex hull of $\St$. Moreover, the dihedral angles between
those faces are the angles between the corresponding circles in
$S^2=\dr_\infty H^3$. Thus, dihedral angle questions on ideal hyperbolic
manifolds can be translated as questions on configurations of circles in $S^2$
having given angles (and invariant under group actions). A precise
relationship with circle packings can be obtained using an idea of
Thurston; this is recalled in section 10.

We can now show that ideal hyperbolic manifolds have a simple
description in terms of the properties of their boundary.
Of course this description does not apply to bent hyperbolic manifolds. 

\bpt \label{pt:ideal}
Let $M$ be an ideal hyperbolic manifold. Then its boundary $\dr M$ is
the union of a finite number of 2-faces, which are ideal polygons with a
finite number of edges in
totally geodesic 2-planes, and which intersect along geodesics. 
\ept

\bpv
Since $E(M)$ is convex co-compact, its convex core $C(M)$ is compact. By
definition, $\dr M$ does not intersect $C(M)$, so $d(\dr M,
C(M))>0$. Therefore, each face of $\dr M$ remains at a positive distance
from $C(M)$. Thus each face of $\dr \Mt$ remains in a compact subset of
the complement of the limit set $\Lambda$ of the action of $\pi_1 M$ on
$H^3$. Therefore, the circles in $S^2=\dr_\infty H^3$ which are the
boundary at infinity of those planes are not tangent to $\Lambda$.

Let $\Omega$ be a fundamental domain with a compact closure
in $S^2\setminus \Lambda$ for the
action of $\pi_1M$. $\Omega$ is contained in a compact subset $K\subset
S^2\setminus \Lambda$, so it intersects a finite number of those
circles; therefore, $\dr M$ has a finite number of faces. 
The proof of the property follows.
\epv



\subsection{Existence of the bent case}

Let $N$ be a complete, convex co-compact hyperbolic 3-manifold. Let
$x_1, \cdots, x_p$ be a family of points in $\dr_\infty N$. If one of
the connected component $\dr_0N$ of $\dr_\infty N$ does not contain any of the
$x_i$, it is not difficult to see that the smallest convex subset of $N$
containing the $x_i$ will be of the bent type; indeed, its boundary will
contain the component of the boundary of the convex core of $N$ facing
$\dr_0N$. We
will sketch here an argument intended to convince the reader that there
are other situations where bent manifolds appear. 

A first remark is that a hyperbolic manifold with horns is bent if and only
if, among the circles in $S^2\setminus \Lambda$ which are the boundary
at infinity of its faces, one has non-empty intersection with
$\Lambda$. Indeed it was already proved in property \ref{pt:ideal} that
this does not happen for an ideal hyperbolic manifold, while the
converse is clear since, for such a circle, there should be a face of
$\dr M$ which is at distance $0$ from $C(M)$. 

Thus, to show the existence of a bent hyperbolic manifold, it is
simplest to search for one among hyperbolic manifolds with  horns with only
one ideal point in each boundary component (those manifolds could
legimitly be called "unicorn manifolds"). Moreover we can restrict our
attention to e.g. the quasi-fuchsian case --- we will see below
that fuchsian hyperbolic manifolds can not be bent.

Fix one of the boundary
components of $M$, say $\dr_1M$, and let $x_1\in \dr_\infty M$ be the
corresponding ideal point. One of the faces of $\dr_1 M$ has non-empty
intersection with $\Lambda$ if and only if there is a circle $C_0$ in
$S^2$, whose interior is in $S^2\setminus \Lambda$, which intersects
$\Lambda$, and whose interior contains no point of the orbit $(\pi_1
M).x_1$.  

So, to prove the existence of a bent hyperbolic manifold, it is
enough to prove the existence of a closed disk $C_0$ which:
\begin{enumerate}
\item has its interior in $S^2\setminus \Lambda$.
\item intersects $\Lambda$.
\item contains no fundamental domain for the action of $\pi_1M$ on
  $S^2$. 
\end{enumerate}
If such a circle exists, there will be a point $x_1$ whose orbit $(\pi_1
M).x_1$ does not intersect the interior of $C_0$. 

To simplify a little the picture, suppose that $M$ has a convex core
whose pleating locus contains a closed geodesic $\gamma$, with a non-zero
pleating angle. In $H^3$, $\gamma$ lifts to a geodesic $\gammat$ with
endpoints $p_1, p_2\in \Lambda$; both $p_1$ and $p_2$ then correspond to
"spikes" of $\Lambda$. Moreover, $\gammat$ lies in the boundary of
$\tilde{C(M)}$, so there is a support plane $P$ of $\tilde{C(M)}$ along
$\gammat$. Let $C_0$ be the corresponding disk in $S^2=\dr_\infty
H^3$. By construction, the interior of $C_0$ does not intersect
$\Lambda$ (while $C_0\cap\Lambda\supset \{ p_1,p_2\}$).

Call $\Omega_1$ the connected component of $S^2\setminus \Lambda$ which
corresponds to $\dr_1M$. 
Let $g_1$ be the hyperbolic metric on $\Omega_1$, conformal to the
canonical metric $g_0$ on $S^2$ and invariant under the
action of $\pi_1M$. Then, by a classical result of conformal geometry
(see e.g. \cite{ahlfors}), the
conformal factor between $g_0$ and $g_1$ is bounded between $c/r$ and
$C/r$, where $c$ and $C$ are two positive constants and $r$ is the
distance to $\Lambda$ in the metric $g_0$.

Therefore, a simple computation shows that the interior of $C_0$ is
contained in a neighborhood of a geodesic in $(\Omega_1, g_1)$, and also
--- by taking the quotient by $\pi_1M$ --- in $(\dr_1M, g_1)$. Thus,
after taking a finite cover of $\dr_1M$, the interior of $C_0$ contains
no fundamental domain for the action of $\pi_1M$. Taking the
corresponding finite cover $\Mb$ of $M$, we see that $\Mb$ is a bent
hyperbolic manifold. 

\medskip

On the other hand, this argument directly shows that bending laminations
can not occur in the case of fuchsian hyperbolic manifolds. 

\brk \label{rk:fuchsian-bent}
Fuchsian hyperbolic manifolds with horns are ideal hyperbolic manifolds.
\erk

\bpv
When $M$ is a fuchsian manifold with horns, $\Lambda$ is a circle in
$S^2=\dr_\infty H^3$. Circles tangent to $\Lambda$ have interiors which
are isometric --- for the hyperbolic metric on the interior of $\Lambda$
--- to horoballs in $H^2$. Therefore they always contain fundamental
domains for any co-compact action on $H^2$.
\epv

\subsection{Necessary conditions on convex surfaces}
\label{ss34}

We recall here some well-known properties of the induced metric and
third fundamental form of convex surfaces in $H^3$; they are classical
except for the length condition on the geodesics of the third
fundamental form, which was understood more recently.

To understand the third fundamental form of non-smooth convex surfaces,
it is helpful to know how the duality between $H^3$ and $S^3_1$ works,
so we will recall it rapidly here; see e.g. \cite{thurston-notes,RH,shu}
for more details. Both spaces can be seen as quadrics
in the Minkowski 4-space $\R^4_1$, with the induced metric:
$$ H^3 = \{ x\in \R^4_1 ~ | ~ \langle x, x\rangle =-1 ~ \wedge ~ x_0>0
\}~, $$
$$ S^3_1 = \{ x\in \R^4_1 ~ | ~ \langle x, x\rangle =1\}~. $$
For $x\in H^3$, let $D$ be the line in $\R^4_1$ going through $0$ and
$x$, and let $D^\perp$ be its orthogonal for the Minkowski inner
product, so that $D^\perp$ is a space-like plane. Then define the dual
$x^*$ of $x$ as the intersection of $D^\perp$ with $S^3_1$, which is a
totally geodesic space-like plane in $S^3_1$. Similarly, the dual of a
point in $S^3_1$ is an oriented plane in $H^3$. The dual of a convex
polyhedron $P$ in $H^3$ is the polyhedron in $S^3_1$ whose vertices are
the duals of the faces of $P$, and whose faces are the duals of the
vertices of $P$ (note that there are other approaches of this duality,
which might actually be more illuminating; see e.g. \cite{shu}).

Given a locally convex surface $S$ in $H^3$, we can define its dual
$S^*$ as
the set of points in $S^3_1$ which are duals of an (oriented) support
plane to $S$. It happens to be another locally convex surface, which is
not necessarily smooth. If $S$ is smooth, then $S^*$ is smooth when $S$
is locally strictly convex. 

\bpt
Let $S$ be a smooth locally strictly convex surface in $H^3$; its third
fundamental form is the induced metric on the dual surface.
\ept

It is tempting to speak of the "third fundamental form" of a non-smooth
surface in the sense "the metric induced on its dual". Some care is
needed, however. For instance, if $S$ is a connected component of the
boundary of the convex core of a convex co-compact manifold, its dual is
a graph; more generally, when $S$ is a convex, developable surface, like
the boundary of a horned hyperbolic manifold, its
dual is one-dimensional. 

On the other hand, this notion of third fundamental form works perfectly
well for compact polyhedra in $H^3$ --- and thus also for objects which
locally look like them. For ideal polyhedra it also works quite well; in
some cases (see below in part \ref{III-angles}) it is helpful to "glue"
a hemisphere in each of the length $2\pi$ circles corresponding to the
ideal vertices.

\blm \label{lm:necessary}
Let $S$ be a smooth (resp. polyhedral) locally convex surface in $H^3$,
and let $I$ and $\III$ be 
its induced metric and third fundamental forms respectively. Then:
\begin{enumerate}
\item $I$ has curvature $K\geq -1$ (resp. has curvature $-1$, except at
  the vertices, where the singular curvature is positive);
\item $\III$ has curvature $K\leq 1$ (resp. has curvature $1$, except at the
  dual vertices, where the singular curvature is negative);
\item the closed geodesics of $\III$ have length $L\geq 2\pi$.
\end{enumerate} 
\elm

In statements (1) and (2), the equality is attained only in the
degenerate cases, i.e. when the surface is not locally strictly
convex. In statement (3) it also corresponds to a very degenerate case,
as we will see below.

\bpv
The first and second point are consequences of the Gauss formula in the
smooth case and can be checked locally in the polyhedral cases. For the
last point the reader is refered to \cite{RH,CD,bao-bonahon,shu,cpt} for
different approaches of the polyhedral case, and e.g. to \cite{these}
for the smooth case.
\epv

\medskip

We will also prove here lemma \ref{lm:topolo}. The proof is based on a
number of propositions. In all this subsection we consider a cellulation
$\sigma$ of $\dr M$, along with a function $w$ on the edges of $\sigma$
verifying the hypothesis of theorem \ref{tm:dihedral}.

\bprop \label{pr:aretes}
No 2-cell of $\sigma$ can have two edges sent to the same segment in
$\dmt$. The intersection of two distinct 2-cells in $\dmt$ can not
contain more than one 1-cell. 
\eprop

\bprop \label{pr:sommet-1}
No 2-cell can have two vertices sent to the same point in $\dmt$. 
\eprop

\bprop \label{pr:sommet-arete}
The intersection of two 2-cells in $\dmt$ can not contain both a vertex
and an edge (which are disjoint). 
\eprop

\bprop \label{pr:sommet-2}
The intersection of two 2-cells in $\dmt$ can not contain two vertices. 
\eprop

The proof of lemma \ref{lm:topolo} clearly follows from those
propositions. 

\bpv[Proof of proposition \ref{pr:aretes}]
Suppose that two edges of a 2-cell $C$ are sent to the same segment in
$\dmt$. This segment would then constitute a circuit in the 1-skeleton
of $\sigma$, on which the sum of the values of $w$ is strictly less than
$\pi$; this would contradict the hypothesis of theorem \ref{tm:dihedral}.

Similarly, if two 2-cells $C$ and $C'$ in $\dmt$ have two edges in
common, then they constitute a circuit on which the sum of the values of
$w$ is strictly less than $2\pi$. 
\epv

\bpv[Proof of proposition \ref{pr:sommet-1}]
Let $C$ be a 2-cell having two vertices which are sent to the same point
$x_0$ in $\dmt$. Since the edges of $C$ are sent to distinct segments by the
proposition \ref{pr:aretes}, $C$ separates $\dmt$ in two parts, one
compact and the  other non-compact, whose closure intersect at $x_0$. 

Let $c$ be the elementary circuit made of the edges adjacent to
$x_0$. Then $c$ is the union of two sequences of edges, one, say $c_i$,
made of the edges contained in the closure of the compact domain of the
complement of $C$, and and other, say $c_o$, made of the other edges. 
Both sequence constitute a non-elementary circuit, and the sum of the
values of $w$ is less than $2\pi$ on both. This again contradicts the
hypothesis of theorem \ref{tm:dihedral}.
\epv

\bpv[Proof of proposition \ref{pr:sommet-arete}]
Suppose that two cells $C$ and $C'$ have in common both a vertex, $x_0$,
and an edge $e$. Let $c$ be the elementary circuit made of the edges
adjacent to $x_0$. Again, $c$ can be written as the union of $c_i$ and
$c_o$, where $c_i$ is the sequence of edges contained in the bounded
domain in the complement of $C\cup C'$. Adding $e$ to either $c$ or $c'$
leads to a non-elementary circuit. But, according to the hypothesis of
theorem \ref{tm:dihedral} the sum of the values of $w$ on
the edges of $c$ is $2\pi$, so that the sum of the values of $w$ has to
be less than $2\pi$ either on the edges of $c_i$ and $e$, or on the
edges of $c_o$ and $e$.
\epv

\bpv[Proof of proposition \ref{pr:sommet-2}]
Suppose now that $C$ and $C'$ share two vertices $x_0$ and $x_1$. Let
$c$ be the elementary circuit made of the edges adjacent to $x_0$, and
let $c'$ be the elementary circuit made of the edges adjacent to
$x_1$. Again, both $c$ and $c'$ decompose into two sequences of edges,
those in the closure of the bounded domain in the complement of $C\cup
C'$ (call them $c_i$ and $c'_i$, respectively) and the others (let them
be $c_o$ and $c'_o$, respectively). 

Then $c_i\cup c'_i$, $c_i\cup c'_o$, $c_o\cup c'_i$ and $c_o\cup c'_o$
each is a non-elementary circuit, while the sum of the values of $w$ on
at least one of them is clearly at most $2\pi$, again contradicting the
hypothesis of theorem \ref{tm:dihedral}. 
\epv

\subsection{Third fundamental form versus dihedral angles} 
\label{III-angles}

One of the points which should be clear from the introduction is that we
want to consider the manifolds with polyhedral boundary together with those
having smooth boundary. It is in that respect necessary to understand what
the relationship between the third fundamental form and the dihedral
angles is. This should clear up how question concerning the third
fundamental form --- as in theorem \ref{tm:hmcb-III} --- are related to
questions on the dihedral angles --- as in theorem \ref{tm:dihedral}. 
We will work out the relationship here; it is mostly
well-known, the most interesting case is that of ideal polyhedra.

First note that for manifolds with a boundary that looks locally like a
compact hyperbolic polyhedron one can consider
the universal cover of $M$, and then the surface $\St^*$ in $S^3_1$ which is
dual to the universal cover of a connected component of $\dr M$. $\St^*$
is a polyhedral surface --- locally like a convex space-like polyhedron in
$S^3_1$, so that it carries a spherical cone-metric with negative
singular curvature at the singularities. The dihedral angles are then
just the lengths of the dual edges.

This metric --- which we will still call the third fundamental form of
the boundary --- lifts to a $\cat(1)$ metric on the boundary of the universal
cover of $M$. Indeed, this splits into a local curvature condition ---
which is satisfied by the local convexity, because the singular
curvature is negative at each vertex --- and a global condition on the
length of the closed geodesics, which is also true here because of lemma
\ref{lm:necessary}.

A rather important point is that, while knowing the dual metric (and the
combinatorics of the polyhedral surface) determines the dihedral
angles, the converse is not true --- the dihedral angles
determine the length of the edges of the dual surface, but it does not
determine the shape of the dual faces with more than 4 edges. This is
already the source of interesting questions for convex polyhedra in
$H^3$; for instance, it is still an open problem to know whether a
convex hyperbolic polyhedron can be infinitesimally deformed without
changing its 
dihedral angles, see e.g. \cite{dap}, or \cite{stoker} for an analogous
(and also open) problem in the Euclidean case. 

For ideal polyhedra, there are two related ways of defining the third
fundamental form. If one considers the dual surface to the universal
cover of one of the components of the boundary, one obtains the induced
metric on the dual of a pleated surface, which is a tree. The third
fundamental form therefore reduces to the lengths of the edges of 
a graph. The length of the edges are the
(exterior) dihedral angles of the corresponding edges of the ideal
polyhedron. 

To each ideal vertex corresponds a face of the dual graph, with the sum
of the edge lengths equal to $2\pi$ (since the sum of the exterior dihedral
angles at an ideal vertex is $2\pi$). One can glue in each of those faces
a hemisphere (with its canonical metric). The result is a metric space
on $\dr M$, which obviously has negative singular curvature at its
singular points, because the singular points correspond to the vertices
of the graph, and the total angle around those points is $\pi$ times
the number of faces. We call $\III$ the corresponding metric on $\dr
M$. Note that $\III$ is the ``natural'' third fundamental form of $\dr
M$ for instance in a limit sense, as follows:

\bpt
Let $(\Omega_n)_{n\in \N}$ be an increasing sequence of open subsets of
$M$ with smooth, convex boundary, such that $\cup_n \Omega_n=M$. Then
the third fundamental forms of $\dr \Omega_n$ converge to $\III$.
\ept

We leave the proof to the reader.

The third fundamental form defined in this way has the important
properties below. The second strongly contrasts with the situation for
compact polyhedra.

\bpt \label{pt:314}
\begin{enumerate}
\item $\III$ lifts to a $\cat(1)$ metric on each
boundary component of the universal cover of $M$.
\item  There is a simple way to recover the  cellulation of $\dr M$ from
$\III$, and therefore also the dihedral angles.
\end{enumerate}
\ept

\bpv
The first point is again a consequence of lemma \ref{lm:necessary}.

For the second point, note that the dihedral angles of ideal polyhedra
are in $(0, \pi)$, and so are the edge lengths of the dual surfaces. 
Now any geodesic segment that enters a hemisphere can exit it only after
a path of length $\pi$. Since the dihedral angles are less than $\pi$,
the edges of the dual cellulation can not
enter the hemispheres. On the other hand all the segments in the boundary
of the hemispheres must be edges, and this recovers the dual cellulation,
and thus also the cellulation of $\dr M$.
\epv

\section{The fuchsian case}

We will investigate in this section some properties of manifolds with
polyhedral boundary in the fuchsian case; that is, we consider a metric
$g$ on $M$ such that $\partial M$ is polyhedral, with an isometric
involution $s$ of $(M, g)$ which fixes a compact surface. We also
suppose that $\dr M$ have two connected components, which are exchanged
by $s$. We refer the reader to \cite{leibon1,leibon2} for some recent
results on those manifolds. 

Another way to
consider such manifolds is to take the universal cover $S$ of one the two
connected components of $\dr M$; it is a convex surface in $H^3$, which
moreover is equivariant under the action of a surface group fixing a
totally geodesic 2-plane. 

We will first prove an existence and uniqueness result for the third
fundamental forms of such surfaces, in the case where they locally look
like a compact (rather than ideal) polyhedron in $H^3$. This is done using a
deformation argument as in the case of hyperbolic polyhedra (see
\cite{Al}), and is also similar to what can be done in the smooth case
\cite{iie}.

This first result will then be used to prove a similar existence result
for the dihedral angles of surfaces which locally look like
an ideal polyhedron; the main idea is to approximate this case by the
previous one, using a compactness result to prevent degeneracies from
occuring. The proof will also give a result concerning the manifolds with a
polyhedral boundary with some ideal and some "non-ideal" vertices.

The result which is obtained in this way --- or at least its ideal part
--- might look like a weak and
partial version of some of the results stated in the introduction. It is
proved in a very different way, however, and turns out to be necessary
for the more general cases, because it implies a technical statement ---
on the connectedness of some spaces of metrics on the boundary --- which
I do not know how to obtain directly. 

Mathias Rousset \cite{rousset1} has recently achieved anoter related
result, concerning the dihedral angles of fuchsian hyperideal polyhedra,
which includes the case of ideal fuchsian polyhedra. He reduces the
study of hyperideal polyhedra to that of finite polyhedra.

In all this section, we fix a surface $S$ of genus $g\geq 2$; we will be
interested in equivariant embeddings of $S$, or, in other terms, in
fuchsian hyperbolic metrics on $S\times [-1,1]$. 

\subsection{Infinitesimal rigidity of finite fuchsian polyhedra}

We first consider equivariant polyhedra which look locally like compact
polyhedra. First recall the definition of a polyhedral embedding.

\bdf \label{df:polyhedral}
A {\bf polyhedral embedding} of a surface $S$ into $H^3$ (resp. $\R^3_1$)
is a couple 
$(\sigma, \phi)$, where $\sigma$ is a cellulation of $S$ and $\phi$ is a
map from $S$ to $H^3$ (resp. $\R^3_1$) which:
\begin{itemize}
\item is injective.
\item sends each edge of $\sigma$ to a segment in $H^3$ (resp. $\R^3_1$).
\item sends each 2-face of $\sigma$ to the interior of a compact, 
convex polygon
in a totally geodesic 2-plane in $H^3$ (resp. $\R^3_1$).
\item is locally convex at each vertex.
\end{itemize}
\edf

The last condition is not always necessary; it is included here
since all the polyhedral objects that we will consider are convex.

\bdf
A {\bf finite equivariant polyhedron} is a couple $(\phi, \rho)$, where
$\phi$ is a polyhedral embedding of the universal cover of a surface $S$
into $H^3$ and $\rho$ is a group morphism from $\pi_1S$ into
$\isom(H^3)$, such that:
$$ \forall x\in \St, \forall \gamma\in \pi_1S, \phi(\gamma
x)=\rho(\gamma)\phi(x)~. $$ 
\edf

We are specially interested in the equivariant polyhedra which are the
boundary of the universal covers of the fuchsian hyperbolic manifolds
mentioned above.

\bdf \label{df:fuchsian}
A finite equivariant polyhedron $(\phi, \rho)$ is {\bf fuchsian} if
$\rho(\pi_1S)$ is contained in the subgroup of elements which leave
invariant a given plane $P\subset H^3$.
\edf

The definition of ideal equivariant polyhedra, and of semi-ideal
equivariant polyhedra, is similar. We first define the notion of
polyhedral map.

\bdf
An {\bf ideal polyhedral embedding} of a surface $S$ into $H^3$ is a couple
$(\sigma, \phi)$, where $\sigma$ is a cellulation of $S$ and $\phi$ is a
map from $S$ to $H^3\cup \dr_\infty H^3$ which:
\begin{itemize}
\item is injective.
\item sends each vertex of $\sigma$ to an ideal point (in $\dr_\infty
H^3$). 
\item sends each edge of $\sigma$ to a geodesic in $H^3$, which
connects the ideal points corresponding to the vertices.
\item sends each 2-face of $\sigma$ to the interior of an ideal  convex
polygon in a totally geodesic 2-plane in $H^3$.
\item is locally convex .
\end{itemize}
A {\bf semi-ideal polyhedral embedding} is defined likewise, except that
the vertices can be sent either to ideal points or to ``usual'' points
of $H^3$.  
\edf

\bdf
An {\bf ideal equivariant polyhedron} is a couple $(\phi, \rho)$ where
$\phi$ is an ideal polyhedral embedding of the universal cover of a
surface $S$ into $H^3$ and $\rho$ is a group morphism from $\pi_1S$ into
$\isom(H^3)$, such that:
$$ \forall x\in \St, \forall \gamma\in \pi_1S, \phi(\gamma
x)=\rho(\gamma)\phi(x)~. $$ 
A {\bf semi-ideal equivariant polyhedron} is defined in the same way,
but with $\phi$ a semi-ideal polyhedral embedding.
\edf

Note that the induced metric on equivariant (semi-)ideal polyhedra is defined
not only on the universal cover of the underlying surface, but also, by
equivariance, on the quotient surface. For ideal equivariant polyhedra
it is a hyperbolic metric with cusps corresponding to the vertices
(i.e. a hyperbolic metric of finite area on the quotient surface
$S$). For semi-ideal polyhedra, the induced metric is a cone-metric,
with singular points corresponding to the vertices which are not
ideal. It also has cusps corresponding to the ideal vertices. 

\bdf \label{df:ideal-fuchsian}
An ideal (resp. semi-ideal) equivariant polyhedron $(\phi, \rho)$ is
{\bf fuchsian} if $\rho(\pi_1S)$ is contained in the subgroup of
elements which leaves invariant a given plane $P\subset H^3$.
\edf

The first result we need to mention is an infinitesimal rigidity result
for finite equivariant polyhedral embeddings in Minkowski 3-space.

\blm \label{rig-minko}
Let $(\phi, \rho)$ be a convex equivariant space-like polyhedron in $\R^3_1$,
such that the representation $\rho$ fixes the origin. There is no
non-trivial infinitesimal deformation of $(\phi, \rho)$, preserving the
condition that $\rho$ fixes the origin, which does not change the
induced metric by $\phi$ to the first order.
\elm

This statement can be found as theorem B in Igor Iskhakov's thesis
\cite{iskhakov}, where it is proved by an extension of Cauchy's ideas on
the rigidity of polyhedra (see \cite{Cauchy,stoker}) to surfaces of
genus $g\geq 2$.  

I guess that an alternate proof could be given, along the approach given
in \cite{iie} for smooth surfaces; the key point would be to replace the
integration by part which works for the smooth case by a discrete
version.

There is yet another way to prove this infinitesimal rigidity result, as
well as the infinitesimal rigidity in \cite{iie}, and more general
results concerning for instance hyperideal fuchsian polyhedra. It is
based on the Pogorelov map, which is also used for instance in
\cite{hmcb}, to bring the problem in $\R^3$. There, the infinitesimal
rigidity can be proved using the elementary fact that, if an
infinitesimal deformation of a (smooth or polyhedral) convex surface is
isometric, then the graph of its coordinates are saddle surfaces. 

As a consequence, we find that the same result holds in the de Sitter
space:

\blm \label{rig-sitter}
Let $(\phi, \rho)$ be a convex equivariant space-like polyhedron in $S^3_1$,
such that the representation $\rho$ fixes a point $x_0$. There is no
non-trivial infinitesimal deformation of $(\phi, \rho)$, preserving the
condition that $\rho$ fixes $x_0$, which does not change the
induced metric by $\phi$ to the first order.
\elm

The proof of lemma \ref{rig-sitter} from lemma \ref{rig-minko} is
related to a remarkable trick invented by Pogorelov
\cite{Po}, which allows one to take an infinitesimal rigidity problem
from a space-form to another. The crucial point is that this can be done
in this setting, i.e. for equivariant objects when the representation
fixes a point, as shown in \cite{iie} for smooth surfaces. Moreover, the
polyhedral case works just like the smooth case, details on this stand
in \cite{dap,shu,cpt}.

Note that the same could be done also in the anti-de Sitter space
$H^3_1$, and this is indeed done in \cite{iie} for smooth surfaces. We
leave this point to the reader, however, since it will not be necessary
below.

Using the duality between $H^3$ and the de Sitter space $S^3_1$, we
immediately find a translation of this lemma in terms of fuchsian
surfaces in $H^3$:

\bcr \label{cr:fuchsian}
Let $(\phi, \rho)$ be a convex equivariant polyhedron in $H^3$,
such that the representation $\rho$ fixes a plane $\pi_0$. There is no
non-trivial infinitesimal deformation of $(\phi, \rho)$, preserving the
condition that $\rho$ fixes $\pi_0$, which does not change the
third fundamental form of the image of $\phi$ to the first order.
\ecr

\subsection{Compactness of fuchsian polyhedra}

The main technical tool of this subsection is a compactness result,
which is necessary to obtain the existence and uniqueness result for
fuchsian polyhedra explained above, and stated below as theorem
\ref{tm:fuchsian}. It is stated in a more general context, however, so
as to be used also later in this section, to prove a result
for ideal or semi-ideal fuchsian polyhedra. 

We now fix two integers $N\geq 1$ and $g\geq 2$; $g$ will be the genus
of the surface $S$ considered, and $N$ will be the
number of vertices of the 
polyhedral surfaces and the number of singular points of the metrics
considered. 

\bdf \label{df:poly}
We call:
\begin{itemize}
\item $\cP^C$ the set of finite fuchsian polyhedra of genus $g$ with $N$
  vertices in 
  $H^3$. 
\item $\cP^I$ the set of semi-ideal fuchsian polyhedra of genus $g$ with $N$
  vertices in $H^3$.
\end{itemize}
\edf

\bdf \label{df:third}
We call:
\begin{itemize}
\item $\cM^C$ the set of spherical cone-metrics on $S$ with $N$
  cone-points where the singular curvature is negative, and such that all
  contractible closed geodesics have length $L>2\pi$ (up to isotopy). 
\item $\cM^I$ the set of spherical cone-metrics with $N$ singular
  points, 
  where the singular curvature is negative, and such that contractible closed
  geodesics have length $L>2\pi$, except when they bound a hemisphere
  (again up to isotopy). 
\end{itemize}
\edf

It is clear (using lemma \ref{lm:necessary}) that the third fundamental
forms of the elements of $\cP^C$ are in $\cM^C$. 
Similarly, the third fundamental forms of the elements of
$\cP^I$ are in $\cM^I$, after one makes a simple surgery: gluing a
hemisphere on each of the circles, of length $2\pi$, which correspond to
the ideal vertices. 

\blm \label{lm:comp-fuchsian}
Let $(\phi_n, \rho_n)$ be a sequence of finite fuchsian polyhedral
embeddings in
$S^3_1$, with representations $\rho_n$ fixing a point $x_0$. Let
$(h_n)\in (\cM^C)^{\N}$ be the 
induced metrics. Suppose that $(h_n)$ converges, as $n\rightarrow
\infty$, to a metric $h\in \cM^I$. 
Then, after taking a subsequence and renormalizing, $(\phi_n, \rho_n)$
converges to 
a convex, fuchsian polyhedron $(\phi, \rho)$. Moreover $\rho$ fixes $x_0$. 
\elm

In this statement, the ``renormalization'' is by composition on the left
by an isometry. 

Note that $(\phi, \rho)$ might have some faces which are tangent to the
boundary at infinity of $H^3$ in the projective model of $H^3$ and
$S^3_1$. More precisely, this happens exactly when $h$ has closed
geodesics of length $2\pi$ bounding hemispheres, and the faces tangent
to the boundary at infinity are precisely those hemispheres. 

The proof of lemma \ref{lm:comp-fuchsian} depends on the following
propositions. The first is well known and its proof is easy, so we leave
it to the reader.

\bprop \label{pr:h}
Let $H_0$ be a space-like totally geodesic plane in $S^3_1$. The
function $v$ on $S^3_1$ defined as the hyperbolic sine of the oriented
distance to $H_0$ satisfies:
$$ \hess(v) = - v g_0~, $$
where $g_0$ is the metric of $S^3_1$. 
\eprop

This function will be used here --- as in other similar problems, see
e.g. \cite{these,iie} --- to control the lengths of closed geodesics of
surfaces or the total extrinsic curvature of the surface along those
geodesics.

The next proposition describes another way to write the canonical metric
of the de Sitter space. The proof is again elementary.

\bprop \label{pr:metrique}
The restriction to the future cone $C_+(x_0)$ of $x_0$ of the canonical
metric of $S^3_1$ can be written as:
$$ v^2 g_{H^2} - \frac{dv^2}{1+v^2}~, $$
where $g_{H^2}$ is the canonical metric on $H^2$.
\eprop

\bpv
For $t>0$, let $\Sigma_t$ be the set of points $x$ in the future cone of
$x_0$ such that there is a time-like segment going from $x_0$ to $x$ of
length $t$. Then a simple computation shows that $\Sigma_t$ has second
fundamental form $\II_t=\coth(t)I_t$, where $I_t$ is the induced metric
on $\Sigma_t$. Moreover the surfaces $\Sigma_t$ are equidistant, so that
an integration shows that the induced metrics are:
$$ I_t = \sinh^2(t) g_{H^2}~,$$
so that the metric on the future cone of $x_0$ can be written as:
$$ \sinh^2(t) g_{H^2}-dt^2~. $$
The proposition follows by setting $v=\sinh(t)$.
\epv

We then need a simple proposition about the solution of an ordinary
differential inequality; it will be applied below to the function $v$
restricted to geodesic segments in the metrics $h_n$.

\bprop \label{pr:edo}
For any $L>0$ and $u_0>0$, there exists $c>0$ such that if $u:[0,
L]\rightarrow \R_+$ is a function which:
\begin{itemize}
\item is smooth and satisfies $u''=-u$ except at $N$ singular points
$x_1, \cdots, x_N\in (0,L)$;
\item has a positive jump in its derivative at the $x_i$;
\item is bounded from below by some  constant $u_0>0$,
\end{itemize}
then:
$$ \int_0^L \sqrt{\frac{1}{u^2} + \frac{u'^2}{u^2(1+u^2)}} ds \leq
c~. $$
\eprop

\bpv
First note that:
\begin{eqnarray*}
\sqrt{\frac{1}{u^2} + \frac{u'^2}{u^2(1+u^2)}} & \leq &
\frac{1}{|u|} + \frac{|u'|}{|u|\sqrt{1+u^2}} \\ 
& \leq & \frac{1}{u_0} + \left|\left(\frac{1}{u}\right)'\right|~.
\end{eqnarray*}
It is therefore enough to prove that the total variation over $[0,L]$ of
$1/u$ is bounded from above by a constant.

Now let $0\leq x_1\leq \cdots \leq x_p\leq L$ be the sequence of local
minima of $u$, and let $0\leq y_1\leq \cdots \leq y_q\leq L$ be its
local maxima. The properties of $u$ clearly imply that there is indeed a
finite number of minima and maxima. Then, if $y_j$ is a local maximum
which immediately follows the local minimum $x_i$:
$$ \int_{x_i}^{y_j} \left|\left(\frac{1}{u}\right)'\right| ds =
\frac{1}{u(y_j)}-\frac{1}{u(x_i)}~. $$ 
But, since $y_j$ is a local maximum, and since $u'$ has a positive jump
at the singular points, $u'(y_j)=0$, so that, for all $s\in [x_i, y_j]$:
$$ u(y_j)\geq u(s)\geq u(y_j)\cos (y_j-s)~. $$
We now consider two cases:
\begin{enumerate}
\item $|y_j-x_i|\geq \pi/4$. Then:
$$ \left|\frac{1}{u(y_j)}-\frac{1}{u(x_i)}\right|\leq \frac{2}{u_0}~. $$
\item $|y_j-x_i|\leq \pi/4$. Then:
$$ \left|\frac{1}{u(y_j)}-\frac{1}{u(x_i)}\right|\leq
\frac{1}{u(y_j)}\left(\frac{1}{\cos(y_j-x_i)}-1\right) \leq
\frac{1}{u(y_j)} 4 (y_j-x_i)^2~. $$
\end{enumerate}
An elementary symmetry argument shows that the same estimates apply when
$x_j$ is a local minimum which immediately follows a local maximum
$y_i$. As a consequence, we find that:
$$ \int_0^L \left|\left(\frac{1}{u}\right)'\right| ds \leq
\frac{2}{u_0}\frac{4L}{\pi} + \frac{4L^2}{u_0}~, $$
and the proposition follows.
\epv

We will also need the following more geometric estimate.

\bprop \label{pr:cones}
For any $r>0$ and any integer $N>0$, there exists $c>0$ as follows. Let
$\sigma$ be a cellulation of the disk $D$ with at most $N$ vertices, and let 
$\phi$ be a polyhedral space-like embedding of $D$ in $S^3_1$ with
combinatorics given by $\sigma$, such
that:
\begin{itemize}
\item the boundary of $\phi(S)$ is convex for the metric induced by $\phi$.
\item the boundary is at distance at least $r$ from the center $x_0$ in
  the induced metric.
\item $\phi$ remains in the future cone $C_+(x_1)$ of a point $x_1$.
\end{itemize}
Then the absolute value of the distance between $\phi(x_0)$ and $x_1$ is at
least $c$.
\eprop

\bpv
Since $\phi(D)$ is space-like with convex boundary, it remains outside the
future cone of each of its points, in particular outside the future cone
of $x_0$. 
Since it also remains in $C_+(x_1)$, it is not difficult to check that,
if $x_1$ was too close to $x_0$, $\phi(D)$ would have to remain 
in an arbitrarily small neighborhood of a cone; this is not possible for
a polyhedral 
surface having a fixed number of singular points. 
\epv

We will also need below a basic result in the theory of hyperbolic
surfaces; see e.g. \cite{FLP} for a proof. Its content is that, to
prevent a sequence of hyperbolic metrics from degenerating, one only
needs to bound from above the lengths of a finite set of closed
geodesics. 

\blm \label{lm:surfaces}
Let $\Sigma$ be a surface of genus $g\geq 2$. There exists a finite
subset $E$ of $\pi_1\Sigma$ such that, for any $\epsilon>0$, the set of
hyperbolic metrics on $\Sigma$ such that the closed geodesics
corresponding to the elements of $E$ have length at most $1/\epsilon$ is
compact. 
\elm

We can now state a proposition showing that, with the hypothesis of
lemma \ref{lm:comp-fuchsian}, the representations of the equivariant
polyhedra $(\phi_n, \rho_n)$ do not diverge; the last part of the proof
will be to show that, under the "length $2\pi$" condition, isolated
vertices can not escape to infinity. 

\bprop \label{pr:hmin}
In the setting of lemma \ref{lm:comp-fuchsian}, the sequence of
representations $(\rho_n)$ converges (after one takes a subsequence).
\eprop

\bpv
First we fix an integer $n\in \N$, and consider an equivariant
polyhedral embedding $\phi_n:\St\rightarrow S^3_1$.
Let $x_n$ be a point of $\St$ where the minimum of $v$
is attained. Let $v_n:=v\circ \phi_n$. Proposition \ref{pr:cones} shows
that $v_n(x_n)$  is 
bounded from below by a strictly positive constant. Let $E\subset
\Gamma=\pi_1S$ be a finite generating set, on which more details will be
given below. 

Let $\gamma\in E$. There is a minimal geodesic segment $c_{n,\gamma}$ going
from $x_n$ to $\gamma x_n$ in $(S, h_n)$. Since $E$ is finite, and since
$h_n\rightarrow h$ as $n\rightarrow \infty$, the length $L(c_{n,\gamma})$ of
$c_{n,\gamma}$ is bounded from above by a constant $L_0$ for each $x\in E$
and each $n\in \N$.  

Since $\phi_n(\St)$ is space-like and equivariant under the action of $\Gamma$,
which fixes $x_0$, $\phi_n(\St)$ remains in the future cone $C_+(x_0)$ of
$x_0$. Calling 
$u_n$ the restriction of $v$ to $c_{n,\gamma}$, proposition  \ref{pr:h} shows
that $u_n$ satisfies $u_n''=-u_n$, except when $c_{n,\gamma}$ crosses an
edge of $\phi_n(\St)$, and then  $u_n'$ has a positive jump. 

Proposition \ref{pr:metrique} implies that $C_+(x_0)$ has a natural
submersion $\rho:C_+(x_0)\rightarrow H^2$, such that the restriction of
$\rho$ to each surface $\{ v=\mbox{const} \}$ in $C_+(x_0)$ is a
homothety. 
Let $s$ be the length element induced by $\phi_n$ on $S$, and $t$ be the
length element of the hyperbolic metric induced on $S$ by $\rho\circ
\phi_n$. 
Proposition \ref{pr:metrique} then indicates that:
$$ ds^2 = u_n^2 dt^2 -\frac{du^2}{1+u_n^2}~, $$
$$ dt^2 = \frac{ds^2}{u_n^2} + \frac{du_n^2}{u_n^2 (1+u_n^2)}~, $$
so that the length
of the image of $c_{n,\gamma}$ by $\rho$ is:
$$ L_s(c_{n,\gamma}) = \int_0^{L(c_{n,\gamma})}
\sqrt{\frac{1}{u_n^2}+\frac{u_n'^2}{u_n^2 (1+u_n^2)}} ds~. $$
According to propositions \ref{pr:cones}, $u_n$ is bounded from below by
a constant $u_0>0$; proposition \ref{pr:edo} therefore shows that, as
$n\rightarrow \infty$, the lengths of the $s(c_{n,\gamma})$ remain
bounded from above by a constant. 
 
Now for each $n\in \N$, $\Gamma$ acts on $S^3_1$ fixing $x_0$, and
therefore $\Gamma$ has an action on $C_+(x_0)$ which leaves globally
invariant all the surfaces $\{ v=\mbox{const}\}$. Thus it acts 
by isometries on $H^2$ through $\rho$. The previous argument shows that the
translation distance of each element of $E$ remains bounded from above
as $n\rightarrow 
\infty$. Lemma \ref{lm:surfaces} then implies that the sequence $\rho_n$
remains in a compact subset of Teichm{\"u}ller space, and therefore that one
of its subsequences converges.
\epv

To prove that the sequence of equivariant polyhedra $(\phi_n, \rho_n)$
actually converges --- and not only the representations --- it is helpful
to consider a projective model of the part of the de Sitter space
$S^3_1$ which stands on one side of a totally geodesic space-like
plane. One can be constructed as follows. Remember that $S^3_1$ is 
isometric to a quadric in Minkowski 4-space with the induced metric:
$$ S^3_1 \simeq \{ x\in \R^4_1 | \langle x,x\rangle =1\}~. $$
Let $P_0$ be the affine hyperplane of equation $x_0=1$ in $\R^4_1$, and
let:
$$ S^3_{1,+} := \{ x\in S^3_1 | x_0>0\}~. $$
There is a natural map from $S^3_{1,+}$ to $P_0$ sending a point $x\in
S^3_{1,+}$ to the intersection with $P_0$ of the line going through $0$
and $x$. By construction it is projective, i.e. it sends geodesics to
geodesics; indeed, an elementary argument using the action of
$\mbox{SO}(3,1)$ shows that the geodesics of $S^3_{1,+}$ are the
intersections with $S^3_{1,+}$ of the 2-planes of $\R^4_1$ containing
$0$. Note that, since the $\phi_n(\St)$ remain in the future cone of a 
point, the projective model can be chosen so that a compact subset of
$\R^3$ contains the image of the surfaces $\phi_n(\St)$ for each $n$.

\bprop \label{pr:comp-proj}
After taking a subsequence and renormalizing, the sequence $(\phi_n(\St))$
converges in the 
projective model described above to an equivariant polyhedron.
\eprop

\bpv
This is a direct consequence of the convergence (after taking a
subsequence) of the representations, as stated in proposition
\ref{pr:hmin}. Indeed, we can renormalize the sequence so that, for a
given vertex $x\in \St$, $\phi_n(x)$ is constant. Then an elementary
compactness argument shows that, after taking a subsequence, the
vertices adjacent to $x$ also have converging images. Going to the
vertices adjacent to those and applying the same compactness argument
shows that they also have converging images, and this can be done until
all vertices in a fundamental domain of $\St$ have converging
image. Proposition \ref{pr:hmin} then implies the result.
\epv

We now have to exclude some cases, corresponding to a limit metric $h$
which has a closed geodesic of length $2\pi$. A similar assertion was
used --- and stated in a more general setting, including higher
dimensions --- in \cite{shu,cpt}.

\bprop \label{pr:2pi}
Suppose there is a finite set of vertices which converge to the same point
$x_\infty$ in $\dr_\infty H^3$, while the other vertices do not. 
Then the limit metric $h$ has a closed geodesic of length $2\pi$ which
does not bound a hemisphere containing no vertex. 
\eprop

\bpv
We use proposition \ref{pr:comp-proj} and suppose that, in the
projective model described above, the sequence $(\phi_n(\St))$ converges to a
limit $P$. Then $x_\infty\in P$, and, by convexity and the fact that the
polyhedra $\phi_n(\St)$ remain in the exterior of the ball $B^3$
corresponding to 
$H^3$, $P$ also contains a neighborhood of $x_\infty$ in the plane $\pi_0$
tangent to $B^3$ at $x_\infty$. More precisely, since $P$ is a polyhedron,
one of its faces is the interior of a convex polygon $Q$ in $\pi_0$.

We will show that the boundary polygon, $\dr Q$, has length
$2\pi$ and is a geodesic for the limit induced metric $h$. The fact that
it has length $2\pi$ is an elementary fact of Lorentz geometry, since it
lies in a degenerate plane in $S^3_1$. To show that it is a geodesic of
$h$, we have to show that both sides of $P$ are concave for $h$.

Consider first the interior of the polygon $Q$. Since it carries a
degenerate metric, it is not difficult to realize that the metrics
induced on the corresponding faces of the $\phi_n(\St)$ converge to the metric
of a hemisphere. Therefore, this part of $P$ is concave for $h$.

For the other side, simply note that, for each $n$, the metric $h_n$ has
a point of negative singular curvature at the vertices of $\dr Q$,
i.e. the limit total angle at those vertices is at least $2\pi$. Since
the limit metric is isometric to a hemisphere in the interior of $Q$, it
means that the limit total angle at each vertex of $\dr Q$ of the
complement of $Q$ in $P$ is at least $\pi$, i.e. that the complement of
$Q$ in $P$ is concave for $h$, as needed.
\epv

\bpn{of lemma \ref{lm:comp-fuchsian}}
Proposition \ref{pr:hmin} shows that the sequence of representations
$(\rho_n)$ converges (after taking a subsequence), while proposition
\ref{pr:2pi} indicates that no vertex can escape to infinity. Therefore
$(\phi_n, \rho_n)$ converges. 
\epn

\subsection{Induced metrics on finite polyhedra}

We can now consider the maps $\Phi^C:\cP^C\rightarrow \cM^C$ and
$\Phi^I:\cP^I\rightarrow \cM^I$ sending a fuchsian polyhedron to its
third fundamental form. We will prove that $\Phi^C$ is a homeomorphism,
and that $\Phi^I$ is bijective. To show that for $\Phi^C$, we will apply a
deformation method, which I believe was invented by Aleksandrov
\cite{Al} to study the induced metrics on hyperbolic polyhedra, although
it was later used in many very different fields. 

First, it is clear that choosing an element of $\cP^C$ is equivalent to
choosing a hyperbolic metric on $S$ (or equivalently a fuchsian action
of $\Gamma:=\pi_1(S)$ on $H^3$) along with $N$ points $x_1, \cdots, x_N$
in $H^3/\Gamma$, under some conditions, i.e. that the $x_i$ all lie on
the boundary of their convex hull. Therefore, $\cP^C$ is a connected
manifold with boundary, of dimension $6g-6+3N$.

It is also clear, using the results of Troyanov \cite{troyanov}, that
choosing an 
element of $\cM^C$ is the same as choosing an element of the Teichm{\"u}ller
space of $S$ with $N$ marked points, along with the singular curvature
at each of the marked points. Thus $\cM^C$ is also a manifolds with
boundary of dimension $6g-6+3N$. 

Moreover, corollary \ref{cr:fuchsian} shows that $\Phi^C$ is locally
injective --- 
and therefore a local diffeomorphism --- between $\cP^C$ and $\cM^C$. To
prove that it is a homeomorphism, we need to prove that $\Phi^C$ is
proper (this is a consequence of lemma \ref{lm:comp-fuchsian}), 
that $\cP^C$ is connected, and that $\cM^C$ is simply connected.

The actual proof below is slightly more complicated than the outline
here; since I do not know how to prove directly the connectedness of
$\cM^C$, we will use a trick invented in \cite{Ri,RH}: one uses the
connectedness of a space of smooth metrics --- which is easy to prove ---
to check that two metrics $g_0$ and $g_1$ with $N$ singular points can
be connected by a 
path of metrics with at most $N'$ singularities, where $N'$ is a (large)
integer depending on $g_0$ and $g_1$.

\blm \label{lm:connected}
\begin{enumerate}
\item 
Let $g_0$ and $g_1$ be elements of the space  $\cM(S, N)$ of
spherical cone-metrics on $S$, with at most
$N$ singular points where the singular curvature is negative, and
contractible closed
geodesics of length $L>2\pi$. There exists an integer $N'$ depending on
$g_0$ and $g_1$ such that $g_0$ and $g_1$ can be connected in
$\cM(S, N')$. 
\item Let $c:S^1\rightarrow \cM(S,N)$. There exists $N'\geq N$ and
  a disk $D\subset \cM(S, N')$ with $\dr D=c(S^1)$.
\end{enumerate}
\elm

\bpvs
We do not give a full proof, since it is almost the same as the one
given by Rivin and Hodgson \cite{RH}; the only difference is that the
surfaces considered here have genus $g\geq 2$, rather than $0$, but this
does not appear in the proof.

The main point is that the space of smooth metrics on $S$ with
curvature $K<1$ and contractible closed geodesics of length $L>2\pi$ is
connected. Indeed, given two such metrics $h_0$ and $h_1$, one can take
any path connecting them in the space of Riemannian metrics on $S$,
and then "scale up" the middle part, to make sure that the curvature
remains small and the closed geodesics remain large (this is also
explained in \cite{cras,these}). 

To prove the first part of the
lemma, one first shows that $g_0$ and $g_1$ can be
approximated by smooth metrics $h_0$ and $h_1$ satisfying the curvature
and geodesic length conditions. One thus obtains a path $(h_t)_{t\in
  [0,1]}$ of metrics satisfying the same conditions. One then proves
that there is an integer $N'$ such that the metrics in $(h_t)_{t\in
  [0,1]}$ can be approximated by polyhedral metrics in $\cM(S, N)$
in a continuous way.

The second point can be proved in an analogous way, using the fact that
the space of smooth metrics on $S$ with curvature $K<1$ and contractible
closed geodesics of length $L>2\pi$ is simply connected.
\epvs

\btm \label{tm:fuchsian}
Let $S$ be a surface of genus $g\geq 2$, and let $h$ be a
spherical cone-metric on $S$, with negative singular curvature at
the singular points. Suppose that all contractible closed geodesics of
$(S, h)$ have length $L>2\pi$. 
Then there is a unique fuchsian polyhedral
embedding of $(S, h)$ into $H^3$ whose third
fundamental form is $h$.
\etm

Note that the uniqueness here is of course up to global isometries of
$H^3$. Another remark is that the length condition is necessary by lemma
\ref{lm:necessary}, because the curvature conditions at the vertices
imply that the image is convex. By the way, it would be interesting to
know whether a similar result also holds with $S^3_1$ replaced by the
anti-de Sitter space.  

\bpv
As mentioned above, we already know that $\Phi^C$ is a local
diffeomorphism. Lemma \ref{lm:comp-fuchsian} shows that $\Phi^C$ is proper,
so it is a covering of the connected components of $\cM^C$ which
intersect its image. The first part of lemma \ref{lm:connected} shows
that all of $\cM^C$ is in the image, while the second part, along with
the fact that $\cP^C$ is connected, indicates that each point of
$\cM$ has a unique inverse image.
\epv

\subsection{Fuchsian ideal manifolds}

We have already mentioned above that the dihedral angles of ideal
polyhedra --- and of ideal fuchsian polyhedral embeddings, etc --- is the
analog of 
the third fundamental form of finite polyhedra (and of finite fuchsian
polyhedral embeddings, etc). Theorem \ref{tm:fuchsian} should therefore
have an analog for ideal fuchsian polyhedral
embeddings in terms of dihedral angles. 
We will prove first the existence part of this statement; the
more general result concerning ideal fuchsian polyhedra will be a
consequence of other 
results proved here (specifically, of theorem \ref{tm:dihedral}) but the
existence result given here will be necessary to prove the more general
statement. Moreover, we will consider here the case of semi-ideal
polyhedral embeddings, which is not covered by theorem \ref{tm:dihedral}.

\blm \label{lm:fuchsian-ideal}
Let $\Gamma$ be the 1-skeleton of a cellulation of a surface $S$ of
genus $g\geq 
2$. Let $w:\Gamma_1\rightarrow (0,\pi)$ be a function on the set
$\Gamma_1$ of edges of $\Gamma$ such that:
\begin{enumerate}
\item for each elementary circuit in $\Gamma$, the sum of the values of
$w$ is equal to $2\pi$;
\item for each non-elementary circuit, the sum of the values of $w$ is
strictly above $2\pi$. 
\end{enumerate}
Then there is a fuchsian ideal 
embedding of $(S, h)$ into $H^3$ whose combinatorics is given by
$\Gamma$, with exterior dihedral angles given by $w$.
\elm

The key point for the remainder of the paper is that, as a consequence,
the space of possible dihedral angle assignations is connected ---
because the space of possible ideal fuchsian polyhedral embeddings is
connected. 

This lemma is actually a consequence, using property \ref{pt:314},
of the more general statement
below, so we don't prove it separately.

\blm \label{lm:semi-ideal}
Let $S$ be a surface of genus $g\geq 2$, and let $h$ be a
spherical cone-metric on $S$, with negative singular curvature at
the singular points. Suppose that all contractible closed geodesics of
$(S, h)$ have length $L>2\pi$, except when they bound a hemisphere. 
Then there is a fuchsian polyhedral
embedding of $(S, h)$ into $H^3$ whose third
fundamental form is $h$.
\elm

\bpv
We choose a sequence of spherical cone-metrics $(h_n)_{n\in \N}$ such
that:
\begin{enumerate}
\item $h_n$ converges to $h$.
\item for each $n$, $h_n$ is a spherical cone-metric on $S$, with
  negative singular curvature at the singularities. 
\item for each $n$, the contractible closed geodesics of $(S, h_n)$ have
  length strictly above $2\pi$.
\end{enumerate}
It is not difficult to find such an approximating sequence; one has to
decrease slightly the length of some of the edges, and to replace the
hemispheres by interiors of convex polygons in $S^2$. 

Then apply theorem \ref{tm:fuchsian} to obtain, for each $n$, a fuchsian
finite polyhedral embedding of $S$ in $S^3_1$. Finally, lemma
\ref{lm:comp-fuchsian} shows that, after renormalizing this sequence and
taking a subsequence, it converges to a semi-ideal fuchsian polyhedral
embedding inducing the metric $h$ on $S$.
\epv

Mathias Rousset \cite{rousset1} recently remarked that the uniqueness
part of this statement can be obtained in a rather straightforward way,
by using the infinitesimal Pogorelov map to show an infinitesimal
rigidity result for semi-ideal fuchsian polyhedra; the uniqueness for
semi-ideal polyhedra then follows from studying the map sending a
fuchsian polyhedron to its third fundamental form in the neighborhood of
semi-ideal polyhedra. He thus obtained that:

\btm[Rousset \cite{rousset1}] \label{tm:semi-ideal}
Let $S$ be a surface of genus $g\geq 2$, and let $h$ be a
spherical cone-metric on $S$, with negative singular curvature at
the singular points. Suppose that all contractible closed geodesics of
$(S, h)$ have length $L>2\pi$, except when they bound a hemisphere. 
Then there is a unique fuchsian polyhedral
embedding of $(S, h)$ into $H^3$ whose third
fundamental form is $h$.
\etm

\section{Triangulations}

This section deals with questions concerning triangulations of a given
ideal hyperbolic manifold. Although the existence of a triangulation
inducing a given triangulation of the boundary might appear natural at
first sight, it is not easy to prove --- at least this is not proved
here. This is similar to the situation concerning finite volume
hyperbolic manifolds, where an ideal triangulation would be helpful but
is not known to exist in general; see \cite{porti-petronio}. We will
only prove that any \ihm $M$ is ``almost triangulable'' in the sense
that it has a finite cover $\Mb$ which does admit an ideal
triangulation.

Note that alternative approaches could perhaps be followed. One is based
on the fact that the main property of the ideal simplices which is used
here, namely that the volume is a concave function, remains valid for
more general ideal polyhedra (this is a result of \cite{Ri2}). Another
(which was pointed out by Francis 
Bonahon) uses a triangulation which might include some degenerate
simplices. Lemma \ref{lm-vol} does not apply for those simplices, since
the volume function is concave but not strictly concave in those
cases. It might however be possible to prove that the sum of the volumes
of the simplices in a triangulation remains strictly concave, which is
basically what one needs in the next sections.

Although those other approaches could presumably lead to shorter proofs,
we chose the method described here because, once some technical points
are proved, it gives rather simple picture of what goes on; and also
because the method it contains might be useful in other settings. 

\pg{Triangulations, cellulations}
We first give more details about what we call a triangulation here.

\bdf \label{df-cell}
Let $M$ be an \ihm. A {\bf cellulation} $C$ of $M$ is a finite family
$C_1, C_2, \cdots, C_n$ of non-degenerate, closed, convex, ideal polyhedra
isometrically embedded in $M$, such that:
\begin{enumerate}
\item for $i\neq j$, the interiors of $C_i$ and $C_j$ are disjoint;
\item the union of the $C_i$ is all of $M$;
\item for $i\neq j$, if $C_i\cap C_j\neq \emptyset$, then it is a face
of both $C_i$ and $C_j$.
\end{enumerate}
The $C_i$ are the ``cells'' of the cellulation $C$. 
\edf

\begin{figure}[h]
\centerline{\psfig{figure=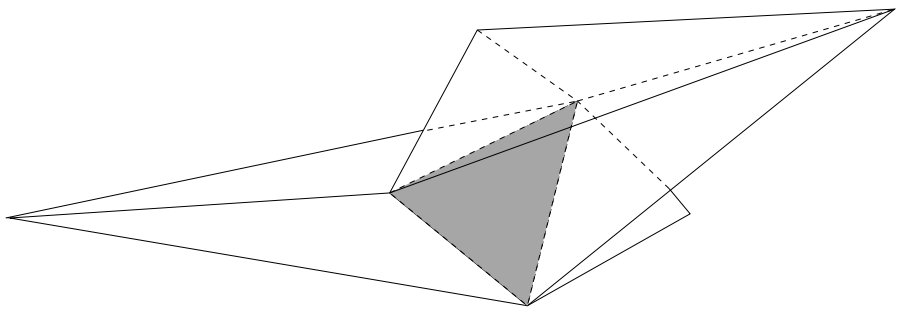,height=4cm}}
\caption{This is forbidden in our cellulations}\label{fig:faces}
\end{figure}

The third condition excludes some ``bad'' configuration, like the one
depicted in figure \ref{fig:faces}, where two simplices have an
intersection which is not a face in any of them.

\bdf \label{df-tri}
A triangulation of $M$ is a cellulation whose cells are all simplices.
\edf

\bdf \label{df-atri}
An \ihm $M$ is {\bf triangulable} if it admits a triangulation. It is
{\bf almost triangulable} if it has a finite cover which is triangulable.
\edf

The main goal of this section will be to prove the:

\blm \label{lm-atri}
Any \ihm is almost triangulable.
\elm

The proof will proceed in several steps. The first point is the:

\bprop \label{pr-cell}
Any \ihm admits a cellulation.
\eprop

\bpv
It is done along the ideas of Epstein and Penner \cite{epstein-penner};
the situation here is simpler since the action of $\pi_1M$ on
$S^2\setminus \Lambda$ is discrete.

Let $M$ be an ideal hyperbolic manifold. Then $M$ is isometric to the
convex hull of a set $\{ x_1, \cdots, x_N\}$ of points in $\dr_\infty
E(M)$, where $E(M)$ is the unique convex co-compact hyperbolic manifold
in which $M$ admits an isometric embedding which is surjective on the
$\pi_1$.  

For each $i\in \{ 1,2,\cdots,N\}$, choose a ``small'' horoball
$b_i\subset E(M)$ with ideal point $x_i$; we suppose that the $b_i$ are
small enough to be disjoint. Let $B_i$ be the lift of $b_i$
to $H^3$, which is $\pi_1M$-invariant collection of disjoint horoballs
in $H^3$. 

Now we want to use this action of $\pi_1M$ on the horoballs to produce a
cellulation of $M$. This is done in \cite{epstein-penner} by considering
the action of $\pi_1M$ on the light cone in Minkowski 4-space, which
contains $H^3$ as a quadric. We will use here a similar, slightly more
complicated but
maybe a little more geometric, approach. It is based more explicitely on
the action of $\pi_1M$ on the space of horospheres, with an explicit
model from \cite{horo}.

We consider the projective model of $H^4$ and half of $S^4_1$, the de
Sitter space of dimension 4, corresponding to the 3-dimensional models
already described above. It can be obtained as follows. $H^4$ and
$S^4_1$ are both isometric to quadrics in Minkowski 5-space, with the
induced metric:
$$ H^4 \simeq \{ x\in \R^5_1 | \langle x,x\rangle = -1 \wedge x_0>0 \}~, $$
$$ S^4_1 \simeq \{ x\in \R^5_1 | \langle x,x\rangle = 1\}~. $$
Let $P_0$ be the affine hyperplane of equation $x_0=1$ in $\R^5_1$;
consider the map sending a point $x\in H^4$ (resp. $x\in S^4_1$ with
$x_0>0$) to the intersection with $P_0$ of the line going through $x$
and $0$. It is not difficult to check that $\phi$ is projective; it maps
$H^4$ to the interior of the radius one ball, and the part of $S^4_1$
standing on one side of a space-like hyperplane to its
exterior. We now consider this model only, with $P_0$ identified with
$\R^4$. 

Using the classical Poincar{\'e} model of $H^3$, we can map conformally
$H^3$ to the interior of a geodesic ball $B_0$ in $S^3$, e.g. to the
hemisphere:
$$ S^3_+:=\{ x\in \R^4 | \langle x,x\rangle =1 \wedge x_1\geq 0 \}~. $$ 
Horospheres in $H^3$ are then mapped to spheres in $S^3$ which are
interior to $S^3_+$ and tangent to its boundary. 
Those spheres are the boundaries at
infinity of the totally geodesic 3-planes in $H^4$ which are asymptotic
to a given 3-plane $H_1$, with $\dr_\infty H_1=\dr S^3_+$. Their dual
points in $S^4_1$ (using the hyperbolic-de Sitter duality,
see subsection \ref{ss34}) form the vertical cylinder $C_0^+$ 
which is tangent to $S^3$ along $\dr S^3_+$. In de Sitter terms, $C_0^+$
is the future (or past, depending of the orientation) light cone of a
point $H_1^*$ which is at infinity in the projective model of $S^4_1$
which we use. 

\begin{figure}[h]
\centerline{\psfig{figure=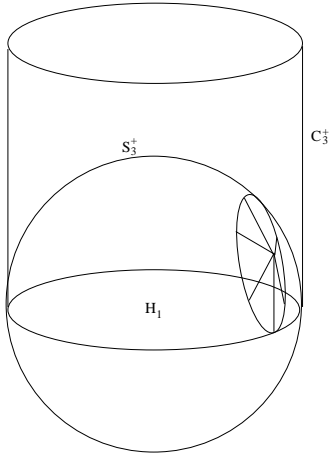,height=6cm}}
\caption{The space of horospheres as a cone}\label{fig:cone}
\end{figure}

By construction, the action of $\pi_1M$ on $H^3$ extends to a conformal 
action on $S^3$, and thus to a projective action on $\R^4$. This action
leaves invariant $S^3_+$, and thus also $H_1$ and $H_1^*$, and therefore
also $C_0^+$. For each $i\in \{ 1,\cdots, N\}$, the horoballs in $B_i$
corresponds to the points of an orbit $O_i$ of the action of $\pi_1M$ on
$C_0^+$. Since the horoballs in $B_i$ are disjoint, it is easy to see
that $O_i$ is discrete, with no accumulation point outside $\dr S^3_+$. 

To finish the proof, we proceed almost
as in \cite{epstein-penner}; consider the
convex hull (in $\R^4$) of $O_1\cup \cdots\cup O_N\cup H_1$, which by the
property just pointed out is locally finite outside $\dr S^3_+=\dr H_1$,
with faces which are polyhedra with a finite number of edges and
vertices. Then take the ``projection'' of this polyhedral surface on
$H^3\subset S^3$ in the vertical direction to obtain the required
cellulation. By construction it is invariant under the action of
$\pi_1M$ on $H^3$.

The cellulation of $M$ obtained in this way has a finite number of
cells. Otherwise, there would exist an edge $e$ meeting an infinite
number of fundamental domains of the action of $\pi_1M$ on $\Mt$. So $e$
would connect a vertex $v$ of $\Mt$ to the limit set $\Lambda$ of the
action of $\pi_1M$ on $H^3$. Now note that, by construction, $M$ is
covered by a finite set of (non-disjoint) ideal simplices $T_1, \cdots,
T_N$. So there would exist $i\in \{ 1, \cdots, N\}$ such that $e$
intersects $\gamma T_i$ for an infinite set $S$ of elements $\gamma$ of
$\pi_1M$. Going back to the projective model used above, for each
$\gamma\in S$, the segment $\gamma e$ should intersect $T_i$ --- and it
would also be
in the boundary of the convex hull constructed above. 

Now $\gamma e$ goes from $\gamma v$ to a point $\gamma v'$ of $\Lambda$; as
$\gamma\rightarrow \infty$, the vertical coordinates of both $\gamma v$
and $\gamma v'$ go to zero, so that the whole segment $\gamma e$ goes to
the horizontal hyperplane containing $0$. Therefore, $\gamma e$ lies
"below" $T_i$, so that the intersection $\gamma e\cap T_i$ can not be in
the boundary of the convex hull. So such an $e$ can not exist, and the
cellulation obtained has a finite number of cells. 
\epv

\subsection{From a cellulation to a triangulation}

We now consider a cellulation $C$ of an \ihm $M$, and call $F_0, F_1,
F_2$ and $F_3$ the sets of its faces of dimension $0,1,2$ and $3$
respectively. We also consider the universal cover $\Mt$ of $M$; $C$
lifts to a cellulation $\Ct$ of $\Mt$, and we call $\Ft_0, \Ft_1, \Ft_2$
and $\Ft_3$ the sets of its faces of the 4 possible dimensions.

From now on, and until section 8, we will consider a fixed triangulation
$\Sigma$ on $M$.  

\bdf \label{df-large}
Let $C$ be a cellulation of an \ihm $M$. $C$ is {\bf large} if, for any
cell $\sigma\in \Ft_3$, any vertex $s$ of $\sigma$ and any $\gamma\in
\pi_1(M)$, if $\gamma s$ is again a vertex of $\sigma$, then
$\gamma=1$. 
\edf

\bprop \label{pr-cover}
For any \ihm $M$ and any cellulation $C$ of $M$, $M$ has a finite cover
$\Mb$ such that $C$ lifts to a large cellulation.
\eprop

\bdf \label{df-polar}
Let $C$ be a cellulation of an \ihm $M$. A {\bf polarization} of $C$ is
a map $:\Ft_3\rightarrow \Ft_0$ which is equivariant under the action of
$\pi_1M$ and such that:
\begin{enumerate}
\item for any $\sigma\in \Ft_3$, $\rho(\sigma)$ is a vertex of $\sigma$;
\item if $\sigma\in \Ft_3$ and $\sigma'\in \Ft_3$ are adjacent
(i.e. $\sigma \cap \sigma'\in \Ft_2$), and if $\rho(\sigma)$ and
$\rho(\sigma')$ are both vertices of $\sigma\cap \sigma'$, then
$\rho(\sigma)=\rho(\sigma')$. 
\end{enumerate}
\edf

Let's pause to remark that, although a $(\pi_1M)$-equivariant map from
$\Ft_3$ to $\Ft_0$ induces a map from $F_3$ to $F_0$, the equivariant
map contains much more information. Indeed this is already apparent in
the simple case where one considers the manifold $S^1$, triangulated with
only one edge and one vertex. Its universal cover is $\R$, triangulated
with vertices at the integers. 
There are two $\Z$-equivariant maps sending a
segment $[k,k+1]$ to one of its endpoints: $\rho_1:[k,k+1]\mapsto k$,
and $\rho_2:[k,k+1]\mapsto k+1$. There is however only one map from
$F_3$ to $F_0$, since $F_0$ has only one element.

\bprop \label{pr-polar}
Any large cellulation $C$ of an \ihm $M$ admits a polarization.
\eprop

\bprop \label{pr-triang}
Any polarized cellulation of an \ihm $M$ can be subdivided to obtain a
triangulation. If the cellulation is large, so is the triangulation
obtained. 
\eprop

The proof of lemma \ref{lm-atri} clearly follows from the three
propositions above so it remains only to prove them. 

\bpn{of proposition \ref{pr-cover}}
Let $\sigma\in \Ft_3$. First note that, if $\gamma \in \pi_1M$ fixes the
vertices of $\sigma$, then $\gamma=e$; indeed, no non-identity element
of $\isom(H^3)$ has more than $2$ fixed points at infinity. Therefore,
the elements of $\pi_1M$ leaving $\sigma$ invariant are determined by
their actions on its vertices, so that there is a finite subset
$E_\sigma$ of $\pi_1M$ of elements leaving $\sigma$ invariant. 

Now $C$ has a finite number of cells, which we can call $\sigma_1,
\cdots, \sigma_N$. For each $i$, let $\sigmab_i$ be a cell of $\Ct$
whose projection on $M$ is $\sigma_i$. Then 
$E:=\cup_{i}E_{\sigmab_i}$ is finite, where $E$ is the set of elements
of $\pi_1M$ 
leaving one of the $\sigmab_i$ invariant. Since there is a finite number of
cells $\sigmab'$ sharing a vertex with some given cell $\sigmab$, the set
$F_{\sigmab}$ of elements $\gamma$ of $\pi_1M$ such that $\sigmab$ shares a
vertex with $\gamma \sigmab$ is finite. Therefore, the set
$F:=\cup_{i} F_{\sigmab_i}$ is finite. 

$\pi_1M$ is linear, and therefore residually finite (see
e.g. \cite{lyndon-schupp}, chapter III, 7.11). Thus
there exists a normal subgroup $\Gamma$ of $\pi_1M$ of finite
index, such that $\Gamma\cap F=\{ 1\}$. The corresponding finite cover
$\Mb$ of $M$ has the required property. 
\epn

\bpn{of proposition \ref{pr-polar}}
We will construct the required polarization $\rho$ as the endpoint of a
sequence of partially defined equivariant functions
$\rho_i:\Ft_3\rightarrow \Ft_0$ (that is, functions defined on a subset
of $\Ft_3$ only).

First choose $s_0\in \Ft_0$; since $C$ is large, for any $\sigma\in
\Ft_3$, at most one of the vertices of $\sigma$ is in
$(\pi_1M)s_0\subset \Ft_0$. Define $\rho_0$ on a cell $\sigma\in \Ft_3$
as follows: 
\begin{itemize}
\item if there exists $\gamma \in \pi_1M$ such that $\gamma s_0$ is a
vertex of $\sigma$, then set $\rho_0(\sigma):=\gamma s_0$;
\item otherwise, leave $\rho_0$ undefined at $\sigma$.
\end{itemize}
It is clear that this partially defined map is equivariant.

Now choose $s_1\in \Ft_0$ such that some $\sigma \in \Ft_3$ on which
$\rho_0$ is not defined has $s_1$ as a vertex, and define $\rho_1$ as
follows: 
\begin{itemize}
\item if there exists $\gamma \in \pi_1M$ such that $\gamma s_1$ is a
vertex of $\sigma$, then $\rho_1(\sigma):=\gamma s_1$;
\item otherwise, $\rho_1(\sigma):=\rho_0(\sigma)$.
\end{itemize}
The second case includes the possibility that $\rho_0$ is undefined at
$\sigma$, then $\rho_1$ remains undefined at $\sigma$.

Then repeat this construction with $s_2$ to obtain a map $\rho_2$,
etc. The number of cells of $F_3$ on which $\rho_i$ is not defined
decreases by at least one unit at each step, and $C$ has a finite number
of cells, so after a finite number of steps we obtain an equivariant map
$\rho:=\rho_N:\Ft_3\rightarrow \Ft_0$ which is everywhere defined.

We now want to prove that $\rho$ is a polarization. It is clear by
construction that, for any $\sigma\in \Ft_3$, $\rho(\sigma)$ is a vertex
of $\sigma$. Let $\sigma'\in \Ft_3$ be another cell, such that
$\sigma\cap\sigma'\in \Ft_2$. Let:
$$ i_0 := \max \{ i\in \{ 1,\cdots ,N\} \; | \; \exists \gamma\in
\pi_1M, \gamma s_i \; \mbox{is a vertex of} \; \sigma\}~, $$
$$ j_0 := \max \{ j\in \{ 1,\cdots ,N\} \; | \; \exists \gamma'\in
\pi_1M, \gamma' s_j \; \mbox{is a vertex of} \; \sigma'\}~. $$
We consider two cases:
\begin{enumerate}
\item $\gamma s_{i_0}$ is a vertex of $\sigma$ but not of $\sigma'$;
then $\rho(\sigma)=\gamma s_{i_0}$ is not a vertex of $\sigma'$, and
condition (2) of definition \ref{df-polar} is satisfied. The same
applies if $\gamma's_{j_0}$ is not a vertex of $\sigma$.
\item $\gamma s_{i_0}$ and $\gamma's_{j_0}$ are both vertices of both
$\sigma$ and $\sigma'$. But then, by definition of $i_0$ and $j_0$,
$i_0=j_0$ and $\rho(\sigma)=\gamma s_{i_0}=\gamma's_{j_0}=
\rho(\sigma')$, so that condition (2) of definition \ref{df-polar} again
applies. 
\end{enumerate}
\epn

\bpn{of proposition \ref{pr-triang}}
Let $C$ be a cellulation of $M$, with a polarization $\rho$. We first
built a $\pi_1M$-invariant triangulation of $\Ft_2$ as follows.

Let $f\in \Ft_2$ be such that some vertex $s$ of $f$ is
$\rho(\sigma)$, where $\sigma$ is one of the cells bounded by $f$. 
Then $s$ is the unique vertex of $f$ with this property, because of
condition (2) of definition \ref{df-polar}. Define a triangulation of
$f$ by adding the edges going from $s$ to all the other vertices of
$f$. Repeat this for all the 2-faces with this property.

Then subdivide all the remaining non-triangular 2-faces of $F_2$, so as
to obtain an equivariant triangulation of $\Ft_2$. 

Finally, define a triangulation of $M$ by subdividing each cell $\sigma$
of $\Ft_3$ by adding triangles containing $\rho(\sigma)$ and any edge of
$\sigma$ not containing $\rho(\sigma)$. It is clear that:
\begin{enumerate}
\item this defines an equivariant decomposition of $\Mt$ into simplices,
which is
obtained by subdividing each cell into simplices, and thus a
decomposition of $M$ into a finite number of simplices with disjoint
interior;
\item the simplices are non-degenerate;
\item if $\sigma\in \Ft_3$ and $\sigma'\in \Ft_3$ are two adjacent
simplices, then $\sigma\cap\sigma'$ is a face in both of them, because
it has to be one of the triangles of the triangulation of $\Ft_2$ obtained
above. 
\end{enumerate}
The definition then directly shows that the triangulation obtained in
this manner from a large cellulation is itself large. 
\epn

\subsection{Some elementary combinatorics}

We now fix a triangulation $\Sigma$ of $M$, with which we will stick
until section 8. We call $f$ the number of its 3-simplices, $t$ the
number of its 2-faces, $e$ the number of its edges, $e_i$ and $e_b$ the
number of interior and boundary edges respectively, and $v$ the number
of vertices. We will need later on the following easy consequence of the
Euler formula.

\blm  \label{euler}
$2 f = 2e_i + e_b - v$. 
\elm

\bpv
Consider the closed triangulated manifold obtained by gluing two copies
of $(M, \Sigma)$ along their boundary by the identity map. This
triangulated manifold has $\fb:=2f$ simplices, $\tb$ 2-faces,
$\eb:=2e_b+e_i$ edges, and $\vb:=v$ vertices. 

Since the Euler characteristic is $0$ in odd dimensions: 
$$ \fb - \tb + \eb - \vb = 0~. $$
Moreover, each 2-face bounds two simplices, and each simplex has 4
faces, so that:
$$ \tb = 2\fb~. $$
Therefore:
$$ -\fb + \eb - \vb = 0~, $$
and the result follows. 
\epv

\section{Hyperbolic structures on triangulated manifolds}

This section contains the definitions and basic properties concerning
some simple notions of singular hyperbolic structures. The idea is to
construct such structures by gluing ideal simplices, and then to show
that the set of those structures with some constraints actually contains
a smooth hyperbolic metric. This will be done in the next section using 
a variational argument.

The ideas used here were mostly developed previously for ideal polyhedra
in $H^3$. Their history is interesting. The first results were obtained
by Andreev \cite{andreev-ideal}, and then developed by Thurston
\cite{thurston-notes}. Colin de Verdi{\`e}re \cite{CdeV} then noted that the
results could be recovered using a variational approach, while Br{\"a}gger
\cite{bragger} identified the functional as the volume. Rivin
\cite{Ri2} then further developed the theory.

\subsection{Sheared hyperbolic structures}

We consider here a triangulation $\Sigma$ of $M$.

\bdf \label{df-struct}
A {\bf sheared hyperbolic structure} on $(M, \Sigma)$ is the choice, for each
simplex $S$ in $\Sigma$, of a diffeomorphism from $S$ to an ideal
simplex. We denote by $\Hsh$ the set of hyperbolic structures on $(M,
\Sigma)$. 
\edf

This vocabulary is justified by the elementary remark that, since all
ideal triangles in $H^2$ are isometric, there is a unique way of gluing
the hyperbolic simplices which are given by a sheared hyperbolic
structure along their common faces. One then obtains a hyperbolic metric
on the complement of the interior edges of $\Sigma$. But this metric
does not extend over the edges; rather, the model of what happens along
an edge is obtained by taking the quotient of the universal cover of
$H^3$ minus a geodesic $g$ by the group generated by the composition of
a rotation of angle $\theta$ around $g$ and a translation of length
$\delta$ along $g$. We then call $\theta$ the {\bf angle} around the
edge, and $\delta$ the {\bf shear} along the edge. 

Those hyperbolic structures will be considered up to diffeomorphisms
acting on the simplices; therefore, the choice of a sheared hyperbolic
structure 
is equivalent to the choice of the three dihedral angles of each
simplex, subject to the condition that their sum is $\pi$. So $\cH$ can
be identified with the product of $f$ 2-simplices, and it has a natural
affine structure.

\bdf \label{df-singular}
A sheared hyperbolic structure on $(M, \Sigma)$ is a {\bf singular
hyperbolic structure} if the 
shear of all the interior edges vanishes. It is {\bf smooth} if, in
addition, the singular angle of all interior edges is $2\pi$. The set of
singular hyperbolic structures is denoted by $\Hsi$, the set of smooth
hyperbolic structures by $\Hsm$.
\edf

In other words, a singular hyperbolic structure on $(M, \Sigma)$ defines
a hyperbolic cone-manifold structure on $M$, which is singular on the
edges of $\Sigma$. A smooth hyperbolic structure defines a hyperbolic
metric on $M$. In both cases the boundary is piecewise totally
geodesic. 

The shear at an interior edge of $\Sigma$ can be understood in the
following elementary way. 

\bdf \label{df-shear}
Let $s$ be an ideal simplex, and let $e$ be an edge of $s$. Choose an
orientation of $e$, and let $f_1$ and $f_2$ be the faces of $s$
containing $e$, in the order defined by the orientation of $e$. Let
$x_1$ and $x_2$ be the orthogonal projections on $e$ of the vertices of
$f_1$ and $f_2$ respectively which are not in $e$. The {\bf shear} of
$s$ at $e$ is the oriented distance between $x_1$ and $x_2$. 
\edf

Note that the shear of a simplex at an edge is clearly independent of
the orientation chosen.

\brk \label{rk-shear}
Let $h\in \Hsh$. The shear of $h$ at an interior edge $e$ is
the sum, over the simplices containing $e$, of their shears at $e$.   
\erk

\subsection{Angles}

For each edge $e$ of an ideal triangulation of a sheared hyperbolic
structure, we define the {\bf angle} at $e$ to be the sum of the dihedral
angles at that edge of the simplices containing it -- this applies to
interior as well as to boundary edges. For boundary edges this will also
be called the {\bf interior dihedral angle}, and the {\bf exterior
  dihedral angle} is
$\pi$ minus the interior dihedral angle. For interior edges, the {\bf excess
angle} is the angle minus $2\pi$, and the {\bf singular curvature}
around the corresponding edge is minus the excess angle.

Note that, for $h\in \Hsh$, the sum of the exterior angles at a vertex is
equal to $2\pi$ plus the sum of the excess angles at the interior
angles. This is checked by applying the Gauss-Bonnet formula to the link
of the vertex, which is piecewise Euclidean manifold. 

\bdf \label{df-angles}
Let $\theta:\Sigma_1\rightarrow \R_+$ be an assignment of ``angles'' to
the edges of $\Sigma$. We will say that $\theta$ is ``ideal'' if:
\begin{itemize}
\item the angles assigned to boundary edges are in $(0,\pi)$;
\item the angles assigned to interior edges are in $(0,2\pi)$;
\item at each vertex, the sum of the exterior dihedral angles of the
boundary edges equals $2\pi$ plus the sum of the angle excess at the
interior angles. 
\end{itemize}
The set of ideal angle assignments is denoted by $\Theta$.
For each $\theta\in \Theta$, we denote by $\Hsh(\theta)$ the set of
sheared hyperbolic structures on $\Sigma$ such that the angles
associated to each edge is given by $\theta$; then
$\Hsi(\theta):=\Hsh(\theta)\cap \Hsi$, and
$\Hsm(\theta):=\Hsh(\theta)\cap \Hsm$. 
\edf

\blm \label{lm-Theta}
If the triangulation $\Sigma$ is large, then $\Theta$ corresponds to the
interior of a polytope of dimension $e-v$ 
in $\R^e$.
\elm

\bpv
The only point is to prove that the constraints on the vertices are
linearly independent. So let $C$ be linear combination of those
constraints which is zero. In other terms, $C$ is a function
$C:\Sigma_0\rightarrow \R$ 
such that, for each oriented edge $e$, $C(e_-)+C(e_+)=0$. Moreover we
supposed that $\Sigma$ is large, so it is quite
obvious that the values of $C$ at the vertices of a triangular face of
$\Sigma$ have to be $0$, and therefore that $C=0$. 
\epv

\bdf \label{df-Thetasm}
We say that an ideal angle assignment is {\bf smooth} if the angle
assigned to each 
interior edge is $2\pi$. The set of those angle assignments is denoted
by $\Thetasm$. 
\edf

\section{First order variation of the volume}

The main goal of this section is to use elementary properties of the
volume --- seen as a functional on the space $\Hsh$ of sheared hyperbolic
structures --- to prove the existence of singular hyperbolic structures
with given angles on the interior and boundary edges. In the next
section, the results and some similar arguments will be used to
understand the smooth hyperbolic 
structures with given dihedral angles at the boundary edges.

\subsection{Definitions and first properties} 

We consider here again an ideal hyperbolic manifold, along with a
triangulation $\Sigma$, and an angle assignation
$\theta\in \Theta$.
We first define the volume of a singular hyperbolic structure in the
most obvious way. 

\bdf \label{df-vol}
For $h\in \Hsh$, the {\bf volume} of $h$, $V(h)$, is the sum of the
hyperbolic volumes of the simplices of $\Sigma$.
\edf

As an immediate consequence of lemma \ref{lm-vol}, we find that:

\blm \label{lm-cvx}
$V$ is a concave function on $\Hsh$ (with the affine structure coming
from the parametrization by the dihedral angles of the simplices).
\elm

We now need to understand how the volume varies when one deforms a
sheared hyperbolic structure. Unfortunately the Schl{\"a}fli formula
(\ref{schlafli}) does not apply directly to sheared hyperbolic
structures, and does not even make sense in this case; indeed there is no
way to choose a horosphere centered at a given vertex, since the holonomy
around an edge $e$ would act on it by translation along the edge, with a
translation distance equal to the shear at $e$.

To understand this point better, we fix a sheared hyperbolic structure $h\in
\Hsh(\theta)$, and a first order variation $\hbu\in
T_h\Hsh(\theta)$. Suppose given an ideal triangulation of $h$. $\hbu$
determines a first order variation of the dihedral angles of the
simplices, to which the Schl{\"a}fli formula (\ref{schlafli}) can be
applied, once a horosphere around each vertex is chosen for each
simplex. To get a better understanding of the first order variation
of the
volume, we can choose, for each vertex $v$, a horosphere centered on $v$
for each simplex containing $v$. We call this collection of horospheres
a {\bf choice of horospheres} at $v$.

\bdf \label{df:var-vol}
Let $e$ be an edge of $\Sigma$, with vertices $e_-$ and $e_+$. Let $H_-$
and $H_+$ be choices of horosphere at $e_-$ and $e_+$ respectively. The
{\bf volume variation at $e$} associated to $\hbu$ is:
$$ \sum_i L_i \alphab_i~, $$
where the sum is over the simplices containing $e$, $L_i$ is the
oriented length of the part of $e$ which is between the horospheres in
$H_-$ and $H_+$, and $\alphab_i$ is the first-order variation of the
dihedral angle at $e$.
\edf

\subsection{Volume differential and shears}

Of course the point is that, once horospheres are chosen,  the volume
variation at $e$ can be seen as the contribution coming from $e$ to the
first order variation of the volume; indeed, the Schl{\"a}fli formula
(\ref{schlafli}) indicates that:

\brk
If a horosphere choice is given for all vertices of $\Sigma$, the first
order variation of the volume under the deformation $\hbu$ is the sum of
the volume variations at the edges.
\erk

\bdf \label{df:coherent}
Let $T$ be a 2-face of the triangulation $\Sigma$, and let $v$ be a vertex
contained in $\Sigma$. A choice of horospheres at $v$ is {\bf coherent} at
$T$ if the horospheres on both simplices containing $T$ have the same
intersection with $T$.
\edf

The point is that, as explained above, if $h$ has a non-zero shear at an
edge $e$ containing $v$, then it is not possible to find a horosphere
choice at $v$ which is coherent on all 2-faces containing $e$.

\bprop \label{var-vol}
Let $e$ be an oriented interior edge of $\Sigma$, with endpoints $e_-$
and $e_+$. Let $f_1, \cdots, f_N$ be the 2-faces containing $e$,
oriented in cyclic order. Choose $i,j\in \{ 1,\cdots, N\}$. Let $H_-$ be
a choice of horospheres at $e_-$ around $e$ which is coherent except at
$f_i$, and let $H_+$ be a choice of horospheres at $e_+$ around $e$
which coherent except at $f_j$. Let $\hbu$ be an infinitesimal  variation
of $h$ which does not change the total angle around the interior edges. The variation of the volume at $e$
associated to $\hbu$ is:
$$ \Vb_e = s_e \sum_{k=j}^{i-1} \alphab_k~, $$
where $s_e$ is the shear of $h$ at $e$, and $\alpha_k$ is the angle
between $f_k$ and $f_{k+1}$.
\eprop

\bpv
The definition of $H_-$ and $H_+$ shows that there are two numbers $L,
L'\in \R$ such that:
\begin{itemize}
\item for each $k\in \{ i, i+1, \cdots, j-1\}$, the distance along $e$
  between the horospheres of $H_-$ and $H_+$ in the simplex having faces 
  $f_k$ and $f_{k+1}$ is $L$.
\item for each $k\in \{ j, j+1, \cdots, i-1\}$, the same distance is
  $L'$. 
\end{itemize}
Then $L'-L=s_e$, again by definition of $H_-, H_+$ and $s_e$. So, by
definition \ref{df:var-vol}:
$$ \Vb_e = \sum_{k=i}^{j-1} L\alphab_k + \sum_{k=j}^{i-1}
L'\alphab_k = \sum_{k=j}^{i-1} L\alphab_k + \sum_{k=j}^{i-1}
(L+s_e)\alphab_k~. $$ 
But the total angle around $e$ remains constant in the variation $\hbu$,
which means that:
$$ \sum_k \alphab_k =0~. $$
Subtracting $L$ times this equation to the previous one leads to the
statement. 
\epv

As a consequence we see that the first-order variation of the volume,
under a deformation which does not change the total angle around the
interior angles or the dihedral angles, has
a remarkably simple form, more precisely it can be expressed only in
terms of the shears at the interior edges. 

\blm \label{dV-shear}
$dV(\hbu)$, as a linear form on $\Hsh(\theta)$, depends only on the shear
at the interior edges of $\Sigma$. 
\elm

\bpv
This is an immediate consequence of the previous proposition, and of the
fact that the dihedral angle at the exterior edges is constant. 
\epv

This means in particular
that the volume --- seen as a function on $\Hsh(\theta)$ --- is critical
when the shear vanishes at all interior edges. We will see below that
the converse is true too. 
This will use a special type of deformations, defined as follows.

\bdf \label{df-local}
Let $h\in \Hsi$, and let $e\in \Sigma_1$ be an oriented interior edge of
$\Sigma$. A {\bf local} deformation $\hbu_0$ of $h$
at $e$ is determined as follows. Choose a face $f$ containing $e$. Then
let $s_0, s_1, s_2, \cdots, s_N=s_0$ be the simplices of $\Sigma$
containing $e$, 
in cyclic order. Let $a_i^+$ and $b_i^+$ be the edges of $s_i$ containing
$e_+$ (other than $e$), and $a_i^-$ and $b_i^-$ the edges of $s_i$
containing $e_-$ (again other than $e$),
ordered so that $a_i^+=b_{i+1}^+$ and that $a_i^-=b_{i+1}^-$. Then, in
the deformation $\hbu_0$:
\begin{itemize}
\item the angles of the simplices $s_i$ at $e$ do not change;
\item the angles of $s_i$ at $a_i^+$ and $b_i^-$ vary at speed $1$;
\item the angles of $s_i$ at $a_i^-$ and at $b_i^+$ vary at speed $-1$.
\end{itemize}
\edf

\begin{figure}[h]
\centerline{\psfig{figure=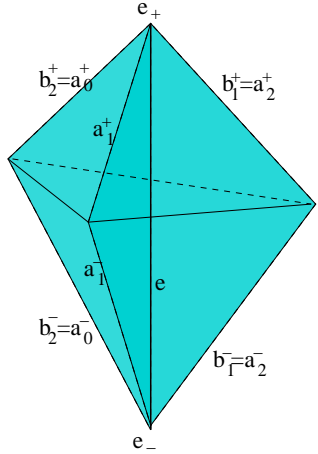,height=6cm}}
\caption{Local deformations}\label{fig:defos}
\end{figure}

It is a simple matter to check the following properties of local
deformations: 

\brk \label{rk:local}
The above definition indeed defines a deformation of $h$ in
$\Hsh$. Moreover, it does not change the total angles at the edges, i.e. if
$h\in \Hsh(\theta)$, then $\hbu\in T_h\Hsh(\theta)$.
\erk

The local deformations can then be used to prove that, among the sheared
hyperbolic structures, those which have zero shear are exactly the
critical points of the volume, seen as a functional on the space of
sheared hyperbolic structures having a given total angles on all edges. 

\blm \label{lm-dV0}
Let $\theta\in \Theta$, and let $h\in \Hsh(\theta)$; then the
restriction of $dV$ to $T_h\Hsh(\theta)$ is zero if and only if the
shear of $h$ at all interior edges is zero, that is, if and only if
$h\in \Hsi$.
\elm

\bpv
Lemma \ref{dV-shear} shows that the volume --- as a functional on
$\Hsh(\theta)$ --- is critical when the shear is zero at all interior
edges. Remark \ref{rk:local} shows that the local deformations are
tangent to $\Hsh(\theta)$. 

Suppose that the shear at an edge $e$ is non-zero, and consider the
simplices containing $e$, as in figure 3. It is possible to choose
horospheres at the vertices of those simplices such that:
\begin{itemize}
\item at $e_+$ and $e_-$, the choices are coherent at the triangles
  $(e,a^+_i,a^-_i)$, except for $i=0$; 
\item at $e_+$ and $e_-$, the distance between the two horospheres on
  each side of the triangle $(e,a^+_0,a^-_0)$ is equal to the shear at
  $e$;
\item for each $i\neq 0$, the two horospheres at the vertex $a_i^-\cap
  a_i^+$ are coherent.
\end{itemize}
Consider a local deformation at $e$. The only contributions to the
variation of the volume comes from the $a_i^-$ and the $a_i^+$, since
they are the only edges where the angles vary. Moreover they add up to zero
except at $a_0^-$ and $a_0^+$, because of the choices of horospheres
described above (they are coherent on both sides of the other
edges). Finally, the contributions from $a_0^-$ and $a_0^+$ are non-zero
and of the same sign. 

So proposition \ref{var-vol} shows that
the shears have to be zero at all interior edges at critical points
of $V$.
\epv

\subsection{Consequences}

The results of the previous subsection now easily lead to interesting
results on singular hyperbolic structures (i.e. those without shears).

\bcr \label{cr:unique-si}
Let $\theta_0\in \Theta$ and $g_0\in \Hsi(\theta_0)$. There exists a
neighborhood $U$ of $\theta_0$ in $\Theta$ and a neighborhood $V$ of
$g_0$ in $\Hsm$ such that, for any $\theta\in U$, there exists a unique
$g\in \Hsi(\theta)\cap V$. 
\ecr

\bpv
By corollary
\ref{cr:concavity}, $V$ is a strictly concave function on
$\Hsh$. Moreover, the restriction of $V$ to $\Hsh(\theta_0)$ has a
critical point at $g_0$. Therefore, for $\theta$ close enough to
$\theta_0$, $V$ has a unique critical point of $V$ on $\Hsh(\theta)$
close to $g_0$; by lemma \ref{lm-dV0} it is a point in $\Hsi(\theta)$.  
\epv

\bcr \label{cr:dim}
For each $\theta\in \Theta$ such that $\Hsh(\theta)$ is non-empty,
$\Hsh(\theta)$ is a submanifold of $\Hsh$, with dimension equal to the
number $e_i$ of interior edges of the triangulation. 
\ecr

\bpv
The space of ideal simplices in $H^3$ has dimension $2$, so that the
dimension of $\Hsh$ is equal to twice the number of simplices in
$\Sigma$. Thus, by lemma \ref{euler}, $\dim\Hsh=2e_i+e_b-v$, where $e_b$
and $v$ are the number of boundary edges and of vertices of $\Sigma$
respectively. 

Specifying the dihedral angles and the excess angle at the interior
edges adds $e_i+e_b$ constraints; they are not linearly independent,
however, since they satisfy at least one linear condition for each
vertex. The loss of dimensions due to those constraints is therefore at
most $e_i+e_b-v$, so that, for $\theta\in \Theta$, $\dim
\Hsh(\theta)\geq e_i$. 

But lemma \ref{dV-shear} indicates that $dV$ depends only on the $e_i$
shears at the interior edges, while the strict concavity of $V$ shows
that the Jacobian of $dV$ is non-degenerate. Therefore, its restriction
to $\Hsm(\theta)$ is also non-degenerate, so that
$\dim\Hsh(\theta)\leq e_i$, and the result follows.
\epv

Finally, using the non-degeneracy of the volume functional also leads to
an infinitesimal rigidity result.

\bcr \label{cr:dim-Hsi}
Let $h\in \Hsi$, and let $\hbu\in T_h\Hsi$ be a first-order deformation
of $h$. Then $\hbu$ corresponds to a non-zero first-order deformation of
either the dihedral angles of the boundary edges, or the excess angle
at the interior edges. $\Hsi$ is, in the neighborhood of $h$, a
manifold of dimension $e_i+e_b-v$.
\ecr

\bpv
This is a direct consequence of corollary \ref{cr:dim}, which shows that
$\Hsh$ is foliated by the $\Hsh(\theta)$, which are submanifolds of
dimension $e_i$. We already know that $\Hsh$ is a manifold of dimension
$2f$, and, by lemma \ref{euler}, $2f=2e_i+e_b-v$, where $e_b$ is the
number of boundary edges and $v$ is the number of vertices. But the
dimension of the space $\Theta$ of possible angle assignations
is at most $e_b+e_i-v$, since the dihedral angles at the boundary edges
satisfy one linear condition for each vertex. Thus this dimension has to
be exactly $e_b+e_i-v$, and each non-trivial deformation in $\Hsi$ ---
which is transverse to the foliation of $\Hsh$ by the $\Hsh(\theta)$ ---
induces a non-trivial deformation of some dihedral angle on the boundary
or of some total angle on an interior edge.
\epv

Finally, an elementary dimension-counting argument, along with the
previous corollary, shows that any infinitesimal deformation of the
dihedral angles on the boundary --- and of the angle excess around the
interior edges --- can be realized by some $\hb\in T_h\Hsi$, if it
satisfies the condition that, at each vertex, the sum of the exterior
dihedral angles remains equal to $2\pi$ plus the sum of the angle excess
at the interior edges. 
Restricting to the case of smooth hyperbolic structures,
this can be formulated as follows. We call {\bf admissible} an
infinitesimal deformation of the dihedral angles such that the sum of
the exterior angles of the edges at any vertex remains $2\pi$.

\blm \label{lm:defos-possibles}
Let $(M,h)$ be an ideal hyperbolic manifold. For any admissible infinitesimal
variation $\thetabu$ of the dihedral angles on the boundary, there
exists a (unique) infinitesimal deformation of $h$ in $\Hsi$ which
induces the variation $\thetabu$. 
\elm

\bpv
By lemma \ref{lm-atri}, there exists a finite cover $\Mb$ of $M$ so that
the hyperbolic structure $\hb$ lifted to $\Mb$ from $h$ admits an ideal
triangulation, say $\Sigma$. Let $\thetab$ be the dihedral angles of
$\hb$. The previous corollary shows that there is
no non-trivial deformation of $\hb$ which does not change the dihedral
angles, and that the dimension of $\Hsm$ in the neighborhood of $\hb$ is
$e_b-v$. Thus, for each admissible admissible infinitesimal deformation
$\thetabbu$ of the 
dihedral angles of $\hb$, there is a unique infinitesimal deformation
$\hbbu$ of $\hb$ inducing $\thetabbu$. 

If now $\thetabu$ is an infinitesimal deformation of $\theta$, it lifts
to an infinitesimal deformation $\thetabbu$ of the dihedral angles of
$\hb$, which is induced by a unique admissible infinitesimal deformation
$\hbbu$ of 
$\hb$. Since $\thetabbu$ is equivariant, 
$\hbbu$ is obviously equivariant (by uniqueness) and therefore
defines an admissible infinitesimal deformation $\hbu$ of $h$ inducing
$\thetabu$, which is unique since $\hbbu$ is.
\epv

Note that it is not necessary to suppose that the hyperbolic ideal
manifold which we consider has a boundary which is triangulated,
i.e. this lemma is also valid when some of the faces have more than $3$
edges. In this case one can add edges to those faces, so as to obtain
triangulation, and then choose the first-order variation of the dihedral
angles at those new edges. In particular it is possible to choose
first-order variations which "lose" the convexity of the boundary.

\section{Dihedral angles}

We consider here an ideal hyperbolic manifold $M$. The first point will
be a compactness result for 
sequences of ideal structures on $M$ with a given triangulation of the
boundary, when the dihedral angles converge. 
We will then consider a fixed triangulation of $M$, then see what happens
when one considers only a fixed cellulation $\sigma$ of the boundary,
and finally how the affine structures corresponding to different
triangulations of the boundary can be glued together. 

\subsection{Compactness of ideal manifolds}

We consider a fixed triangulation $\sigma$ of $\dr M$; we will give in
this subsection a compactness lemma for ideal structures with boundary
combinatorics corresponding to $\sigma$, then use it in the next
sub-sections to obtain a description of the possible dihedral angles. 

\blm \label{lm:comp-ideal}
Let $(h_i)$ be a sequence of ideal structures on $M$ with boundary
combinatorics given by $\sigma$. Suppose that the dihedral angles of the
boundary vertices for the $h_i$ converge as $i\rightarrow \infty$ to
limits in $(0, \pi)$. Then:
\begin{itemize}
\item either there is a non-elementary circuit in $\sigma$ on which the
  sum of the limit dihedral angles is $2\pi$;
\item or there exists a sub-sequence of $(h_i)$ which converges to an
  ideal structure on $M$. 
\end{itemize}
\elm

The intuitive idea behind this lemma is simple. Consider the induced
metrics on one of the boundary components. For each $n$, it is a
complete hyperbolic metric of finite area. So, by a classical result in
Teichm{\"u}ller theory, there is a subsequence which either converges to a
hyperbolic metric of finite area, or has a closed geodesic whose length
goes to $0$. In the second case, the short geodesic will --- by a simple
convexity argument --- correspond to a non-elementary circuit on which
the sum of the exterior dihedral angles goes to $2\pi$. If the first
case happens for all boundary components, then both the induced metric
and the dihedral angles converge, so that the universal cover in $H^3$
of each boundary component converges, and thus the sequence of metrics
on $M$ has to converge. 

The "formal" proof, however, is a little more complicated since one has
to keep track of the relationship between the induced metrics and the
triangulation; it will use a couple of definitions and propositions. 
First there is a natural notion of collapsing of vertices for the
sequence $(h_i)$. We consider the triangulation $\sigma$ of $\dr M$, and
its lift to a triangulation $\sigmat$ of the universal cover of
$\dr M$. Note that each hyperbolic structure $h_i$ on $M$ defines a
hyperbolic structure on the universal cover $\Mt$ of $M$, and thus a
polyhedral embedding of $\tilde{\dr M}$ in $H^3$.

\bdf \label{df:collapse}
Let $S$ be a finite subset of the set of vertices of $\sigmat$, with
cardinal at least $2$. We say that $S$ {\bf collapses} if there exists
a subsequence of $(h_i)$ and a sequence of conformal maps
$\rho_i:H^3\mapsto B_0(1)\subset \R^3$ such that:
\begin{itemize}
\item the set $\rho_i(S)$ converges to a point $x_0$ in $S^2$;
\item no vertex of $\sigmat$, adjacent to a vertex in $S$, converges to
  $x_0$. 
\end{itemize}
\edf

Here $B_0(1)$ is the radius $1$ ball centered at $0$ in $\R^3$.
This notion is useful because it implies the existence of a
non-elementary circuit in $\sigma$ on which the sum of the dihedral
angles goes to $2\pi$. 

\bprop \label{pr:collapse-pts}
If a finite subset $S$ of vertices of $\sigma$ (with cardinal at least $2$)
collapses, then, maybe after taking a subsequence of $(h_i)$, there is a
non-elementary circuit in the 1-skeleton of $\sigma$ on which the sum of
the dihedral angles goes to $2\pi$. 
\eprop

\bpv
The argument is the same as in section 4 (and that in \cite{shu, cpt}):
the polyhedron dual to $\tilde{\dr M}$ in de Sitter space, for the limit
hyperbolic structure, has a face tangent to the boundary at infinity at
the point where the collapse occurs, and the metric induced on this face
is the metric of a hemisphere. Therefore the length of its boundary is
$2\pi$, and this provides a non-elementary circuit in $\sigmat$ with sum
of the dihedral angles converging to $2\pi$. 
\epv

Note that the non-elementary circuit on which the sum of the dihedral
angles goes to $2\pi$ is simply the sequence of edges between (a
connected component of) $S$ and its complement. 

We now have to define the kind of metrics that appear on the boundary
of an \ihm $M$. Here $S$ is a compact orientable surface of
genus $g\geq 2$, with $N$ points $v_1, \cdots, v_N$ removed ($N\geq
1$), along with a triangulation $\sigma$ whose vertices are the
$v_i$ --- for instance, $S$ can be $\dr M$ with the vertices removed. 
We call $a_1, a_2, \cdots, a_e$ the edges of $\sigma$, and
$\tau_1, \cdots, \tau_t$ its triangles.

\bdf \label{df-hyp}
A {\bf hyperbolic structure} on $(S, \sigma)$ is a family of gluings
between the adjacent triangles of $\sigma$, 
each being endowed with the metric of the ideal triangle in $H^2$.
The set of hyperbolic structures on $(S,\sigma)$ will be denoted by
$\cM$. 
\edf

Obviously a hyperbolic structure in this sense defines a hyperbolic
metric on $S$, but it is in general not complete; 
indeed, there might be a ``shift'' at some of the edges, in the sense
described below.

\bdf \label{shift}
Let $g$ be a hyperbolic structure on $(S, \sigma)$, and let $E$ be an
edge of $\sigma$. Choose an orientation of $E$, and let $T_+$ and $T_-$
be the triangles of $\sigma$ standing on the ``right'' and on the
``left'' of $E$ respectively. Let $u_+$ and $u_-$ be the orthogonal
projections on $E$ of the vertices opposite to $E$ in $T_+$ and $T_-$
respectively. The {\bf shift} $\sh(E)$ of $g$ at $E$ is the oriented
distance, along $E$, between $u_-$ and $u_+$. 
\edf 

Note that $\sh(E)$ does not depend on the orientation chosen for 
$E$. This definition is somehow related to the notion 
of ``shear'' defined above for an ideal simplex at an edge, but it is
better to use a different name to avoid confusions.

\bprop \label{pr:shifts-lim}
Suppose that no finite set of vertices of $\sigma$ collapses. Then,
after taking a subsequence, the shifts of the edges of $\sigma$ for the
$h_i$ converge. 
\eprop

\bpv
Choose a "fundamental domain" in $\sigmat$, i.e. a connected subgraph
$\sigma_0\subset \sigmat$ such that each edge in $\sigma$ has a unique
inverse image in $\sigma_0$ for the canonical projection $\sigmat\mapsto
\sigma$. 

We consider again the polyhedral embeddings $\phi_i:\tilde{\dr
M}\rightarrow H^3$ associated to the ideal hyperbolic structures $h_i$
on $M$. 
Since $\sigma_0$ has a finite number of edges, an elementary compactness
argument shows that there exists a sequence $(\rho_n)$ of conformal maps
from $H^3$ to the radius $1$ ball in $\R^3$ such that the images by
$\rho_i\circ \phi_i$ of the vertices of $\sigma_0$ converge in $S^2$. 

Since by hypothesis there is no collapsing subset of vertices, there
exists a constant $C>0$ such that the edges of the image in $S^2$ of
$\sigma_0$ by the limit have length between $1/C$ and $C$. An elementary
argument then shows that the distance on $S^2$ between the vertices of
$\sigma_0$ and the orthogonal projections on the edges of the triangles
of the opposite vertices remain bounded between $1/C'$ and $C'$, where
$C'>0$ is another constant.

This in turns indicates that (after taking a subsequence) the shifts of
$\sigma$ converge.  
\epv

\bpn{of lemma \ref{lm:comp-ideal}}
By hypothesis the dihedral angles of the $(h_i)$ converge, moreover the
metrics induced on the boundary converge by proposition
\ref{pr:shifts-lim}. Thus the polyhedral embeddings of $\tilde{\dr M}$ in
$H^3$ converge, and therefore the hyperbolic structures $(h_i)$ also do
so. 
\epn

\subsection{Dihedral angles with a given triangulation}

We now have all the tools necessary to describe the dihedral angles of
ideal manifolds with a given triangulation of the boundary. We consider
a triangulation $\Sigma$ of $M$. 

\blm \label{lm:triangulation-int}
Consider a triangulation $\Sigma$ of the interior of $M$. Let
$\Hsm^\Sigma$ be the set of smooth hyperbolic structure obtained by
gluing ideal simplices according to $\Sigma$. Then:
\begin{enumerate}
\item the space of possible dihedral angles forms the interior of a
  convex polyhedron $\Theta_\Sigma$ in $\R^{e_b}$.
\item for any $\theta\in \Theta_\Sigma$, there is a unique smooth
  hyperbolic structure on $M$ with those dihedral angles; $\Hsm^\Sigma$
  therefore inherits the affine structure of $\Theta_\Sigma$.
\item on each point of the boundary of $\Theta_\Sigma$, one of the
  following happens:
  \begin{enumerate}
    \item one of the simplices of $\Sigma$ has a dihedral angle which
    goes to $0$;
    \item the dihedral angle of one of the boundary edges goes to $0$ or
    to $\pi$;
    \item there is a non-elementary circuit in the 1-skeleton of
    $\Sigma$ for which 
    the sum of the dihedral angles goes to $2\pi$.
  \end{enumerate}
  \item the volume $V$ is a smooth, strictly concave function on
    $\Theta_\Sigma$. 
\end{enumerate}
\elm

The proofs of points (1) and (2) are consequences of the previous
section. Indeed the condition that the dihedral angles are given by
$\theta$ defines a family of linear conditions on the 
angles, so an affine subspace of the space of angle assignations of the
simplices of $\Sigma$. Adding the condition that the excess angle at
each interior edge is $2\pi$ defines an affine subspace, and thus a smaller
polyhedron in the space of angle assignations on the simplices of
$\Sigma$. Projecting to the space of possible dihedral angles on the
boundary edges thus defines a convex polyhedron. 

Point (2) follows from a deformation argument. Consider the map $\Phi:
\Hsm^\Sigma\rightarrow \Theta_\Sigma$ sending an ideal hyperbolic structure to
its dihedral angles.
Corollary \ref{cr:unique-si} shows that $\Phi$ is locally injective (its
differential is injective), while lemma \ref{lm:comp-ideal} shows that
$\Phi$ is proper --- it is therefore a covering. But $\Hsm^\Sigma$ is
connected, while $\Theta_\Sigma$ is contractible, so $\Phi$ is a
homeomorphism. 
Therefore each assignation of dihedral angles in $\Theta_\Sigma$ is
indeed realized by a smooth hyperbolic structure. 

Point (3) is a consequence of lemma 8.1, which indicates that no
degeneration other than (a), (b) and (c) can occur.

Point (4) is also a consequence of the definitions, and of the following
remarks, which applies directly to this case.

\brk \label{rk:concavity}
Let $\Omega\in \R^N$ be a convex subset, and let $f:\Omega\rightarrow
\R$ be a smooth, strictly concave function. Let
$\rho:\R^N\rightarrow \R^p$ be a 
linear map, with $p<N$, and let $\Omegab:=\rho(\Omega)$. Define a
function:
$$
\begin{array}{cccc}
\fb: & \Omegab & \rightarrow & \R \\
& y & \mapsto & \max_{x\in \rho^{-1}(y)} f(x) 
\end{array}
$$
Then $\Omegab$ is convex, and $\fb$ is a smooth, strictly
concave function on $\Omegab$.
\erk

\bpv
It is quite obvious that $\Omegab$ is convex, and also that $\fb$ is
smooth since $f$ is strictly concave.

Let $\cb:[0,1]\rightarrow \Omegab$ be a geodesic segment, parametrized
at constant speed. By definition of $\fb$, there exist points $x_0,
x_1\in \Omega$ such that:
$$ \cb(0)=\rho(x_0), ~ \cb(1)=\rho(x_1), ~ \fb\circ \cb(0)=f(x_0), ~
\fb\circ \cb(1)=f(x_1)~. $$
Let $c:[0,1]\rightarrow \Omega$ be the geodesic segment parametrized at
constant speed such that $c(0)=x_0$ and $c(1)=x_1$. Since $\rho$ is
linear, $\cb=\rho\circ c$. 

Moreover, since $f$ is strictly concave
$$ \forall t\in (0,1), ~ f\circ c(t)>t f\circ c(0) + (1-t) f\circ
c(1)~. $$
Therefore, the definition of $\fb$ shows that:
$$ \forall t\in (0,1), ~ \fb\circ \cb(t)\geq f\circ c(t)>t f\circ c(0) +
(1-t) f\circ c(1)= \fb\circ \cb(0)+\fb\circ \cb(1)~. $$
This shows that $\fb$ is strictly concave.
\epv

\subsection{Fixed triangulations of the boundary}

Using the ideas and results of the previous section, we can now give a
description of the set of dihedral angles which can be achieved on the
ideal hyperbolic manifolds with a given cellulation of the
boundary. Studying the possible degenerations and a global argument
using the results of section 3 (on the fuchsian case) will then lead to
theorem \ref{tm:dihedral}. So we now consider an ideal hyperbolic
manifold $M$, which does not necessarily admit an ideal
triangulation. Of course the point is that, by lemma \ref{lm-atri}, a
finite cover of $M$ admits one.

\blm \label{lm:triangulation-bord}
Consider a cellulation $\sigma$ of the boundary $\dr M$, and
let $\Hsm^\sigma$ be the set of ideal hyperbolic structures on $M$ such
that the boundary is triangulated according to $\sigma$. Then each
connected component of 
$\Hsm^\sigma$ has a natural affine structure $\cA_\sigma$, obtained by
gluing the $\Hsm^\Sigma$ --- for $\Sigma$ a triangulation of the interior
of $M$ inducing the triangulation $\sigma$ of the boundary --- in the
natural way on their intersections. $(\Hsm^\sigma, \cA_\sigma)$ is
affinely equivalent to the disjoint union of the interiors of a finite
set of convex polyhedra in
$\R^{e_b}$, and $V$ is a concave function on each.
\elm

\bpv
Let $\theta\in\Theta_\sigma$, and let $h\in \Hsm^\sigma(\theta)$. By
lemma \ref{lm:defos-possibles}, the infinitesimal deformations of $h$
are parametrized by the admissible infinitesimal deformations of the
dihedral angles at the exterior edges. 

Let $\Theta_\Sigma$ be the set of dihedral angles that can be attained
by deformation from $h$. By proposition \ref{pr:degen}, the boundary is
made of angle assignations satisfying one of a number of possible affine
equalities, so that $\Theta_\Sigma$ is the interior of a convex
polyhedron in  $\R^{e_b}$. Moreover, its boundary can be decomposed into
the union of two components, $\dr\Theta_\Sigma = \dr_t\Theta_\Sigma\cup
\dr_r\Theta_\Sigma$, where:
\begin{enumerate}
\item $\dr_r\Theta_\Sigma$ is the set of boundary points where
  the sum of the dihedral angles on the edges of a non-elementary
  circuit goes to $2\pi$, or the dihedral angle at some boundary edge
  goes to $0$ or to $\pi$.
\item $\dr_t\Theta_\Sigma$ is the set of boundary points where the
  previous condition does not apply, but one of the simplices has an
  angle which goes to $0$.
\end{enumerate}
 
By lemma \ref{lm:necessary}, $\dr_r\Theta_\Sigma$ corresponds to the
boundary of 
$\Hsm^\sigma$. On the other hand, let $\hb$ correspond to a point of
$\dr_t\Theta_\Sigma$.  
Lemma \ref{lm:defos-possibles} shows that $\hb$ can be deformed so that
its dihedral angles go beyond $\dr_t\Theta_\Sigma$. So there is another
ideal triangulation $\Sigma'$ of $(M, \hb)$ such that $\hb\in
\Hsm^{\Sigma'}$. 
Since the affine structure on $\Hsm^{\Sigma}$ is defined in terms of the
variations of the dihedral angles, it extends naturally to an affine
structure on $\Hsm^{\Sigma'}$. This shows that $\Hsm^\sigma$ carries an
affine structure $\cA_{\sigma}$, 
obtained by gluing the affine structures on the
$\Hsm^{\Sigma}$. 

Moreover, at a point of
$\dr_t\Theta_\Sigma\cap \dr_r\Theta_\Sigma$, the faces of
$\Hsm^{\Sigma'}$ corresponding to $\dr_r\Theta_\Sigma$ are the
extensions of the faces of $\Hsm^{\Sigma}$, so 
that $\Hsm^{\sigma}$ is locally convex for the affine structure. It is
therefore affinely equivalent to the
interior of a convex polyhedron in $\R^{e_b}$. 
\epv

In addition, lemma \ref{lm:triangulation-int} allows for a simple
description of the possible boundary behavior; since the cases
corresponding to the collapse of a simplex correspond to the boundary of
the cells $\Hsm^\Sigma$ rather than to the boundary of $\Hsm^\sigma$,
the boundary of $\Hsm^\sigma$ is characterized by one of the other
possible limit cases. 

\bprop \label{pr:degen}
Let $(h_n)_{n\in \N}$ be a sequence of elements of $\Hsm^\sigma$
converging to a limit $h\in \dr \Hsm^\sigma$. After taking a
subsequence, one of the following occurs:
\begin{enumerate}
\item some boundary edge has a dihedral angle which goes to $0$ or
  $\pi$;
\item there exists a non-elementary circuit in $\Sigma_{|\dr M}$ on
which the sum of the dihedral angles goes to $2\pi$.
\end{enumerate}
\eprop

\subsection{Spaces of angle assignations}

The results already obtained give us a good understanding of the local
properties of the map sending an ideal hyperbolic manifold to the angle
assignation on the boundary edges. 
To complete the proof of theorem \ref{tm:dihedral}, we need some global
topological information on the space of angle assignations on the
boundary satisfying the 
conditions of the theorem. 

We need some additional notations. We will call 
$(\dr_iM)_{1\leq i\leq N}$ the connected components of $\dr M$. 
For each $n\in (\N\setminus \{ 0\})^N$, we call $\cA_n$ the space of
couples $(\Gamma, w)$, where $\Gamma$ is the 1-skeleton of a cellulation
of $\dr M$ having $n_i$ vertices in each $\dr_iM$, and
$w:\Gamma_1\rightarrow (0, \pi)$ is a function satisfying the conditions
of theorem \ref{tm:dihedral}; and we call $\cM_n$ the space of ideal
hyperbolic manifolds having $n_i$ vertices in $\dr_iM$, for each
$i$. $\Phi_n$ will then be the natural map sending an element of $\cM_n$
to the corresponding dihedral angles assignation on a cellulation of
$\dr M$, which is in $\cA_n$.

We will prove that if $M$ has incompressible boundary, then, for each
$n$, $\cA_n$ is connected. Thus $\Phi_n$ is a covering.
To prove that it is actually a homeomorphism, we will consider some
special dihedral angle assignations for which the uniqueness follows
from the Mostow rigidity. 

\blm \label{lm:An}
Suppose that $M$ has incompressible boundary. Then,
for each $n\in (\N\setminus \{ 0\})^N$, $\cA_n$ is connected.
\elm

This is the only point where we use the fact that $M$ has incompressible
boundary. 

Recall that, for $g\geq 2$ and $r\geq 1$, $\cA_{g,r}$ is the space
of dihedral angles assignations for fuchsian manifolds of genus $g$,
with $r$ vertices. From now on, we call $g_i$ the genus of
$\dr_iM$. With the topological assumptions which we have made on $M$,
$g_i\geq 2$ for each $i$. 

\brk
Suppose that $M$ has incompressible boundary. Then,
for each $n\in (\N\setminus \{ 0\})^N$, $\cA_n$ is homeomorphic to
$\Pi_{i=1}^N \cA_{g_i, n_i}$. 
\erk

\bpv
The conditions in theorem \ref{tm:dihedral}, which
define both $\cA_n$ and the $\cA_{g_i, n_i}$, describe each connected
component of $\dr M$ independently from the others. Since $M$ has
incompressible boundary, the curves in each connected component of $\dr
M$ which are contractible in $M$ are those which are contractible in
$\dr M$; so the condition on the dihedral angles on $\dr_iM$ are exactly
those describing $\cA_{g_i, n_i}$.
\epv

\bprop \label{pr:agn}
For each $g\geq 2$ and each $r\geq 1$, $\cA_{g,r}$ is homeomorphic to
the space $\cT_{g,r}$ of conformal structures on a surface of genus $g$,
with $r$ marked points (up to isotopy). 
\eprop

\bpv
Let $\Sigma_g$ be a compact surface of genus $g$, and let $t$ be a
conformal structure on $\Sigma_g$. There is a unique hyperbolic metric
on $\Sigma_g$, say $h$, in the conformal class defined by $t$. Taking
the warped product:
$$ M_g := (\Sigma_g\times \R, dt^2 + \cosh^2(t) h) $$
determines a complete fuchsian hyperbolic manifold. If $p_1, \cdots,
p_r$ are points in $\Sigma_g$, they define points $p_{1, \pm }, \cdots,
p_{r, \pm }$ in $\dr_\infty M_g$, with $p_{i_+}$ in the ``upper''
boundary and $p_{i,-}$ in the ``lower'' boundary for each $i$, and
$p_{i_+}$ exchanged with $p_{i,-}$ by the isometric involution on
$M_g$. 

Let $\Mb$ be the convex hull of the $p_{i, \pm }$. By remark
\ref{rk:fuchsian-bent}, $\Mb$ is an 
ideal hyperbolic manifold. Moreover, each fuchsian ideal hyperbolic
manifold is obtained in this way, and thus is associated to a $t\in
\cT_{g,r}$. 

Finally, theorem \ref{tm:semi-ideal} shows that $\cA_{g,r}$ is
homeomorphic to the space of fuchsian ideal hyperbolic manifolds with
each boundary surface of genus $g$ with $r$ vertices, and the result
follows. 
\epv

\bpv[Proof of lemma \ref{lm:An}]
For each $g\geq 2$ and $r\geq 1$, $\cT_{g,r}$ is connected. So the
result follows from the previous proposition.  
\epv

\subsection{Unique realization for some dihedral angle assignations}

The content of the previous subsection is sufficient to ensure that the
map $\Phi_n$ is a covering for each choice of $n$, and that the number
of inverse images of all elements of $\cA$ is the same. The next lemma
states that this number is $1$ for some cellulations of $M$ and choices
of dihedral angles of the edges. 

\blm \label{lm:ex-unique}
Let $\sigma$ be a cellulation of $\dr M$ such that:
\begin{itemize}
\item the faces of $\sigma$ can be separated into two sets $F_-$ and
$f_+$, with each edge bounding a face of $F_-$ and one of $F_+$.
\item each vertex is adjacent to $4$ faces. 
\end{itemize}
Let $w:\sigma_1\rightarrow (0, \pi)$ be the function assigning the value
$\pi/2$ to each edge.
Then $w$ satisfies the hypothesis of theorem \ref{tm:dihedral}, and
there is at most one ideal hyperbolic manifolds with combinatorics given
by $\sigma$ and dihedral angles given by $w$. 
\elm

\bpv
Let $g$ be a ideal hyperbolic metric on $M$ with boundary combinatorics
$\sigma$ and all angles $\pi/2$. We consider $4$ copies of $M$, say
$M_1, M_2, M_3$ and $M_4$, each with the metric $g$. We identify the
corresponding faces of the $M_i$ according to the following table, where
a ``$+$'' on line $i$ and column $j$ means that each face of $M_i$ in
$F_+$ is glued to the corresponding face on $M_j$, and the same for $-$. 
\begin{center}
\begin{tabular}{|c|c|c|c|c|}
\hline
& 1 & 2 & 3&4 \\
\hline
1&&--&&+\\
\hline
2&--&&+&\\
\hline
3&&+&&--\\
\hline
4&+&&--&\\
\hline
\end{tabular}
\end{center}
Let $N$ be the result of those gluings. It is easy to check that:
\begin{enumerate}
\item the topology of $N$ is determined by the topology of $M$ and of
$\sigma$, and does not depend on the ideal hyperbolic metric on $M$.
\item $N$, with the metric coming from the gluing of the $M_i$, is a
complete, finite volume hyperbolic manifold.
\item two non-isometric ideal hyperbolic metrics on $M$ would result in
two non-isometric hyperbolic metrics on $N$.
\end{enumerate}
By Mostow rigidity there is at most
one complete, finite volume hyperbolic metric on $N$. The result
follows. 
\epv

\subsection{Proof of the main theorem}

We can now prove theorem \ref{tm:dihedral}, concerning the description
of dihedral angles of ideal hyperbolic manifolds. 

\btm \label{tm:dihedral}
Suppose that $M$ has incompressible boundary. Then the set $\Theta_\sigma$
of dihedral angles of ideal hyperbolic structures whose boundary is
cellulated according to $\sigma$ is given by the condition that:
\begin{enumerate}
\item on each elementary circuit in $\sigma$, the sum of the
exterior dihedral angles is $2\pi$;
\item on each non-elementary circuit, the sum of the dihedral angles 
is strictly larger than $2\pi$. 
\end{enumerate}
Each dihedral angles assignation is obtained on a unique ideal
hyperbolic manifold.
\etm

\bpv
Let $\sigma$ be a cellulation of $\dr M$. By lemma
\ref{lm:triangulation-bord}, there exists a number $\nu(\sigma)\in \N$
such that each dihedral assignation on the edges of $\sigma$, satisfying
the hypothesis of the theorem, is realized by exactly $\nu(\sigma)$
ideal hyperbolic manifolds with boundary combinatorics $\sigma$. 

Let $n\in (\N\setminus \{ 0\})^N$. By lemma \ref{lm:An}, $\cA_n$ is
connected; therefore, for all cellulations $\sigma$ of $\dr M$ with
$n_i$ vertices in $\dr_iM$, the value of $\nu(\sigma)$ is the same, and
is equal to a number $\nu(n)\in \N$.

Let $\sigma$ be a cellulation of $\dr M$. Let $\sigmab$ be the
cellulation of $\dr M$ obtained by ``splitting'' a vertex $v$ of $\sigma$,
i.e. replacing $v$ by two vertices $v_1$ and $v_2$ connected by an edge $e$,
and replacing two edges $e_1$ and $e_2$ containing $v$ by two edges
each, say $e_{11},e_{12}$ and $e_{21},e_{22}$. If $v\in \dr_iM$, let
$\nb=(n_1, n_2, \cdots, n_{i+1}, \cdots, n_N)$. 

Each dihedral
angle assignation $\theta$ to the edges of $\sigma$ can be identified
with a ``degenerate'' dihedral angle assignation $\thetab$ on $\sigmab$,
with a dihedral angle $0$ at $e$, dihedral angles at $e_{i1}$ and $e_{i2}$
equal to half the dihedral angles at $e_i$, and the other dihedral
angles equal to those of $\theta$. If $\theta$ satisfies the condition
of theorem \ref{tm:dihedral}, so does $\thetab$ (except of course for
the fact that one of its angles is $0$). Let $(\theta_n)$ be a
sequence of dihedral angle assignations on the edges of $\sigmab$ which
converges to $\thetab$. 

For each $i\in \{ 1, \cdots, \nu(\nb)\}$, there is a sequence of
ideal hyperbolic metrics $(g_{i,n})_{n\in \N}$ with dihedral angles
equal to $\theta_n$. As $n\rightarrow \infty$, $(g_{i,n})$ converges to
an ideal hyperbolic metric $g_i$ on $\dr M$ with boundary combinatorics
given by $\sigma$, and dihedral angles by $\theta$. So $\nu(\nb)\geq
\nu(n)\geq 1$. 

But it is easy to check that, for each $n\in (\N\setminus \{ 0\})^N$,
there exists some $p\in (\N\setminus \{ 0\})^N$, with $p_i\geq n_i$ for
all $i\in \{ 1, \cdots, N\}$, such that there exists a cellulation of
$\dr M$, with $p_i$ vertices in each $\dr_iM$, to which lemma
\ref{lm:ex-unique} applies. 

So, for each $n\in (\N\setminus \{ 0\})^N$, $\nu(n)=1$, which proves the
theorem. 
\epv

It might be useful to note that the convexity of the boundary of $M$ is not
really necessary for many steps of the proof of this theorem; for
instance the infinitesimal rigidity results should still hold for at
least some manifolds with non-convex boundary. This contrasts with most
results concerning the induced metric on the boundary, where the
convexity is crucial.

\subsection{An affine structure on some Teichm{\"u}ller spaces}

\begin{figure}[h]
\centerline{\psfig{figure=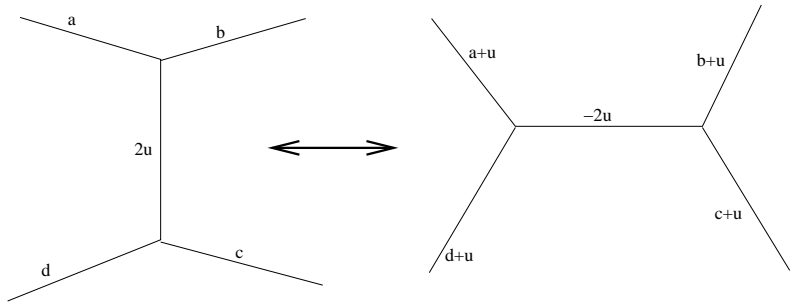,height=4cm}}
\caption{Definition of $\cA$, on the dual graph}\label{fig:transfo}
\end{figure}

A consequence of the previous considerations is that there is a natural
affine structure on the space of ideal hyperbolic structures on $M$,
coming from the parametrization by the dihedral angles.
Its definition
is simple and depends directly only on the dihedral angles. It has
``cone-like'' singularities along the ``cells'' of codimension at least
$2$, corresponding to non-generic polyhedral structures of the
boundary. The remarkable point, however, is that it also has a
well-defined volume element. This, which I can not explain better
than by doing an explicit computation, might indicate that this affine
structure is somehow meaningful. The definition below should be easier
to understand with figure \ref{fig:transfo}.

\bdf 
The affine structure $\cA$ on $\Hsm$ is defined as follows.
\begin{itemize}
\item for each triangulation $\sigma$ of the boundary, the restriction of
$\cA$ to $\Hsm^\sigma$ is defined as $\cA_{\sigma}$ above.
\item if $\sigma$ is a cellulation with one 4-gonal face $f$ and the other
faces triangles, then $\Hsm^{\sigma}$ is a codimension 1 face of $\Hsm$
which bounds two maximal dimension faces $\Hsm^{\sigma_1}$ and
$\Hsm^{\sigma_2}$, where $\sigma_1, \sigma_2$ are triangulations
obtained by refining $\sigma$ by adding edges $e_1$ and $e_2$
respectively to the 4-gonal face of $\sigma$. Then $\Hsm^{\sigma_1}$ and
$\Hsm^{\sigma_2}$ are glued by the map, from a neighborhood of $\sigma$
in the extension of $\sigma_1$ corresponding to having a negative angle
at $e_1$ to $\Hsm^{\sigma_2}$, which sends a configuration with angle
$-2u$ at $e_1$ and $a,b,c,d$ at the 4 edges of $f$, to a
configuration with angle $2u$ at $e_2$, and angles $a-u, b-u,c-u,d-u$ at
the edges of $f$.
\end{itemize}
\edf

Note that this definition makes sense since the transformation preserves
the angle conditions defining the space of possible angle assignations.
It is clear that the lengths conditions on elementary and non-elementary
circuits are satisfied for one of the configurations if and only if
they are satisfied for the other.

\blm
$\cA$ has a holonomy preserving a volume form.
\elm

\bpv
To prove this one should check that the holonomy around all codimension
two faces of $\Hsm$ preserve a volume form. There are two kinds of
codimension two faces, those corresponding to a cellulation of $\dr M$
with two 4-gonal faces, and those corresponding to the cellulations with
one 5-gonal face. In the first case of the holonomy is trivial, so the
only case to check is the second.

\begin{figure}[h]
\centerline{\psfig{figure=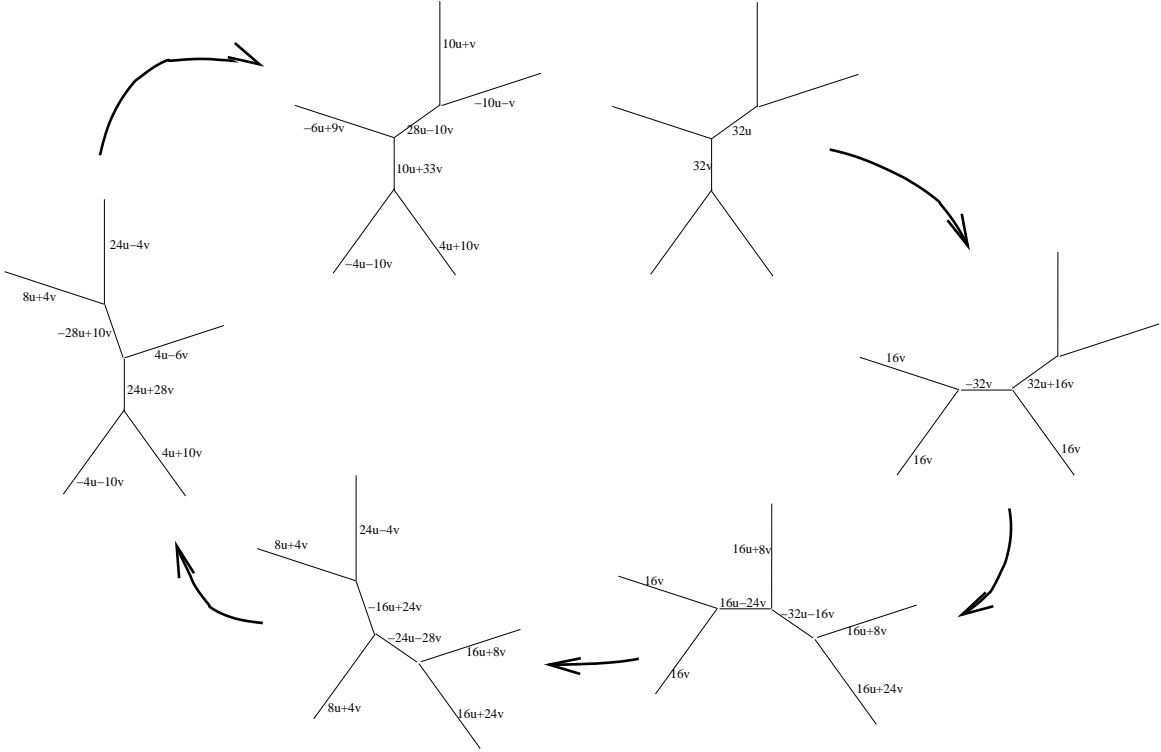,height=10cm}}
\caption{Holonomy computation}\label{fig:holonomie}
\end{figure}

When a cellulation $\sigma$ of $\dr M$ has one 5-gonal face and its
other faces are triangles, the corresponding cell $\Hsm^{\sigma}$ of
$\Hsm$ bounds five maximal dimension and five codimension one faces. The
holonomy around $\Hsm^\sigma$ is computed by considering the
diagram in figure \ref{fig:holonomie}, corresponding to the
transformation of the angles as one goes 
through the five maximal dimension faces --- the angles are marked on the
dual graph, and the constant part of the angles on the five ``exterior''
edges are not marked. 

The only part of the holonomy which is relevant for the volume is the
one corresponding to the 
two ``central'' edges. The computation made on figure
\ref{fig:holonomie} shows that
the holonomy for those two angles is given by the matrix:
$$ M:=\frac{1}{32}\left(
\begin{array}{cc}
28 & -10 \\
10 & 33 
\end{array}
\right) $$
But $\det(M)=1$, so that the holonomy around $\Hsm^\sigma$ preserves a
volume form.
\epv

\brk
$(\Hsm, \cA)$ has locally convex boundary at all points which are in
the closure of a face $\Hsm^{\sigma}$.
\erk

\bpv
We have to prove that if $\sigma$ is a triangulation of $\dr M$ and
$F_0$ is a codimension 1 face of $\Hsm^{\sigma}$ which is in the
boundary of $\Hsm$, then, across each codimension 1 face of $\Hsm^\sigma$
which is not in $\dr\Hsm$, the extension of $F_0$ is either outside
$\Hsm$ or in its boundary.

Let $\Hsm^{\sigma_1}$ be a codimension 1 face of $\Hsm^\sigma$, so that
$\sigma_1$ is obtained from $\sigma$ by removing an edge from $\sigma$,
and $\sigma_1$ has exactly one face $f$ with $4$ edges.
There are 3 cases to consider:
\begin{enumerate}
\item $F_0$ corresponds to a non-elementary circuit on which the sum of
the dihedral angles is $2\pi$. Then the definition of $\cA$ shows
directly that the extension of $F_0$ beyond $\Hsm^{\sigma_1}$ remains in
the boundary of $\Hsm$.
\item $F_0$ corresponds to an edge $e\not\in f$ with dihedral angle equal to
$\pi$. It is then clear again from the definition that the extension of $F_0$
beyond $\Hsm^{\sigma_1}$ remains in $\dr \Hsm$, since the transformation which
occurs in the definition of $\cA$ does not change the angle at $e$.
\item $F_0$ corresponds to an edge $e\in f$ with dihedral angle equal to
$\pi$. 
But then the 
definition of $\cA$ shows that, if one goes beyond $F_0$ on an affine
line starting in $\Hsm^{\sigma_1}$, the angle at $e$ is above $\pi$
after $\Hsm^{\sigma_1}$, so that the extension of $F_0$ beyond
$\Hsm^{\sigma_1}$ is outside $\Hsm$.
\end{enumerate}
\epv

An elementary consequence of the constructions which we have just seen
--- applied to the fuchsian case --- is that there is a natural affine
structure on the Teichm{\"u}ller space of a surface of genus at least 2,
with at least one marked point. The point is that fuchsian manifolds
form an affine submanifold of the ideal hyperbolic manifolds (they
corresponds to the case where the triangulations and dihedral angles are
the same on the two components of the boundary). Moreover, 
if one remains in the
category of fuchsian manifolds, there is no bent case (see remark
\ref{rk:fuchsian-bent}), so that the 
affine structure defined on ideal manifolds spans the whole Teichm{\"u}ller
space. One thus obtain the:

\btm \label{tm:teichmuller-affine}
For each $g\geq 2$ and each $N\geq 1$, there is a natural unimodular
piecewise affine structure $\cA_{g,N}$ on
the Teichm{\"u}ller $\cT_{g,N}$ space of the genus $g$ surface with $N$
marked points.
\etm

There are some natural
questions that remain on this structure; for instance, how many cells
it has, and whether the affine structure is globally convex.

\subsection{Hyperbolic polyhedra}

The previous results apply in a specially simple way to ideal polyhedra
in $H^3$, i.e. the case where $M$ is topologically a ball. Then, for a
fixed number $N\geq 3$ of vertices, there is a fixed number of
combinatorial types of possible Delaunay 
cellulations of $S^2$ with $N$ vertices, so that the space $\Hsm$ of
dihedral angles assignations is the union of a finite number of cells. 
But even here some questions remains, for instance whether the set of
ideal polyhedra with $N$ vertices, with the affine structure coming from
the dihedral angles, is globally convex.

\section{Induced metrics}

We will study in this section some properties of the 
metrics induced on the boundaries of ideal hyperbolic 
manifolds. They are finite area hyperbolic metrics on each connected
component of the boundary of $M$ without its ideal points, and we 
will show that the infinitesimal deformations 
of the interior metric are
parametrized by the infinitesimal deformations of the boundary 
metric; this is again a consequence of the Schl{\"a}fli formula. In
addition, we will prove a global result in the special case of fuchsian
manifolds, then the induced metrics on the boundary are in one-to-one
correspondence with the hyperbolic structures on the interior.

\subsection{Finite area metrics on surfaces}

We will use again here the notion of hyperbolic structure on a
triangulated surface, as defined in the previous section. We consider
again a compact, orientable surface $S$ of genus at least $2$ with $N$
points $v_1, \cdots, v_N$ removed. The following
additional definition is natural. 

\bdf \label{df-sh-som}
Let $g\in \cM$ and let $v$ be a vertex of $\sigma$. The {\bf shift} of
$g$ at $v$ is the sum of the shifts of $g$ at the edges containing
$v$. We will say that $g$ is {\bf complete} if its shift is zero at all
vertices. The set of complete structures will be denoted by $\cM_c$.
\edf

Of course the notion of completeness defined here is the same as the
usual, topological notion. Indeed if the shift of $g$ at a vertex $v$ is
non-zero, it is possible to use this --- and the fact that ideal
triangles are exponentially thin near their vertices --- to attain $v$ in
a finite time, by ``circling'' around it to take opportunity of the
shift. The reciprocal is not difficult to prove either.

Finally, knowing the shifts of a complete metric $g\in \cM_c$
at all edges determines all the gluing along the edges; let $S_0$ be the
set of maps from the set of edges of $\sigma$ to $\R$ such that, at each
vertex, the sum of the adjacent edges is $0$. We
identify $S_0$ with $\R^{e-v}$. Then:

\bprop \label{pr-shifts}
The map $F$ from $\cM_c$ to $S_0=\R^{e-v}$ sending $g$ to $F(g)$
defined, for an edge $E$, by $F(g)(E)=\sh_g(E)$, is a bijection.
\eprop

\bpv
It is quite obvious that the hyperbolic structures on a given
triangulated surface are determined by their shifts. On the other hand
each shift function on the edges corresponds to a possible gluing (under
the condition that 
the sum on the edges at any vertex is zero, by completeness). 
\epv

\subsection{The lengths of the edges of $\sigma$}
Let $g\in \cM_c$. Each vertex $v_i$ of 
$\sigma$ has a neighborhood which is
isometric to a neighborhood of a cusp in the quotient of $H^2$ by a
parabolic isometry. For each such $v_i$, choose a horocycle $H_i$
centered at $v_i$. Then, for each edge $E$ of $\sigma$, going from a
vertex $v_i$ to a vertex $v_j$, let $l(E)$ be the oriented length of $E$
between $H_i$ and $H_j$. The orientation is chosen so that $l(E)$ is
negative when the horoballs bounded by $H_i$ and $H_j$ overlap.
Clearly, replacing $H_i$ by another horocycle centered at $v_i$ changes
the lengths of the edges containing $v_i$ as an end by the addition of a
constant. So $l$ defines a function:
$$ l:\cM_c\rightarrow \R^e /\R^v~. $$

\bprop \label{pr-bij}
$l$ is a bijection between $\cM_c$ and $\R^e /\R^v$.
\eprop

The proof uses the following elementary property of ideal triangles in
$H^2$. 

\bsl
Let $x_1, x_2, x_3, x_4$ be four distinct points on $S^1=\dr_\infty
H^2$, in this cyclic order. For each $i\in \{1,2,3,4\}$, let $h_i$ be a
horocycle centered at $v_i$, and, for $i\neq j$, let $l_{ij}$ be the
distance between $h_i$ and $h_j$ along the geodesic going from $v_i$ to
$v_j$ --- which is negative if $h_i$ and $h_j$ overlap. Let $\pi_2$ and
$\pi_4$ be the orthogonal projections on $(x_3,x_1)$ of $x_2$ and $x_4$
respectively, and let $\delta$ be the oriented distance between $\pi_2$
and $\pi_4$ on $(x_3,x_1)$. Then:
$$ 2\delta = l_{12}-l_{23}+l_{34}-l_{41}~. $$
\esl 

\begin{figure}[h]
\centerline{\psfig{figure=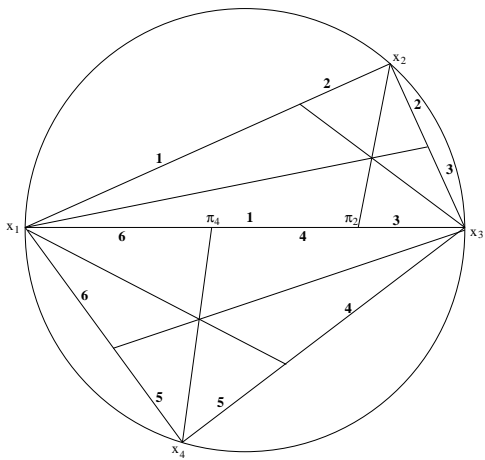,height=6cm}}
\caption{Ideal triangles (in the projective model of $H^2$)} 
\label{fig:triangles}
\end{figure}

\bpv
It follows from figure \ref{fig:triangles}, where numbers from 1 to 6 are
attached to lengths of segments. Elementary properties of the ideal
triangle show that:
\begin{eqnarray}
l_{12}-l_{23}+l_{34}-l_{41} & = & (1+2)-(2+3)+(4+5)-(5+6)\nonumber 
\\
& = & 1-3+4-6\nonumber \\
& = & (1-6) + (4-3) \nonumber \\
& = & 2\delta~. \nonumber
\end{eqnarray}
\epv

\bpn{of proposition \ref{pr-bij}}
Using proposition \ref{pr-shifts}, we are brought to proving that the
map: 
$$ L:S_0\simeq \R^{e-v}\rightarrow \R^e/\R^v $$
sending the $e$-uple of shifts of a metric to the lengths of its edges
is a bijection.
But the sub-lemma above shows that $L$ is linear, and explicitely
describes its inverse.
\epn

\subsection{Rigidity and the Schl{\"a}fli formula}

We now consider an \ihm $M$; thanks to the previous paragraphs, we know
that, if $\dr M$ is triangulated, then its induced metric is determined
by the lengths of the edges (defined up to the addition of a constant
for each vertex). Of course, if $\dr M$ is not triangulated --- i.e. if
some of the faces are polygons with more than three vertices --- one can
subdivide its cellulation to obtain a triangulation.

We will prove in this paragraph and the 
next that $M$ is infinitesimally
rigid, in the following sense. We call $\Phi$ the map from $\Hsm$ to
$\cM_c$ sending $g_0\in \Hsm$ to the metric induced on its
boundary. Then: 

\blm \label{lm-rigid}
For any $g_0\in \Hsm$, $T_{g_0}\Phi$ is an isomorphism.
\elm

The proof will take two steps. We first treat the case when the
boundary of $M$ is triangulated, then the general case. The first case
will rely on the Schl{\"a}fli formula and remark \ref{rk:concavity}. Recall
from lemma \ref{lm:triangulation-int} that, for any ideal hyperbolic
manifold  $(M, g_0)$, the volume $V$ is a strictly concave function on
$\Hsm$ in the neighborhood of $g_0$. 

\bprop \label{pr-cas1}
Suppose that all faces of $(M, g_0)$ are triangles, and that it is
triangulable. Then $T_{g_0}\Phi$ is an isomorphism.
\eprop

\bpv
We know by corollary \ref{cr:dim-Hsi} that $\dim \Hsm = e-v=\dim
\cM_c$. Moreover, the volume $V(g)$ is a strictly concave function on
$\Hsm$ by lemma \ref{lm:triangulation-int}. 


Recall the Schl{\"a}fli formula from section 1:
$$ dV = -\frac{1}{2} \sum_E l(E) d\theta(E)~, $$
where the sum is over all edges, and $\theta(E)$ is the interior
dihedral angle of edge $E$. The strict concavity of $V$ with respect to
the $\theta(E)$ means that the Hessian of $V$ with respects to those
angles is negative definite, so that the differential of the lengths
$l(E)$ --- with respect to the dihedral angles --- 
is non-degenerate. Thus
$T_gl$ is injective, and the result follows.
\epv

Now let $(M,g)$ be an ideal hyperbolic manifold with non-triangular
boundary. By lemma \ref{lm-atri}, $M$ has a finite covering $\Mb$ which admits
an ideal triangulation $\Sigma$. $\Sigma$ induces an ideal triangulation
of $\dr \Mb$, and the argument above --- the convexity of the volume and
the Schl{\"a}fli formula --- still implies that, for each assignation of a
real number $r(e)$ to each edge $e$ of $\dr \Mb$ (defined up to addition
of a constant for each vertex) there is a unique
infinitesimal deformation of $\Mb$ such that the length of $e$ (defined
with respect to a set of horospheres at the vertices)
varies at the rate $r(e)$. Therefore the argument given above shows that
the infinitesimal deformations of $\Mb$ are uniquely determined by the
induced variation of the boundary metric, and that each variation of the
boundary metric is obtained. 

This is in particular true for the variation of the boundary metric which
are invariant under the action of $\pi_1M/\pi_1\Mb$ on $\Mb$, and the
uniqueness implies that they are associated to infinitesimal
deformations of $\Mb$ equivariant under the same action. So we obtain
the proof of lemma \ref{lm-rigid}, which we can reformulate as:

\blm \label{lm:defos-ideales}
For each ideal hyperbolic manifold $M$, and for each infinitesimal
deformation $\hb$ of the boundary metric $h$ on $\dr M$, there is a unique
infinitesimal deformation of $M$ inducing $\hb$.
\elm

\subsection{The fuchsian case}

Lemma \ref{lm:defos-ideales} provides only an infinitesimal deformation
result for the induced metrics on the boundaries of ideal hyperbolic
manifolds. It is then difficult to obtain a global existence and
uniqueness results for the boundary case, in particular because we have
no such infinitesimal rigidity result for bent hyperbolic manifolds,
which necessarily enter the picture. 
But in the fuchsian case, the bent case is excluded by remark
\ref{rk:fuchsian-bent}, and a global result can be achieved.

\btm \label{tm:metrique-fuchsian}
Let $S$ be a surface of genus $g\geq 2$, and let $N\geq 1$. 
For each complete, finite area
hyperbolic metric $h$ on $S$ with $N$ cusps, there is a unique ideal
fuchsian hyperbolic manifold $M$ such that the induced metric on each
component of the boundary is $h$.
\etm

\bpv
We will use a deformation proof, following the original approach of
Aleksandrov \cite{Al} for similar polyhedral question. We choose $g\geq
2$ and $N\geq 1$. We call $\cT_{g,N}$ the Teichm{\"u}ller space of marked
conformal structures on the compact surface $\Sigma_g$ of genus $g$ with
$N$ marked points. There is a natural map $\Phi_{g,N}$ from $\cT_{g,N}$
to itself, defined as follows. For $h\in \cT_{g,N}$, the conformal
class on $\Sigma_g$ (i.e. forgetting the marked points on $\Sigma_g$)
contains a unique hyperbolic metric $h_0$ on 
$\Sigma_g$; taking its universal cover defines a conformal map from
$\Sigmat_g$ with the conformal structure in $h$ to the upper hemisphere
$S^2_+$. The marked points in $h$ then define an equivariant set $S$ of
points in $S^2_+$. Taking the boundary of the convex hull of those
points in $\R^3$ (and subtracting its intersection with the plane
containing the equator) leads to a polyhedral surface, which is
invariant under the natural action of $\pi_1\Sigma_g$ on
$H^3$. The quotient is homeomorphic to $\Sigma_g$, and has $N$
vertices (which are the intersection points with $S^2_+$); it carries a
complete, finite area hyperbolic metric $\hb$ induced by the canonical
metric on $H^3$. 
$\Phi_{g,N}(h)$ is defined as this metric $\hb$, considered as an
element of $\cT_{g,N}$. 

The proof of the theorem is an immediate consequence of the following
points. 
\begin{enumerate}
\item $\Phi_{g,N}$ is locally injective, i.e. its differential is an
isomorphism at each point. 
\item $\Phi_{g,N}$ is proper.
\item $\cT_{g,N}$ is connected and simply connected.
\end{enumerate}
Point (1) is just lemma \ref{lm-rigid}, while point (3) is well known
(since we use the marked Teichm{\"u}ller space). So we only have to prove
point (2). It can be reformulated as a degeneration statement: if
$(M_n)_{n\in \N}$ is a sequence of fuchsian ideal manifolds of genus $n$
with $N$ vertices, and if the sequence of induced metrics $(h_n)$
converges to a limit $h\in \cT_{g,N}$, then $(M_n)$ converges. 

In other terms, we have to prove that if either the conformal structure
on $\dr E(M_n)$ degenerates or the marked points collapse, then the
length of some closed geodesic on $\dr M$ goes to $0$.

In the case of a degeneration of the conformal structure on $\dr
E(M_n)$, the action of $\pi_1\dr M$ on $S^2_+$ also degenerates, and it
is then easy to check that some closed geodesic goes to $0$. Similarly,
if the set of marked points collapses, then the collapsing subset is
bounded, in $\dr M$, by a long thin tube of very small diameter, and
the compactness result follows. 
\epv


\subsection{Ideal polyhedra in $H^3$}

An elementary remark is that, for ideal polyhedra in $H^3$, the approach
above gives a simple proof of the existence and uniqueness of ideal polyhedra
having a given induced metric, see \cite{rivin-comp,rivin-annals}. In
particular this approach does not use the Cauchy method to get the
global uniqueness, since it follows from the global deformation result. 

Rather, the infinitesimal rigidity can be obtained as in lemma
\ref{lm:defos-ideales} using the Schl{\"a}fli formula and the results
concerning the dihedral angles. Applying a deformation argument to
conclude then only requires a compactness result --- namely, that if a
sequence of ideal polyhedra has induced metrics which converge, then it
has a converging subsequence (modulo isometries). Such a compactness
result is easy to obtain directly.

\section{Circles packings}

\subsection{From dihedral angles to circle configurations}

As mentioned in section 2, the polyhedral questions considered here can
be formulated in terms of configurations of circles in $S^2$. More
precisely, let 
$M$ be an ideal hyperbolic manifold. Consider its universal cover $\Mt$
as a subset of $H^3$. Its boundary $\dr \Mt$ is the disjoint union of a
set of convex surfaces, and their boundary is the limit set
$\Lambda$ of the action of $\pi_1M$ on $H^3$.

Each face of $\Mb$ defines an oriented totally geodesic plane in $H^3$, with
boundary at infinity an oriented circle in $S^2\setminus \Lambda$. If two faces
$F, F'$ of $\dr \Mt$ are adjacent, then the corresponding oriented circles $C,
C'$ intersect with angle equal to the exterior dihedral angle between
$F$ and $F'$. Moreover, a simple convexity argument shows that the union
of the interiors of the 
circles corresponding to the faces of $\dr \Mt$, along with the
intersection points, is $S^2\setminus
\Lambda$. So the results concerning the dihedral angles of ideal
hyperbolic manifolds can be formulated in terms of configurations of
circles in 
the complement of $\Lambda$, with given intersection angles. 

Therefore, an ideal hyperbolic manifolds determines a circle
configuration on $\dr M$, in the following sense. 

\bdf 
Consider $\dr M$ endowed with a $\CP^1$-structure, for instance coming
from a complete hyperbolic metric on $M$. A {\bf circle configuration}
on $\dr M$ is a finite set of open disks in $\dr M$ (for the $\CP
1$-structure), such that:
\begin{itemize}
\item the union of the closures of the disks is $\dr M$. 
\item the intersection points of the circles bounding two disks are never
  contained in any of the disks.
\end{itemize}
\edf

Note that this definition implies that any point which is contained in
two of the circles is in at least a third. 

This relationship between ideal hyperbolic manifolds and circle
configurations means that the results concerning the dihedral angles
should be compared to 
those concerning configurations of circles, e.g. the rigidity
results of \cite{he-annals}. 

\medskip

One can also consider circle packings in the more restrictive sense
of sets of disks whose interiors are pairwise disjoint, with the
complement made of disjoint, polygonal regions. 
Translating theorem \ref{tm:dihedral}, we will find the following
result, which reduces to the classical Koebe circle packing theorem as
extended by Thurston (see \cite{koebe,thurston-notes}).

\btm \label{tm:koebe}
Let $\Gamma$ be the 1-skeleton of a triangulation of $\dr M$. There is a
unique couple $(g, c)$, where $g$ is a complete, convex co-compact
hyperbolic metric on $M$, and $c$ is a circle packing on $\dr M$ (for
the $CP^1$-structure defined on $\dr M$ by $g$) whose incidence graph
is $\Gamma$. 
\etm

=== V{\'e}rifier ref Thurston et al. 

It is helpful to use a trick due to Thurston
\cite{thurston-notes}. Consider such a circle packing, such that the
complement of the circles is the disjoint union of
"triangles" bounded by 3 circle arcs. Then add, for each connected
component of the complement of the circles, a circle which is orthogonal
to the 3 circles which it interests. One then obtains a configuration of
circles intersecting at right angles.
The same construction works in the
other way. Consider a configuration of circles intersecting at right
angles, such that:
\begin{itemize}
\item the circles can be separated in two sets, the "white" and the
  "black" circles, such that each "black" circle only intersects
  transversally "white" circles, and conversely.
\item each "black" circle intersects transversally exactly three "white"
  circles.  
\end{itemize}
Then the set of "white" circles make up a circle packing.

Therefore, the proof theorem \ref{tm:koebe} is reduced to finding circle
configurations on $\dr M$, where the intersections between the circles
are always at right angles. Moreover, for the circle configuration
coming from circle packings by adding orthogonal circles, it is easy to
check that all circuits are made of at least 4 edges, and strictly more
than 4 unless they are elementary. Therefore, theorem \ref{tm:dihedral}
applies, and theorem \ref{tm:koebe} follows. 

If one considers graphs which are the 1-skeleton of cellulations (more
general than triangulation) the same construction works; one obtains
circles packings with the added property that, for each intersticial
region, there is a circle orthogonal to all the adjacent circles. The
uniqueness, under this additional condition, is proved by noting that
each such packing gives rise to a circle configuration with right angle
intersections, and thus to an ideal hyperbolic manifold with boundary
faces intersecting at right angle (whose uniqueness is known by theorem
\ref{tm:dihedral}). 

\medskip

Note that, seen in this light, the
statements made here are related to different results on the
rigidity of circle packings, see e.g. \cite{schramm}.

\medskip

In the fuchsian case, the results above can be stated as describing
circle packings --- or configurations of circles with given angles ---
on surfaces of genus $g\geq 2$.

\subsection{Circles packings and $CP^1$-structures}

A slightly different idea is to consider the same circle packings on the
boundary at infinity of $E(M)$, which, as a Riemann surface, can be
endowed with a hyperbolic metric. The quasi-fuchsian case of our statements
concerning dihedral angles should then be compared to recent results of
S. Kojima, S. Mizushima and S. P. Tan \cite{KMT}.

They consider a compact surface $S$ of genus $g\geq 2$, with a triangulation
$\tau$. By a result of Thurston \cite{thurston-notes}, there is then a
unique hyperbolic metric $g_0$ 
on $S$ and a unique circle packing $C_0$ in $g_0$ whose nerve is $\tau$; it is
called the Andreev-Thurston packing. Note that it can be recovered by
applying the "Thurston trick" described above and theorem
\ref{tm:dihedral} in the fuchsian case.

Kojima, Mizushima and Tan are interested in the deformations of the $\C P^1$
structure on $S$ defined by $g_0$, and of the circle packing $C_0$. They
prove among other things that, near $(g_0, C_0)$, there is a
$6g-6$-parameter deformation space. 

Notice that this local deformation result --- which is only the easier
part of the results of Kojima, Mizushima and Tan --- can be recovered from the
results of section 8, applied to the quasi-fuchsian case. Indeed, the
Thurston-Andreev circle packing corresponds to a fuchsian ideal manifold
$M$ with all dihedral angles equal to $\pi/2$. Let $\Sigma_+$ and
$\Sigma_-$ be the two connected components of $\dr M$. Consider the
deformations of the dihedral angles assignations on the edges of $\dr M$
which do not change the angles on $\Sigma_+$. This deformation space has
dimension equal to the number of edges $e$ of the triangulation $\tau$,
minus the number $v$ of its vertices --- because the sum at each vertex
of the exterior dihedral angles has to remain $2\pi$. But the Euler
formula shows that $e-v=6g-6$.

Now each such deformation of the dihedral angles determines a unique
quasi-fuchsian deformation of the hyperbolic structure on $M$ by lemma
\ref{lm:defos-possibles}; this in turn determines a unique deformation
of the $\C P^1$ structures on $\Sigma_+$ and $\Sigma_-$. Since the
dihedral angles remain equal to $\pi/2$ on $\Sigma_+$, the converse of
the Thurston construction above still determines a circle packing, which
is a deformation of $C_0$. This recovers the deformation studied by
Kojima, Mizushima and Tan.

\section{Concluding remarks and questions}

\subsection{Manifolds with cusps}

Most of what was done in this paper could presumably be extended from convex
co-compact manifolds to manifolds with cusps; that is, we consider
manifolds $M$ which remain of finite volume,
with a boundary which is polyhedral, i.e. the union of a finite set of
totally geodesic faces which intersect along edges, with all their
vertices ideal; but we allow those manifolds to have some cusps,
i.e. $M$ has finite volume but 
contains points arbitrarily far from the boundary.

The most delicate technical point to check is the content of section 4,
concerning ideal triangulation. Proposition \ref{pr-cell} extends to the
case with cusps using the ideas of Epstein and Penner
\cite{epstein-penner} in a manner which is closer to what they do in
their paper, i.e. using the properties of the action of $\pi_1M$ on the
parabolic points. The combinatorial arguments used in section 4 to
obtain an ideal triangulation from a cellulation then work without any
difference in the case with cusps.

Doing this in the special case of a manifold with one cusp, one could
presumably recover
some results of Thurston on circle packings on the torus (see
\cite{thurston-notes}, section 13.7). The dihedral angles, however, should
not be restricted to be acute in this approach.

\subsection{Hyperideal polyhedra}

A remarkable point is that some of the properties described in the
introduction survive when one considers polyhedra or surfaces with
complete metrics of 
infinite area. In the polyhedral case this corresponds to hyperideal
polyhedra; in the projective model of $H^3$, they can be described as 
polyhedra having some "usual" vertices in $H^3$ and some
"hyperideal" vertices beyond infinity, but such that all edges meet the
interior of $H^3$. An existence and uniqueness statement for the induced
metric and the third fundamental form on those polyhedra can be found in
\cite{shu}, while a recent work of Bao and Bonahon \cite{bao-bonahon}
describes the set of their dihedral angles --- a problem which is related
to the third fundamental form.

For smooth surfaces, it looks like a possible analog of hyperideal
polyhedra is the class of convex surfaces in $H^3$ whose boundary at
infinity is 
a circle. A statement concerning the induced metrics and the third
fundamental forms of such surfaces can be found in \cite{rsc}, but it
deals only with the special case of surfaces of constant Gauss
curvature. 
I also believe that one could prove a similar existence result for
hyperbolic metrics on a manifold $M$ inducing a given metric, with
constant curvature $K\in [-1, 0)$, on the boundary (see \cite{bkc}). 
A more general result, extending to metrics of non-constant
curvature on the boundary, would be more demanding in terms of
analytical techniques. 

It would be interesting to understand whether those property also remain
valid for hyperbolic manifolds with boundaries, when the metric has a
"boundary at finite distance" and also a "boundary at infinity", thus
extending the notion of hyperideal polyhedron. One
might consider either a "polyhedral" type "finite boundary" --- which
would have to be "hyperideal" in the sense that, for each end, all
edges "converge" to a hyperideal vertex --- or a smooth "finite
boundary", maybe with the additional condition that it meets the
boundary at infinity along a circle.

One basic problem here is that we still lack an
infinitesimal rigidity result here, stating that each infinitesimal
deformation of those objects induces a non-trivial variation of both the
induced metric and the third fundamental form. I believe that 
this could be achieved by using the same line of reasoning as here, but
in the hyperideal instead of the ideal case. More precisely, the
argument given at the end of section 2 to prove that the volume of ideal
simplices is concave, also works for hyperideal simplices, if the volume
of a hyperideal simplex is defined as the volume of the corresponding
truncated simplex. Once this key point is obtained, the same method as
here could presumably be used to understand ``manifolds with
hyperideal-like boundary''. 
A rigidity result for hyperideal polyhedra, or for manifolds
with boundaries locally like hyperideal polyhedra, would follow. I'm
hoping to come back to this matter in the future. 

\subsection{Complete manifolds of finite volume}

Some of the tools used here could also have applications to the
studying the deformations of complete hyperbolic manifolds of finite
volume. 

Those manifold have a decomposition in cells, each isometric to an ideal
polyhedron, by \cite{epstein-penner}. Thus the methods of section 5
might lead to ideal triangulations of finite covers, as in
e.g. \cite{neumann-zagier}. 

One could then study the deformations obtained by changing the dihedral
angles of the ideal triangles; the volume remains a concave functional
on the space of angle assignations. 

An important difference with the case of ideal hyperbolic manifolds,
however, is that the shears along the singular edges is not the only
obstruction to having a complete structure. The link of the
vertices is now a torus --- instead of a disk for ideal hyperbolic
manifolds --- which carries a similarity structure, with singularities
at points corresponding to the singular edges; to have a complete
structure, it is also necessary that the conformal part of the holonomy
of each link vanishes. But this vanishing could also be a necessary
and sufficient condition for the vanishing of the differential of the
volume among deformations which do not change the total angle around
each singular edge (at least it is sufficient).

\subsection{Convex cores}

The questions mentioned in the introduction, concerning the extension of
theorems \ref{tm:hmcb-I} and \ref{tm:hmcb-III} to cases where $\dr M$ is
not smooth, have another natural setting: the case where $M$ is supposed
to be the convex core of a complete, convex co-compact manifolds
$N$. This is equivalent to supposing that $\dr M$ is convex and
developpable, or that it is convex and has a hyperbolic (i.e. constant
curvature $-1$) induced metric. 

The question concerning theorem \ref{tm:hmcb-I} is whether, for any
hyperbolic metric $h$ on $\dr M$, there is a unique hyperbolic metric
$g$ on $M$ such that $\dr M$ is convex and developpable, with induced
metric $h$. The existence part holds, as was proved by Labourie
\cite{L6} (another proof exists, using a degree argument and the known
results from \cite{epstein-marden} on a conjecture of Sullivan). The
uniqueness, however, remains unknown. 

When one considers the third fundamental forms of the boundary, one
encounters a statement on the pleating measures of convex cores. Here
again, satisfactory results on the existence part have been obtained
recently by Bonahon, Otal \cite{bonahon-otal} and Lecuire
\cite{lecuire}, but the uniqueness is unknown. 

Note that the questions on the convex cores can rather naturally be
extended to the "horned hyperbolic manifolds" introduced in section
3. For instance, the question on the induced metric is: 

\bq
Let $h$ be a complete, finite area metric on $\dr M$ minus a finite
number of points. Is there a unique metric $g$ on $M$ such that
$(M, g)$ is a hyperbolic manifold with horns and that the induced metric on
$\dr M$ is $h$ ?
\eq

\subsection{Orbifolds, cone-manifolds}

The results given here, in particular theorem \ref{tm:dihedral}, can
might also be used to construct hyperbolic orbifolds, or hyperbolic
manifolds. This already happens in the fuchsian case, see
\cite{thurston-notes}. 

Other relation with cone-manifolds can be obtained, given an ideal
hyperbolic manifold $M$, by gluing two copies of it along their common
boundary. One obtains a finite volume, non-compact cone-manifold which
is singular 
along a family of geodesics --- corresponding to the edges of $M$ ---
each going between two points at infinity. Moreover, the convexity of
$M$ means that the total angle around the singular geodesics is always
strictly less than $2\pi$. The infinitesimal
rigidity of ideal hyperbolic manifolds could therefore be a consequence
of an infinitesimal rigidity result for such cone-manifolds, which could
perhaps be proved using the ideas of \cite{HK}. Moreover, those
non-compact cone-manifolds could themselves be of interest.

\subsection{Affine structures}

There are several questions related to the affine structures coming from
the dihedral angles of the boundary. It
would be interesting to understand if the affine structure extends also
to bent hyperbolic manifolds, so as to define an affine structure on the
space of hyperbolic manifolds with horns.

Another question is whether the space of ideal hyperbolic manifolds ---
or even of hyperbolic manifolds with horns --- is convex for this affine
structure.
A related question concerns the behavior near the boundary of the
cellulation by the combinatorics of the boundary. For instance, is this
cellulation locally finite at a boundary point corresponding to no
bending and no degeneration of the conformal structure ?

\subsection{Smooth boundaries}

The questions treated in this paper have analogs for hyperbolic
manifolds with smooth boundary. For instance, given a complete metric
$h$ of finite area, with $K>-1$ on $\dr M$
minus a given number of points, which is asymptotically hyperbolic, is
there a unique hyperbolic metric on $M$ with convex boundary, inducing
$h$ on the boundary ? 

The results quoted above on hyperideal manifolds might induce one to
believe that the same question can be asked for metrics of infinite are,
under the condition that the boundary at infinity of each end is a
circle in $\dr_\infty M$. This is the analog of polyhedral condition
that appeared in the previous sub-section. Even when $M$ is a ball,
however, this is proved only for constant curvature metrics on $\dr M$
minus $N$ points (see \cite{rsc}, which also contains a dual result on
the third fundamental forms).

\subsection{Higher dimensions}

There is some hope of finding some properties related to
questions A and B in higher dimension, in the setting of Einstein
manifolds of negative curvature. An elementary step in this direction is
taken in \cite{ecb}, where an infinitesimal deformation result
concerning question A is proved. 

Another similar perspective is given by a the work of Graham and Lee
\cite{graham-lee}, which can be considered as a first step towards a
possible extension of the classical Ahlfors-Bers theorem 
in higher dimension, again in the setting of Einstein manifolds. Similar
ideas form part of the mathematical side of the
"AdS-CFT correspondence", an important conjecture in mathematical
physics relating deep properties of conformally compact Einstein
manifolds to other deep properties of their conformal boundaries; see
e.g. \cite{maldacena,witten,witten-yau,graham-witten} and the references
given there.

\section*{Acknowledgments}

This works own much to remarks made, at various times, by Francis
Bonahon, Fran{\c c}ois
Labourie, Greg McShane and Igor Rivin. I also got a significant help
from Richard Kenyon. Other helpful comments came from Michel Boileau,
Philippe Eyssidieux, Gilbert Levitt, Anne Parreau and Joan Porti. I
would like to thank them all here.

\bibliographystyle{alpha}

\end{document}